\documentclass[journal,draftcls,onecolumn,letterpaper,twoside,11pt]{IEEEtran}

\usepackage[margin=1in]{geometry}
\usepackage{amssymb, mathrsfs,verbatim} 
\usepackage[cmex10]{amsmath}
\usepackage{cite}
\RequirePackage[colorlinks,linkcolor=black,citecolor=black,urlcolor=blue]{hyperref}
\usepackage{enumerate}
\usepackage[dvips]{graphicx}
\usepackage{verbatim}
\usepackage{bbm}
\usepackage{subfig}
\usepackage{epsfig}
\usepackage{epstopdf}
\usepackage{graphicx} % standard LaTeX graphics tool
\usepackage{subcaption}
\usepackage{subfig}

\usepackage{float,latexsym,amsfonts,amstext,times}

\usepackage{chngcntr}
\usepackage{apptools}

\usepackage[lined,linesnumbered,algoruled,noend]{algorithm2e}
\usepackage{cases}

\usepackage{color}
\newtheorem{theorem}{Theorem}[section]

\newtheorem{proposition}{Proposition}[section]
\newtheorem{corollary}{Corollary}[section]
\newtheorem{lemma}{Lemma}[section]

\newcommand{\cR}{\mathcal{R}}
\newcommand{\cS}{\mathcal{S}}

\newcommand{\cC}{\mathcal{C}}
\newcommand{\cF}{\mathcal{F}}

\newcommand{\cL}{\mathcal{L}}
\newcommand{\cN}{\mathcal{N}}

\newcommand{\cH}{\mathcal{H}}
\newcommand{\cE}{\mathcal{E}}

\newcommand{\Pro}{\mathsf{P}}
\newcommand{\Exp}{\mathsf{E}}
\newcommand{\dD}{\mathsf{D}}

\newcommand{\bN}{\mathbb{N}}
\newcommand{\bR}{\mathbb{R}}
\newcommand{\bS}{\mathbb{S}}
\newcommand{\dd}{\mathsf{d}}
\newcommand{\sL}{\mathsf{L}}

\begin{document}

\title{Sequential anomaly identification with sampling constraints for composite hypotheses
 %\thanks{Author  is  supported }
}

\author{Aristomenis Tsopelakos, Georgios Fellouris, Olgica Milenkovic
\thanks{}
}

%\markboth{IEEE Transactions on Information Theory,~Vol.~~, No.~~, ~~~Year (submitted)}%
%{Author: Title}

\maketitle

\begin{abstract}
We consider the problem of sequential anomaly identification with sampling constraints, where multiple independent data sources are observed over time, but only a limited subset can be sampled at each instant. The objective is to identify the anomalous sources while controlling the probabilities of at least one false alarm and at least one missed detection below prescribed levels. An identification policy must jointly determine when to stop sampling, which sources to declare anomalous, and how to allocate observations across sources under the sampling constraints. We seek policies that minimize the expected time to decision to a first-order asymptotic approximation as the error probabilities vanish. The observations from each source follow an unknown density belonging to one of two disjoint sets, corresponding to regular and anomalous behavior. Focusing on a nonparametric framework, which allows greater flexibility in modeling regular and anomalous statistical behavior, we establish general conditions for asymptotic optimality and provide a unified characterization of asymptotically optimal policies for composite hypothesis testing under sampling constraints. Simulation results illustrate the performance of the proposed policies, and the robustness of the nonparametric framework.
\end{abstract}	


\section{Introduction}\label{sec:Intro}
In many engineering applications, the goal is to classify a collection of data sources according to their underlying distributions based on sequentially collected measurements. For example, in navigation monitoring systems, one seeks to identify malfunctioning sensors that introduce noise into the measurements \cite{nikif93}. In computer network intrusion detection, deviations from a standard traffic profile, such as a reduced data transmission rate, may indicate a potential attack \cite{tart2}. In cognitive radio spectrum sensing \cite{cognitive}, the objective is to detect available frequency bands where the signal is of very low amplitude. Such applications motivate the study of sequential multiple testing problems for composite hypotheses \cite{wald2004sequential,wong1968asymptotically,schwarz1962asymptotic,lorden1977nearly,lai1988nearly,robbins1974expected,pavlov91,gafni2023anomaly,hemo_g}.

Composite hypothesis testing extends the classical binary hypothesis
testing framework by replacing two competing distributions with two
disjoint sets of distributions. This generalization introduces
additional complexity because the true distribution is known only
to belong to one of these sets. The study of sequential composite
hypothesis testing dates back to Wald \cite{wald2004sequential},
who considered a simple null hypothesis against a composite
alternative and introduced the Weighted SPRT (WSPRT) and the
Generalized SPRT (GSPRT). While the error probabilities of the
WSPRT can be controlled using techniques similar to those developed
for simple hypotheses, analyzing the GSPRT is more involved.
Lorden \cite{lorden1973open} studied open-ended generalized
likelihood ratio tests of a simple null hypothesis against a
one-sided composite alternative in exponential families,
establishing near-optimal expected sample sizes subject to a
bound on the probability of rejecting the null hypothesis.
Further asymptotic analyses of sequential tests include the
works of Lorden \cite{lorden1977nearly}, Woodroofe
\cite{woodroofe1982nonlinear}, and Chan and Lai \cite{chan_lai00},
the latter focusing on error-probability approximations in
multiparameter exponential families. For separate parametric
families that need not belong to a common exponential family,
Li et al. \cite{li2014generalized} established logarithmic
asymptotic characterizations of the GSPRT's maximal type I and
type II error probabilities, and derived first-order asymptotic
expressions for its expected sample size.

The Bayesian formulation of composite hypothesis testing was
introduced by Schwarz \cite{schwarz1962asymptotic}, who proposed
a first-order asymptotically Bayes-optimal procedure for parametric
families in which the parameter belongs to one of two disjoint
sets. Variants of this procedure were subsequently shown to
preserve asymptotic optimality
\cite{pavlov91,wong1968asymptotically,kiefer1963asymptotically},
and further refinements achieved second-order optimality
\cite{lorden1977nearly}. Related procedures were also developed
in the frequentist setting, where first-order asymptotic
optimality was established under appropriate conditions
\cite{lai1988nearly}. In a related direction, Yilmaz et al. \cite{yilmaz2015sequential} studied sequential joint
detection and estimation, where a decision in favor of the
alternative is accompanied by an estimate of a random parameter.
They showed that optimal joint performance requires accounting
for estimation in the design of the detection rule, rather
than optimizing the two tasks separately.

Optimal sequential testing was further extended to experimental design, where the decision maker determines both when to stop and which observations to collect. Chernoff
\cite{chernoff1959sequential} considered a simple null hypothesis against a finite set of alternatives and proposed an asymptotically optimal procedure based on adaptive experiment selection. This framework was subsequently extended to countably infinite alternatives \cite{Alb61} and multiple competing hypotheses \cite{bessler60,bessler60_2}, with later developments addressing frequentist error criteria \cite{nitinawarat2013controlled}.
More recently, Deshmukh et al.~\cite{deshmukh} studied controlled sequential testing of composite hypotheses in exponential family models, with the unknown parameter belonging to one of several prespecified sets. Assuming that the number of available experiments equals the number of hypotheses, they established first order asymptotic optimality using sampling rules inspired by the multi-armed bandit literature \cite{garivier2016optimal}.

In our work, we consider multiple data sources, each associated with a composite hypothesis testing problem, and adopt a frequentist formulation in which the objective is to minimize the expected stopping time as the error probabilities approach zero. The observations from each source follow a density that belongs to one of two disjoint sets, corresponding to regular and anomalous behavior, respectively. A source is declared anomalous if its underlying density belongs to the latter set. Within exponential family models, the problem of sequential anomaly identification for composite hypotheses has been studied in \cite{song_fel_gener,song_fel_supple}, where all sources are observed at each time and generalized error metrics are considered. The authors propose first-order asymptotically optimal policies and provide explicit characterizations of the optimal performance.

We now consider the nonparametric setting under sampling constraints, where only a subset of sources can be observed at each time due to resource limitations. Nonparametric formulations of sequential hypothesis testing impose minimal assumptions on the underlying distributions and have been studied both in classical settings \cite{darling1968some} and, more recently, through time-uniform, nonasymptotic inference methods \cite{howard2021time}. These works highlight the inherent difficulty of the problem when the statistical model is not fully specified. In contrast to the existing literature, we consider a setting in which multiple sources must be monitored under explicit sampling constraints. In particular, at each time instant only a subset of sources can be observed, and the goal is to design a policy that jointly determines when to stop, which sources to declare anomalous, and how to allocate observations across sources. The stopping and decision rules must ensure that the error probabilities remain below prescribed levels, while the sampling rule must satisfy the imposed constraint. Existing works on sequential anomaly identification typically assume that all sources are observed at each time or focus on restricted settings. For example, when there is a single anomalous source and only one observation can be collected at each time, asymptotically optimal policies have been derived in \cite{oddball,hemo_g}. In our formulation, we allow for multiple anomalous sources, simultaneous observations of multiple sources, and more general sampling constraints, expressed in terms of the ratio between the expected number of observations and the expected stopping time. Examples show that accounting for distributional shape beyond parameters is essential for efficient sampling allocation, further motivating the necessity of a nonparametric framework.

We impose regularity conditions on the nonparametric family
similar in spirit to those introduced by Pavlov
\cite[Section~3]{pavlov91}, who established their validity
for multidimensional exponential families under compactness.
Within this framework, we establish conditions for consistency
and asymptotic optimality of the proposed procedures. In particular,
we show that suitable convergence of the test statistics and
the estimators ensures these properties. In addition, asymptotic optimality is achieved when the long-run sampling frequency of each source exceeds an explicitly computable threshold. Based on this characterization, we design probabilistic sampling rules that attain the optimal performance. Our results generalize our previous work on sequential anomaly identification for the binary (simple vs. simple) case \cite{a_prob21}, and establish connections with the multi-armed bandit literature \cite{lai1985asymptotically,anantharam1987asymptotically,garivier2016optimal,prab22,bubeck2009pure}. In classical bandit problems, the objective is to maximize cumulative reward, which requires identifying the arm with the largest mean as quickly as possible \cite{lai1985asymptotically,anantharam1987asymptotically,garivier2016optimal}. In contrast, in \cite{prab22}, each arm is associated with a composite hypothesis testing problem, and the objective is to minimize the expected stopping time and switching costs. This setting naturally aligns with pure-exploration bandit problems \cite{bubeck2009pure,garivier2016optimal}, where the goal is to identify a subset of arms satisfying a statistical property. Our results provide a unified characterization of asymptotically optimal sampling strategies across these settings. The proposed framework is also closely connected to recent developments in machine learning that integrate generative models with active learning for anomaly detection \cite{liu2019generative,schlegl2019f,ruff2021unifying}, highlighting the role of adaptive data acquisition in improving detection performance. However, in these approaches, data acquisition is typically passive, with all features or data streams observed uniformly over time. In contrast, our formulation explicitly incorporates sampling constraints and enables the sequential allocation of observations to the most informative sources, providing a principled approach to active anomaly detection under resource constraints.

The remainder of the paper is organized as follows.
Section~\ref{sec:formulate} introduces the problem formulation,
Section~\ref{motiv} provides an illustrative example, and
Section~\ref{UALB} derives a universal asymptotic lower bound
(UALB) on the expected stopping time. Section~\ref{prob_rules}
presents a family of policies, and Section~\ref{asyo} establishes
conditions under which an appropriately designed policy attains
the UALB and is therefore asymptotically optimal.
Sections~\ref{test_st} and~\ref{estim} present test statistics
and estimators that satisfy the required asymptotic conditions.
Sections~\ref{simul1} and~\ref{simul2} evaluate the performance
of the proposed policy through simulation studies. Finally,
Section~\ref{concl} concludes the paper and outlines directions
for future research.

We end this section with some notations we use throughout the paper. We use $:= $ to indicate the definition of a new quantity and $\equiv$ to indicate a duplication of notation. We set $\bN := \{1, 2 \ldots, \}$ and $[n] := \{1, \ldots, n\}$ for  $n \in \bN$, we denote by $|A|$ the size and by $2^A$ the powerset of a set $A$,  by $\lfloor  a  \rfloor$ the floor and by $\lceil a  \rceil$ the ceiling of a  positive number $a$, and by $\mathbf{1}$ the indicator of an event. The symbol $\vee$ stands for the maximum, and the symbol $\wedge$ for the minimum of two quantities. We write $x\sim y$ when  $\lim (x/y)=1$, $x\gtrsim y$ when  $\liminf (x/y) \geq 1$, and $x\lesssim y$ when  $\limsup (x/y) \leq 1$, where the limit is taken in some sense that will be specified. We denote by $\|\cdot\|_p$ the $\sL^p$-norm of a function. Moreover, \textit{iid} stands for  independent and identically distributed, and we say that a sequence of positive numbers $(a_n)$ is \textit{summable} if $\sum_{n=1}^\infty a_n< \infty$ and \textit{exponentially decaying} if there are real numbers  $c,d>0$ such that $a_n \leq c \exp\{-d\, n \}$ for every $n \in \bN$. A useful property is that if $(a_n)$ is exponentially decaying, so is  the sequence $(\sum_{m \geq \zeta n} a_m)$, for any  $\zeta \in (0,1]$.

\section{Problem formulation}\label{sec:formulate}
Let $(\bS, \cS)$ be a measurable space, and let $(\Omega, \cF, \Pro)$ be a probability space on which are defined $M$ independent sequences of iid $\bS$-valued random elements, $\{X_i(n): n \in \bN\}$, $i \in [M]$, corresponding to $M$ distinct data sources. In addition, let $\{Z_i(n): n \in \bN\}$, $i \in [M]$, be $M$ independent sequences of iid $\mathrm{Uniform}[0,1]$ random variables, used for randomization purposes, and define $Z(n):=(Z_1(n),\ldots,Z_M(n))$. 

For each $i \in [M]$, the observations $\{X_i(n): n \in \bN\}$ are generated according to a density $f_i \in \Phi_i$, with respect to a $\sigma$-finite measure $\nu_i$, where the class of densities $\Phi_i$ is equipped with the $\sL^1$-norm $\|\cdot\|_1$. We say that source $i$ is \emph{anomalous} if $f_i \in \Phi_i^1$, and \emph{regular} if $f_i \in \Phi_i^0$, where $\Phi_i := \Phi_i^0 \cup \Phi_i^1$ and $\Phi_i^0 \cap \Phi_i^1 = \emptyset$. Let $\boldsymbol{\Phi} := \Phi_1 \times \cdots \times \Phi_M$ denote the product space. For each $A \subseteq [M]$, define
\begin{equation*}
\boldsymbol{\Phi}_A := \left\{(g_1,\ldots,g_M) : g_i \in \Phi_i^1 \text{ for } i \in A, \ \text{and } \ g_j \in \Phi_j^0 \text{ for } j \notin A \right\},
\end{equation*}
where $A$ stands for the subset of anomalous sources. For each $A \subseteq [M]$ and $\boldsymbol{f} := (f_1,\ldots,f_M) \in \boldsymbol{\Phi}_A$, we denote by $\Pro_{\boldsymbol{f}}$ the probability measure and by $\Exp_{\boldsymbol{f}}$ the corresponding expectation, when the underlying true densities are given by $\boldsymbol{f}$. When the dependence on $\boldsymbol{f}$ is not essential for our computations, we simply write $\Pro$ and $\Exp$.

For each $i \in [M]$ and $g \in \Phi_i$, we define the Kullback--Leibler (KL) divergence
\begin{equation*}
	D(f_i,g) := \Exp_{f_i}\left[ \log\left( \frac{f_i(X_i)}{g(X_i)} \right)\right],
\end{equation*}
where $X_i$ denotes a generic observation from source $i$. We impose the following assumptions, which are required for the development of our results:

\begin{enumerate}[(i)]
	
	\item \textit{Separation of hypotheses:}
	\begin{align}
		I_i(f_i) &:= \inf_{g \in \Phi_i^0} D(f_i,g) \in (0,\infty), \quad \forall\, i \in A, \label{inf1} \\
		J_i(f_i) &:= \inf_{g \in \Phi_i^1} D(f_i,g) \in (0,\infty), \quad \forall\, i \notin A, \label{inf2}
	\end{align}
	
	\item \textit{Continuity of the infimum maps:}
	\begin{equation}\label{inf_c}
		\lim_{f \to f_i} I_i(f) = I_i(f_i), 
		\quad 
		\lim_{f \to f_i} J_i(f) = J_i(f_i).
	\end{equation}
	
\end{enumerate}
Both assumptions depend on the true density $f_i$, which is unknown, and thus it is required that they hold locally for every $f_i \in \Phi_i$, i.e., in a neighborhood of $f_i$.

The main objective is to identify the set of anomalous sources $A$, i.e., to determine whether $f_i \in \Phi_i^1$ or $f_i \in \Phi_i^0$ for each $i \in [M]$, based on sequential observations under sampling constraints that preclude observing all sources at every time. To this end, we consider policies that consist of three components:
(i) a stopping time $T$,
(ii) a sequence $\Delta := \{\Delta(n): n \in \bN\}$, where $\Delta(n) \subseteq [M]$ denotes the current estimate of the anomalous set and $\Delta(T)$ is the final decision, and
(iii) a sequence $R := \{R(n): n \in \bN\}$, where $R(n) \subseteq [M]$ is the set of sources sampled at time $n$. We say that $R$ is a \emph{sampling rule} if $R(n)$ is $\cF_{n-1}^R$-measurable for every $n \in \bN$, where the filtration $\{\cF_n^R : n \in \bN \}$ is defined recursively by
\begin{align*}
	\cF_n^R &:= \sigma\!\left(\cF_{n-1}^R,\, Z(n),\, \{X_i(n): i \in R(n)\}\right), \quad n \in \bN, \\
	\cF_0^R &:= \sigma(Z(0)).
\end{align*}

A triplet $(R, T, \Delta)$ is called a \emph{policy} if $R$ is a sampling rule, $T$ is a stopping time with respect to $\{\cF_n^R : n \in \bN\}$, i.e., $\{T=n\} \in \cF_n^R$ for all $n$, and $\Delta(n)$ is $\cF_n^R$-measurable for every $n \in \bN$. In this case, we refer to $T$ as the \emph{stopping rule} and to $\Delta$ as the \emph{decision rule}. Our policies $(R,T,\Delta)$ must satisfy prescribed error and sampling constraints, described next.

\textit{Error constraints:} Let $\alpha, \beta \in (0,1)$. A policy $(R, T, \Delta)$ is said to belong to the class $\cC(\alpha,\beta)$ if it satisfies the error constraints
\begin{equation}\label{error_constr}
	\Pro_{\boldsymbol{f}}\!\left(|\Delta(T) \setminus A| \geq 1\right) \leq \alpha, 
	\quad 
	\Pro_{\boldsymbol{f}}\!\left(|A \setminus \Delta(T)| \geq 1\right) \leq \beta,
\end{equation}
for all $A \subseteq [M]$ and all $\boldsymbol{f} \in \boldsymbol{\Phi}_A$. 
These correspond to controlling the probabilities of at least one \emph{false alarm} and at least one \emph{missed detection}, respectively.

Prior to stating the sampling constraint, we introduce the following notation. For any sampling rule $R$, let $R_i(n)$ denote the indicator of whether source $i$ is sampled at time $n$, i.e.,
\begin{equation*}
	R_i(n) := \mathbf{1}\{ i \in R(n) \}.
\end{equation*}
We define
\begin{align}\label{nrt}
	N^R_i(n) := \sum_{m=1}^n R_i(m), 
	\qquad 
	\pi_i^R(n) := \frac{N^R_i(n)}{n},
\end{align}
where $N^R_i(n)$ is the number of times source $i$ is sampled in the first $n$ time instants, and $\pi_i^R(n)$ is the corresponding sampling proportion.

\textit{Sampling constraint:} Let $K \in (0,M]$. We say that a policy $(R, T, \Delta)$ belongs to the class $\cC(\alpha,\beta,K)$ if, in addition to satisfying the error constraint \eqref{error_constr}, it holds that
\begin{equation}\label{samp_const}
	\sum_{i=1}^{M} \Exp\!\left[ N^R_i(T) \right] \leq K \, \Exp[T].
\end{equation}
That is, the expected total number of observations up to the stopping time does not exceed $K$ times the expected stopping time. Using the identity
\begin{align*}
	\sum_{i=1}^M \Exp \left[ N_i^R(T) \right] 
	= \Exp \left[ \sum_{n=1}^T \Exp \left[ |R(n)| \,\big|\, \cF_{n-1}^R \right] \right],
\end{align*}
we observe that \eqref{samp_const} is satisfied whenever
\begin{equation}\label{samp_const_0}
	\sup_{n \in \bN,\, n \le T} \Exp \left[ |R(n)| \,\big|\, \cF_{n-1}^R \right] \leq K, \quad \text{a.s.}
\end{equation}
Condition \eqref{samp_const_0} is stronger than \eqref{samp_const}, and thus \eqref{samp_const} allows for sampling rules that do not necessarily satisfy \eqref{samp_const_0}.

Our goal is to design a policy $(R,T,\Delta) \in \cC(\alpha,\beta,K)$ whose expected stopping time attains the infimum
\begin{equation} \label{J}
	\mathcal{J}_{\boldsymbol{f}}(\alpha,\beta,K) :=  
	\inf_{(R,T,\Delta) \in \cC(\alpha,\beta,K)} \Exp_{\boldsymbol{f}}[T]
\end{equation}
to a first-order asymptotic approximation as $\alpha,\beta \to 0$, for every $\boldsymbol{f} \in \boldsymbol{\Phi}$. Throughout, we assume that there exists $r \in (0,\infty)$ such that
\begin{equation}\label{r}
	|\log(\alpha)| \sim r |\log(\beta)|, \quad \alpha,\beta \to 0.
\end{equation}

\section{Background results}\label{sec:background}

In our previous work \cite{a_prob21}, we studied sequential
anomaly identification under sampling constraints in the
simple-versus-simple setting. Each source generates independent
and identically distributed observations, independently of
the other sources, according to one of two specified densities:
$f_{0i}$ when source $i$ is regular and $f_{1i}$ when it is
anomalous. The corresponding information numbers,
$I_i:=D(f_{1i},f_{0i})$ and $J_i:=D(f_{0i},f_{1i})$, are assumed
to be positive and finite. The objective is to identify the
unknown anomalous set $A$ while controlling the probabilities
of false positive and false negative errors below prescribed
levels $\alpha$ and $\beta$. Observations are collected under the sampling
constraint \eqref{samp_const}.

The asymptotic optimality theory in \cite{a_prob21} begins
with a universal asymptotic lower bound on the expected
stopping time of any procedure satisfying the error and
sampling constraints. It then identifies target sampling
frequencies $\{c_i^*(A):i\in[M]\}$ and sufficient conditions
under which maintaining these frequencies ensures that
the lower bound is attained. These frequencies depend on
the information numbers, the sampling budget, and the
relative error levels. Specifically, it suffices that
\begin{equation*}
	\sum_{n=1}^{\infty}
	\Pro_A\!\left(
	\pi_i^R(n)<c_i^*(A)-\epsilon
	\right)<\infty,
	\qquad
	\forall\,i\in[M],\quad \forall\,\epsilon>0,
\end{equation*}
where $\pi_i^R(n)$ is the empirical sampling frequency
of source $i$ at time $n$. This condition controls the
probabilities that a source is sampled less frequently
than its target requires.

To satisfy this condition, a probabilistic sampling rule
was proposed using one Bernoulli draw per source.
At time $n+1$, source $i$ is sampled with probability
$c_i^*(\Delta_n^R)$, where $\Delta_n^R$ is the current
estimate of the anomalous set. Thus, the sampling
probabilities adapt to the observations through the
estimated configuration. Together with the stopping
and decision rules of \cite{a_prob21}, this sampling
rule attains first-order asymptotic optimality.

The asymptotically optimal expected stopping time is
determined by the rate at which evidence is accumulated.
Sampling an anomalous source $i$ with frequency $c_i^*(A)$
yields an information rate $c_i^*(A)I_i$; the corresponding
rate for a regular source $j$ is $c_j^*(A)J_j$.
When $\alpha=\beta\to0$, the target frequencies balance
these information rates, and an appropriately designed
procedure satisfies
\begin{equation}\label{eq:bg-performance}
	\begin{aligned}
		\Exp_A[T^R]
		&\sim
		\frac{|\log\alpha|}{c_i^*(A)I_i},
		&& i\in A,\\
		\Exp_A[T^R]
		&\sim
		\frac{|\log\beta|}{c_j^*(A)J_j},
		&& j\notin A.
	\end{aligned}
\end{equation}
These expressions hold for arbitrary $K\in(0,M]$, with
the dependence on $K$ incorporated in the target frequencies.
Thus, optimal sampling accounts for both the information
provided by each observation and the frequency with which
each source is observed.

In the present work, we extend these results to nonparametric
hypothesis testing, where the null and alternative hypotheses
consist of classes of densities. The information numbers and
target sampling frequencies now depend on the unknown
underlying densities. We introduce a generalized criterion
for asymptotic optimality that jointly addresses the sampling
rule, the test statistics, and the density estimators.
Specifically, it requires sampling probabilities that
dominate the target frequencies, test statistics that
accumulate evidence at the corresponding information rates,
and density estimators that converge sufficiently quickly
to the true densities. Together, these conditions guarantee
that the designed procedure attains first-order asymptotic
optimality.

\section{An illustrative example}\label{motiv}

This example illustrates how distributional misspecification
can affect the allocation of limited observations across
sources. Although accurate mean estimation suffices for
consistent hypothesis identification in this setting,
a Gaussian surrogate that matches the true mean and
variance can yield sampling frequencies that differ
substantially from the optimal frequencies determined by
the true densities. We show that the resulting allocation
can cause an arbitrarily large loss in first-order asymptotic
efficiency, even when the true likelihood ratios are used
for testing. This highlights the importance of accounting
for distributional shape in hypothesis testing under sampling
constraints and motivates the use of nonparametric density
estimation when a reliable parametric family is unavailable.

We consider two independent sources, each generating independent
and identically distributed observations. Exactly one source is
observed at each time, i.e.,
\begin{equation*}
	|R(n)|=1,
	\qquad n\geq1,
\end{equation*}
which satisfies the sampling constraint \eqref{samp_const_0} with
$K=1$. For each $i\in[2]$, let $f_i$ denote the underlying true
density, and let $f_0:=\cN(0,\sigma^2)$ denote the common
Gaussian null density. Let $\Phi_i^1$ be the class of all
$L$-Lipschitz densities on $\bR$ with mean in $[\epsilon,B]$,
variance $\sigma^2$, and fourth moment bounded by $C$,
where $L,C>0$ are fixed and sufficiently large.
We consider the hypothesis-testing problem
\begin{equation}\label{eq:motivating-hypotheses}
	\begin{aligned}
		H_i^1 &: f_i\in\Phi_i^1,\\
		H_i^0 &: f_i=f_0,
	\end{aligned}
	\qquad i\in[2],
\end{equation}
where $0<\epsilon<1\leq B$ and $\sigma>0$ are fixed.

Class $\Phi_i^1$ permits densities with the same mean and variance
but different distributional shapes and information numbers.
It therefore allows us to isolate the effect of distributional
shape on the sampling allocation.

Suppose that both sources are anomalous and that $f_i$ is
a Laplace density with location parameter $\mu_i$, and
scale parameter $\sigma/\sqrt{2}$, i.e.,
\begin{equation}\label{eq:motivating-true-densities}
	f_i
	=
	\cL\left(\mu_i,\frac{\sigma}{\sqrt{2}}\right),
	\qquad i\in[2],
\end{equation}
where
\begin{equation}\label{eq:motivating-source-means}
	\mu_i
	:=
	\begin{cases}
		1, & i=1,\\
		\epsilon, & i=2.
	\end{cases}
\end{equation}
Throughout this example, dependence on $\epsilon$ is suppressed
in the notation unless needed when taking limits.

For the purposes of this illustrative example, we use the corresponding
simple-versus-simple testing problems, in which each true alternative
density $f_i$ is treated as known. These provide an oracle benchmark
that reflects the first-order asymptotic behavior of
appropriately designed procedures for the composite problem.
The results needed to establish this behavior for composite
hypotheses are developed later in the paper.

In this simple-versus-simple setting, our previous work
\cite{a_prob21} shows that, under suitable regularity conditions,
an appropriately designed procedure with limiting sampling
frequencies satisfying
\begin{equation}\label{eq:optimal-frequency-balance}
	c_1 I_1=c_2 I_2,
	\qquad
	c_1+c_2=1,
\end{equation}
attains the first-order asymptotically optimal expected stopping time
\begin{equation}\label{eq:motivating-optimal-stopping-time}
	\Exp_{[2]}[T_{\mathrm{or}}]
	\sim
	\frac{|\log\alpha|}{c_1 I_1},
	\qquad \alpha=\beta\to0,
\end{equation}
where $T_{\mathrm{or}}$ denotes the stopping time of the oracle
procedure, and $\alpha$ and $\beta$ are the prescribed
error-probability bounds in \eqref{error_constr}.
Here and below, $\Exp_{[2]}$ denotes expectation when both sources
are anomalous with the particular densities specified in
\eqref{eq:motivating-true-densities}. The corresponding information
numbers are
\begin{equation}\label{eq:motivating-true-information}
	\begin{aligned}
		I_1
		&:=
		\dD\!\left(f_1,f_0\right)
		=
		\frac{1}{2\sigma^2}
		+\frac{\log\pi-1}{2},\\
		I_2
		&:=
		\dD\!\left(f_2,f_0\right)
		=
		\frac{\epsilon^2}{2\sigma^2}
		+\frac{\log\pi-1}{2}.
	\end{aligned}
\end{equation}

For each fixed $\epsilon>0$, sufficiently accurate estimation
of the mean is enough to distinguish the two hypotheses with
high probability. However, the mean alone does not determine
the information numbers needed to compute the optimal sampling
frequencies. Consequently, treating a Laplace alternative as
Gaussian can lead to an inefficient sampling allocation.

To illustrate this effect, suppose that a practitioner replaces
\eqref{eq:motivating-hypotheses} by the simplified Gaussian problem
\begin{equation}\label{eq:motivating-surrogate-hypotheses}
	\begin{aligned}
		H_i^{1,\mathrm{s}}
		&: f_i\in
		\left\{
		\cN(\mu,\sigma^2):\mu\in[\epsilon,B]
		\right\},\\
		H_i^{0,\mathrm{s}}
		&: f_i=f_0,
	\end{aligned}
	\qquad i\in[2].
\end{equation}
Let
\begin{equation*}
	f_i^{\mathrm{s}}
	:=
	\cN(\mu_i,\sigma^2),
	\qquad i\in[2],
\end{equation*}
denote the Gaussian surrogate densities with the correct means.
Under the true densities, the surrogate log-likelihood ratios
have mean increments
\begin{equation*}
	\begin{aligned}
		I_1^{\mathrm{s}}
		&:=
		\Exp_{f_1}\!\left[
		\log\frac{f_1^{\mathrm{s}}(X)}{f_0(X)}
		\right]
		=\frac{1}{2\sigma^2},\\
		I_2^{\mathrm{s}}
		&:=
		\Exp_{f_2}\!\left[
		\log\frac{f_2^{\mathrm{s}}(X)}{f_0(X)}
		\right]
		=\frac{\epsilon^2}{2\sigma^2}.
	\end{aligned}
\end{equation*}
Let $c_1^{\mathrm{s}}$ and $c_2^{\mathrm{s}}$ denote the sampling
frequencies computed from these surrogate information numbers,
so that
\begin{equation*}
	c_1^{\mathrm{s}}I_1^{\mathrm{s}}
	=
	c_2^{\mathrm{s}}I_2^{\mathrm{s}},
	\qquad
	c_1^{\mathrm{s}}+c_2^{\mathrm{s}}=1.
\end{equation*}
For this illustrative comparison, we use the first-order
stopping-time approximation associated with these surrogate
information rates:
\begin{equation}\label{eq:motivating-surrogate-stopping-time}
	\Exp_{[2]}[T_{\mathrm{s}}]
	\sim
	\frac{|\log\alpha|}
	{c_1^{\mathrm{s}}I_1^{\mathrm{s}}},
	\qquad \alpha=\beta\to0,
\end{equation}
for each fixed $\epsilon>0$, where $T_{\mathrm{s}}$ denotes the
stopping time of the surrogate procedure.

Under these stopping-time asymptotics, for each fixed
$\epsilon>0$, the ratio of the expected stopping times satisfies
\begin{equation}\label{eq:motivating-stopping-time-ratio}
	\cR(\epsilon)
	:=
	\lim_{\alpha\to0}
	\frac{\Exp_{[2]}[T_{\mathrm{s}}]}
	{\Exp_{[2]}[T_{\mathrm{or}}]}
	=
	\frac{c_1(\epsilon)I_1}
	{c_1^{\mathrm{s}}(\epsilon)I_1^{\mathrm{s}}},
\end{equation}
where we have made the dependence of the sampling frequencies
on $\epsilon$ explicit. The quantities $I_1$ and
$I_1^{\mathrm{s}}$ are fixed and strictly positive.
Moreover, as $\epsilon\downarrow0$, $c_1(\epsilon)$ converges
to a strictly positive limit, whereas
$c_1^{\mathrm{s}}(\epsilon)$ converges to zero.
It follows that
\begin{equation*}
	\cR(\epsilon)\longrightarrow\infty
	\qquad\text{as }\epsilon\downarrow0.
\end{equation*}

As $\epsilon$ decreases, the Gaussian surrogate treats the
second source as increasingly difficult to distinguish from
the null and allocates nearly all observations to it.
However, the true Laplace alternative remains separated from
the Gaussian null through its distributional shape, so the
optimal allocation continues to sample each source with a
positive limiting frequency. Thus, the first-order
stopping-time comparison illustrates how ignoring distributional
shape can lead to an arbitrarily large loss of efficiency.

\section{Universal asymptotic lower bound}\label{UALB}
In this section, we derive the universal asymptotic lower bound on $\mathcal{J}_{\boldsymbol{f}}(\alpha,\beta,K)$ as $\alpha, \beta \to 0$ under the scaling condition \eqref{r}. For any $A \subseteq [M]$ with $A \neq \emptyset$, we define
\begin{align*}
	I^{*}_{A}(\boldsymbol{f}) &:= \min_{i \in A} I_{i}(f_{i}), 
	\qquad 
	\hat{K}_{A}(\boldsymbol{f}) := \sum_{i \in A} \frac{I^{*}_{A}(\boldsymbol{f})}{I_{i}(f_{i})},
\end{align*}
and for any $A \subseteq [M]$ with $A \neq [M]$,
\begin{align*}
	J^{*}_{A}(\boldsymbol{f}) &:= \min_{i \notin A} J_{i}(f_{i}), 
	\qquad 
	\check{K}_{A}(\boldsymbol{f}) := \sum_{i \notin A} \frac{J^{*}_{A}(\boldsymbol{f})}{J_{i}(f_{i})}.
\end{align*}
The quantities $I^{*}_{A}(\boldsymbol{f})$ and $J^{*}_{A}(\boldsymbol{f})$ represent the minimum Kullback--Leibler divergences over the sets $\{ I_{i}(f_{i}) : i \in A \}$ and $\{ J_{i}(f_{i}) : i \notin A \}$, respectively. For $0 < |A| < M$, we define the ratio
\begin{equation*}
	s^{*}_{A}(\boldsymbol{f}) := \frac{I^{*}_{A}(\boldsymbol{f})}{J^{*}_{A}(\boldsymbol{f})}.
\end{equation*}

In the following theorem, we present the Universal Asymptotic Lower Bound (UALB), expressed in terms of factors $x_{A}(\boldsymbol{f}), y_{A}(\boldsymbol{f}) \in [0,1]$. The terms $\hat{K}_{A}(\boldsymbol{f})$ and $\check{K}_{A}(\boldsymbol{f})$ can be interpreted as the effective sampling effort required over the sources in $A$ and $A^c$, respectively, for achieving $x_{A}(\boldsymbol{f}) = 1$ and $y_{A}(\boldsymbol{f}) = 1$. In addition to characterizing the UALB, the quantities $x_A(\boldsymbol{f})$ and $y_A(\boldsymbol{f})$ play a key role in the design of asymptotically optimal sampling rules, as will be seen in following section.

\begin{theorem}\label{theo:lowerbound}
	Let $A \subseteq [M]$ and $\boldsymbol{f} \in \boldsymbol{\Phi}_A$.
	
	(i) If $|A| = 0$, then as $\alpha, \beta \to 0$,
	\begin{align*}
		\mathcal{J}_{\boldsymbol{f}}(\alpha, \beta, K) 
		\gtrsim \frac{|\log \beta|}{y_A(\boldsymbol{f})\, J^{*}_{A}(\boldsymbol{f})},
	\end{align*}
	where
	\begin{equation*}
		y_A(\boldsymbol{f}) := \left(\frac{K}{\check{K}_A(\boldsymbol{f})}\right) \wedge 1, 
		\qquad 
		x_A(\boldsymbol{f}) := 0.
	\end{equation*}
	
	(ii) If $0 < |A| < M$, then as $\alpha, \beta \to 0$ under \eqref{r},
	\begin{align*}
		\mathcal{J}_{\boldsymbol{f}}(\alpha, \beta, K)
		\gtrsim \frac{|\log \alpha|}{x_A(\boldsymbol{f})\, I_A^*(\boldsymbol{f})}
		\sim \frac{|\log \beta|}{y_A(\boldsymbol{f})\, J_A^*(\boldsymbol{f})},
	\end{align*}
	where
	\begin{equation*}
		\begin{aligned}
			x_A(\boldsymbol{f}) 
			&:= \frac{K}{\hat{K}_A(\boldsymbol{f}) + (s^{*}_{A}(\boldsymbol{f})/r)\,\check{K}_A(\boldsymbol{f})} 
			\;\wedge\; \frac{r}{s^{*}_{A}(\boldsymbol{f})} \;\wedge\; 1, \\
			y_A(\boldsymbol{f}) 
			&:= \frac{K}{\check{K}_A(\boldsymbol{f}) + (r/s^{*}_{A}(\boldsymbol{f}))\,\hat{K}_A(\boldsymbol{f})} 
			\;\wedge\; \frac{s^{*}_{A}(\boldsymbol{f})}{r} \;\wedge\; 1.
		\end{aligned}
	\end{equation*}
	
	(iii) If $|A| = M$, then as $\alpha, \beta \to 0$,
	\begin{align*}
		\mathcal{J}_{\boldsymbol{f}}(\alpha, \beta, K) 
		\gtrsim \frac{|\log \alpha|}{x_A(\boldsymbol{f})\, I^{*}_{A}(\boldsymbol{f}) + y_A(\boldsymbol{f})\, J^{*}_{A}(\boldsymbol{f})},
	\end{align*}
	where
	\begin{equation*}
		x_A(\boldsymbol{f}) := \left(\frac{K}{\hat{K}_A(\boldsymbol{f})}\right) \wedge 1, 
		\qquad 
		y_A(\boldsymbol{f}) := 0.
	\end{equation*}
\end{theorem}

\begin{IEEEproof}
	Fix $A \subseteq [M]$ and $\boldsymbol{f} \in \boldsymbol{\Phi}_A$. We distinguish two cases, depending on whether the infima in \eqref{inf1}--\eqref{inf2} are attained.
	
	\textit{Case 1:} Suppose that for each $i \in A$ (resp. $i \in A^c$) there exists a function $g_i^*$ that attains the infimum in \eqref{inf1} (resp. \eqref{inf2}). Let $\boldsymbol{g}^* := (g_1^*, \ldots, g_M^*)$. Consider the following $M$-binary (simple vs. simple) hypothesis testing problem:
	\begin{align}\label{Msh}
		\{ X_i(n) : n \in \bN \} \sim f_i 
		\quad \text{vs.} \quad 
		\{ X_i(n) : n \in \bN \} \sim g_i^*, 
		\qquad i \in [M].
	\end{align}
	Let $\cC^{\boldsymbol{f},\boldsymbol{g}^*}_{\text{simple}}(\alpha,\beta,K)$ denote the class of policies satisfying the error constraint \eqref{error_constr} and the sampling constraint \eqref{samp_const} for the problem \eqref{Msh}. Any policy $(R,T,\Delta) \in \cC(\alpha,\beta,K)$ for the composite problem also satisfies the same constraints for the simple problem \eqref{Msh}. Therefore,
	\begin{align}\label{sub}
		\cC(\alpha,\beta,K) \subseteq \cC^{\boldsymbol{f},\boldsymbol{g}^*}_{\text{simple}}(\alpha,\beta,K).
	\end{align}
	By \cite[Theorem 5]{a_prob21}, for any $(R,T,\Delta) \in \cC^{\boldsymbol{f},\boldsymbol{g}^*}_{\text{simple}}(\alpha,\beta,K)$, we have $\Exp_{\boldsymbol{f}}[T]$ asymptotically lower bounded by the UALB in Theorem~\ref{theo:lowerbound}. Combining this with \eqref{sub} yields the desired result.
	
	\textit{Case 2:} Suppose that for some $i \in [M]$, the infimum in \eqref{inf1} or \eqref{inf2} is not attained. Under the \textit{separation of hypotheses} assumption, for any $\epsilon > 0$ there exists $g_i^\epsilon \in \Phi_i$ such that
	\begin{equation*}
		\begin{aligned}
			I_i(f_i) < D(f_i,g_i^\epsilon) \le (1+\epsilon) I_i(f_i), &\quad i \in A,\\
			J_i(f_i) < D(f_i,g_i^\epsilon) \le (1+\epsilon) J_i(f_i), &\quad i \in A^c.
		\end{aligned}
	\end{equation*}
	Let $\boldsymbol{g}^\epsilon := (g_1^\epsilon,\ldots,g_M^\epsilon)$ and consider the corresponding simple testing problem,
	\begin{align*}
		\{ X_{i}(n) \, :\, n \in \bN \} \sim f_i \quad \mbox{vs} \quad \{ X_{i}(n) \, :\, n \in \bN \} \sim g^{\epsilon}_{i}, \quad i \in [M].
	\end{align*}
	Repeating the argument of Case~1, we obtain
	\begin{equation*}
		\Exp_{\boldsymbol{f}}[T] \gtrsim (1 - o(\epsilon))\, \text{UALB},
	\end{equation*}
	since $I_i(f_i)$ and $J_i(f_i)$ are replaced by $D(f_i,g_i^\epsilon)$. Letting $\epsilon \to 0$ completes the proof.
\end{IEEEproof}

\section{A family of policies}\label{prob_rules}
In this section, we present a family of policies $(R, T, \Delta)$ that, under an appropriate design, achieve first-order asymptotically optimal performance \eqref{J} by attaining the UALB derived in Section~\ref{UALB}. 

\subsection{Test Statistics and Estimators}

Throughout this work, we denote by $\Lambda_i^R(n)$ the test statistic for source $i \in [M]$, based on the observations collected from source $i$ up to time $n$ under the sampling rule $R$. The test statistic $\Lambda_i^R(n)$ is designed to quantify the evidence in favor of the hypothesis that source $i$ is anomalous. In particular, positive values of $\Lambda_i^R(n)$ indicate that $i \in A$, whereas negative values indicate that $i \notin A$. Thus, $\Lambda_i^R(n)$ serves as a sequential decision statistic, analogous to a log-likelihood ratio in classical sequential testing.

We also denote by $\hat{\boldsymbol{f}}_n$ the estimator of the true density vector at time $n$, constructed from the observations up to time $n$. The estimator approximates the unknown densities governing each source and is particularly important in the nonparametric setting where the density models are not known a priori. 

The test statistics we apply are the Generalized Log-Likelihood Ratio (GLLR) and the Adaptive Log-Likelihood Ratio (ALLR). The key difference is that GLLR is defined via maximization over the hypothesis sets, whereas ALLR is constructed recursively using sequential estimates of the underlying densities. The density estimators we apply are the Maximum Likelihood Estimator (MLE), and the Kernel Density Estimator (KDE). 

\subsection{Stopping and Decision rules}\label{st_dec}
The stopping rule is defined as
\begin{equation}\label{stp0m}
	T^R := \inf\{ n \geq 1 : \Lambda_i^R(n) \notin (-a,b), \ \forall i \in [M] \},
\end{equation}
i.e., sampling continues as long as there exists $i \in [M]$ such that $\Lambda_i^R(n) \in (-a,b)$. At each time $n$, the estimate of the anomalous set is given by
\begin{equation}\label{delta0m}
	\Delta^R(n) := \{ i \in [M] : \Lambda_i^R(n) > 0 \}, \quad n \in \bN,
\end{equation}
and the final decision is $\Delta^R(T^R)$. We assume that the test statistics, the estimator, and the thresholds $a,b$ are chosen so that $(R,T^R,\Delta^R)$ satisfies the error constraint \eqref{error_constr}. In addition, for asymptotic optimality, we require that
\begin{equation}\label{ab}
	a \sim |\log(\beta)|, \qquad b \sim |\log(\alpha)|, \quad \text{as } \alpha, \beta \to 0.
\end{equation}
Accordingly, any proposed choice of test statistic and estimator must be shown to satisfy these conditions.

\subsection{Probabilistic Sampling Rules}

We introduce a family of sampling rules, referred to as
\emph{probabilistic sampling rules}. Under an appropriate design,
these rules guarantee that the policy $(R,T^R,\Delta^R)$, with $T^R$
and $\Delta^R$ defined in \eqref{stp0m}--\eqref{delta0m}, achieves
the first-order asymptotic optimality in \eqref{J} by attaining the
UALB.

A sampling rule $R$ is called \emph{probabilistic} if, for each
$i\in[M]$, there exists a measurable function
\begin{equation*}
	c_i^R:\boldsymbol{\Phi}\to[0,1],
\end{equation*}
where $c_i^R(\boldsymbol{g})$ denotes the probability of sampling
source $i$ at time $n+1$, given that
$\hat{\boldsymbol f}_n=\boldsymbol{g}$. Specifically,
\begin{equation}\label{bern_1}
	R_i(n+1)
	=
	\mathbf{1}\left\{
	Z_i(n)\leq
	c_i^R\bigl(\hat{\boldsymbol f}_n\bigr)
	\right\},
	\qquad
	i\in[M],\quad n\in\bN,
\end{equation}
where $\{Z_i(n):i\in[M],\,n\in\bN\}$ are iid uniform random
variables on $[0,1]$, independent of the observations and the past
sampling decisions.

We next provide a sufficient condition under which a probabilistic
sampling rule satisfies the sampling constraint \eqref{samp_const}.

\begin{proposition}\label{ck}
	If $R$ is a probabilistic sampling rule such that for all $\boldsymbol{g} \in \boldsymbol{\Phi}$,
	\begin{align}\label{sum_less_K} 
		\sum_{i=1}^M c^R_i(\boldsymbol{g}) \leq K,
	\end{align}
	then \eqref{samp_const} holds for any $\{\cF_n^R\}$-stopping time $T$.
\end{proposition}

\begin{IEEEproof}
	For any $n \in \bN$,
	\begin{align*}
		\Exp \left[ |R(n)| \,\big|\, \cF_{n-1}^R \right]
		= \sum_{i=1}^M \Pro \left(R_i(n)=1 \,\big|\, \cF_{n-1}^R \right)
		= \sum_{i=1}^M c^R_i\bigl(\hat{\boldsymbol{f}}_{n-1}\bigr) \leq K.
	\end{align*}
\end{IEEEproof}

The following condition ensures that every source is sampled with a
uniformly positive probability, independently of the current density
estimates. Consequently, the empirical sampling frequency of every
source remains bounded away from zero with high
probability.

\begin{lemma}\label{exp_d}
	Let $A \subseteq [M]$, $\boldsymbol{f} \in \boldsymbol{\Phi}_A$, and let $R$ be a probabilistic sampling rule. If, for each $i \in [M]$,
	\begin{equation}\label{inf_cond}
		\min_{D \subseteq [M]} \inf_{\boldsymbol{g} \in \boldsymbol{\Phi}_D} c_i^R(\boldsymbol{g}) > 0,
	\end{equation}
	then there exists $\rho \in (0,1)$ such that, for every $i \in [M]$, the sequence
	\[
	\left\{ \Pro_{\boldsymbol{f}}\!\left(\pi_i^R(n)<\rho\right) : n \in \bN \right\}
	\]
	is exponentially decaying.
\end{lemma}

\begin{IEEEproof}
	Fix $A \subseteq [M]$, $\boldsymbol{f} \in \boldsymbol{\Phi}_A$, a sampling rule $R$, and $i \in [M]$. Define
	\begin{equation}\label{stq}
		\begin{aligned}
			\widetilde{\pi}_i^R(n)
			:= \frac{1}{n} \sum_{m=1}^{n} \left(R_i(m)- c_i^R\!\left(\hat{\boldsymbol{f}}_{m-1}\right)\right) = \pi_i^R(n) - \frac{1}{n} \sum_{m=1}^{n} c_i^R\!\left(\hat{\boldsymbol{f}}_{m-1}\right)
			\leq \pi_i^R(n) - \zeta,
		\end{aligned}
	\end{equation}
	where
	\[
	\zeta := \min_{D \subseteq [M]} \inf_{\boldsymbol{g} \in \boldsymbol{\Phi}_D} c_i^R(\boldsymbol{g}) > 0.
	\]
	Fix $\epsilon \in (0,\zeta)$. For any $\rho \in (0,\zeta-\epsilon)$, \eqref{stq} implies
	\begin{equation*}
		\begin{aligned}
			\Pro_{\boldsymbol{f}}\left( \pi_i^R(n)<\rho\right)
			\leq \Pro_{\boldsymbol{f}}\left( \widetilde{\pi}_i^R(n)<\rho-\zeta \right)\leq \Pro_{\boldsymbol{f}}\left(
			\sum_{m=1}^{n} \left(R_i(m)-c_i^R(\hat{\boldsymbol{f}}_{m-1})\right)
			< -\epsilon n \right).
		\end{aligned}
	\end{equation*}
	By \eqref{bern_1},
	\begin{equation*}
		\Exp_{\boldsymbol{f}}\!\left[
		R_i(n)-c_i^R(\hat{\boldsymbol{f}}_{n-1})
		\,\big|\, \cF_{n-1}^R
		\right]
		=
		\Pro_{\boldsymbol{f}}\!\left(R_i(n)=1\,\big|\,\cF_{n-1}^R\right)
		-c_i^R(\hat{\boldsymbol{f}}_{n-1})
		=0.
	\end{equation*}
	Therefore,
	\[
	\left\{R_i(n)-c_i^R(\hat{\boldsymbol{f}}_{n-1}) : n \geq 1\right\}
	\]
	is a martingale difference sequence, uniformly bounded by $1$. The result follows from the Azuma--Hoeffding inequality.
\end{IEEEproof}

\section{Asymptotic Optimality}\label{asyo}

In view of the UALB in Theorem~\ref{theo:lowerbound}, we now
provide sufficient conditions under which a policy
$(R,T^R,\Delta^R)$ attains this bound and therefore satisfies
the first-order asymptotic optimality criterion in \eqref{J}.
We first introduce the target sampling frequencies
$\{c_i^*(\boldsymbol{f}):i\in[M]\}$ through
\begin{equation}\label{c_star}
	\begin{aligned}
		c_i^*(\boldsymbol{f}) I_i(f_i)
		&= x_A(\boldsymbol{f})\,I_A^*(\boldsymbol{f}),
		\qquad i\in A,\\
		c_i^*(\boldsymbol{f}) J_i(f_i)
		&= y_A(\boldsymbol{f})\,J_A^*(\boldsymbol{f}),
		\qquad i\notin A,
	\end{aligned}
\end{equation}
where $x_A(\boldsymbol{f})$ and $y_A(\boldsymbol{f})$ are
defined in Theorem~\ref{theo:lowerbound}. By the definition
of $r$ in \eqref{r} and the expressions for
$x_A(\boldsymbol{f})$ and $y_A(\boldsymbol{f})$, whenever
$0<|A|<M$,
\begin{equation*}
	x_A(\boldsymbol{f})\,I_A^*(\boldsymbol{f})
	=
	r\,y_A(\boldsymbol{f})\,J_A^*(\boldsymbol{f}).
\end{equation*}
Together with \eqref{c_star}, this implies
\begin{equation*}
	c_i^*(\boldsymbol{f}) I_i(f_i)
	=
	r\,c_j^*(\boldsymbol{f}) J_j(f_j),
	\qquad i\in A,\quad j\notin A.
\end{equation*}
Thus, the target frequencies balance the information
accumulation rates within each group and relate the rates
of the two groups through the factor $r$. 

The assumptions \eqref{inf1}--\eqref{inf2} ensure that the
information numbers in \eqref{c_star} are strictly positive
and finite. Suppose, in addition, that the target sampling
frequencies satisfy the uniform lower bound
\begin{equation*}
	\min_{D\subseteq[M]}
	\inf_{\boldsymbol{g}\in\boldsymbol{\Phi}_D}
	c_i^*(\boldsymbol{g})>0, \quad \forall\, i \in [M],
\end{equation*}
then any sampling rule satisfying
\begin{equation*}
	c_i^R(\boldsymbol{g})\geq c_i^*(\boldsymbol{g}),
	\qquad
\forall\,	i\in[M],\quad \forall\, D\subseteq[M],\quad \forall\,
	\boldsymbol{g}\in\boldsymbol{\Phi}_D,
\end{equation*}
also has sampling probabilities uniformly bounded away
from zero. According to
Lemma~\ref{exp_d}, it follows that the empirical sampling
frequency of every source remains bounded away from zero
with high probability.

The following theorem establishes asymptotic optimality under
a suitable asymptotic growth condition on the test statistics.
This condition ensures that evidence accumulates at rates
sufficient to attain the UALB, as determined by the target
sampling frequencies $\{c_i^*(\boldsymbol{f}):i\in[M]\}$
and the corresponding information numbers.

\begin{theorem}\label{theo:AO} 
	Let $A \subseteq [M]$, $\boldsymbol{f}\in \boldsymbol{\Phi}_{A}$, and let $R$ be a sampling rule that satisfies \eqref{samp_const} with $T=T^R$, defined in \eqref{stp0m} with thresholds $a,b$ satisfying \eqref{ab}. If for any $\epsilon > 0$ the test statistics satisfy
	\begin{equation}\label{amp1}
		\begin{aligned}
			&\sum_{n=1}^{\infty} \Pro_{\boldsymbol{f}} \left( 
			\frac{1}{n} \Lambda^R_i(n) <
			(c_i^*(\boldsymbol{f})-\epsilon) I_i(f_i) - \epsilon  
			\right) < \infty, \quad \forall i \in A, \\
			&\sum_{n=1}^{\infty} \Pro_{\boldsymbol{f}} \left( 
			-\frac{1}{n} \Lambda^R_j(n) <
			(c_j^*(\boldsymbol{f})-\epsilon) J_j(f_j) - \epsilon  
			\right) < \infty, \quad \forall j \notin A,
		\end{aligned}
	\end{equation}
	then $R$ is asymptotically optimal under $\Pro_{\boldsymbol{f}}$.
\end{theorem}

\begin{IEEEproof}
	We prove the case $0 < |A| < M$, as the cases $|A|=0$ and $|A|=M$ follow similarly. Fix $\epsilon > 0$, and let $a,b > 0$ be the thresholds in \eqref{stp0m} that satisfy  \eqref{ab}. Define
	\begin{equation}
		N_\epsilon(a,b) := \max_{i \in A,\, j \notin A} 
		\left\{
		\frac{a}{(c_j^*(\boldsymbol{f})-\epsilon) J_j(f_j)-\epsilon}, \;
		\frac{b}{(c_i^*(\boldsymbol{f})-\epsilon) I_i(f_i)-\epsilon}
		\right\}.
	\end{equation}
	
	Then
	\begin{equation}\label{expt200}
		\Exp_{\boldsymbol{f}}[T^R] \leq N_\epsilon(a,b) + \sum_{n > N_\epsilon(a,b)} \Pro_{\boldsymbol{f}}(T^R > n).
	\end{equation}
	
	To prove the claim, it suffices to show that
	\begin{equation}\label{tsh2}
		\sum_{n > N_\epsilon(a,b)} \Pro_{\boldsymbol{f}}(T^R > n) < \infty,
	\end{equation}
	then letting $a,b \to \infty$ in \eqref{expt200}, and using \eqref{ab}, we obtain
	\begin{equation*}
		\Exp_{\boldsymbol{f}}[T^R]
		\lesssim
		\max_{i \in A,\, j \notin A}
		\left\{
		\frac{|\log \beta|}{(c_j^*(\boldsymbol{f})-\epsilon) J_j(f_j)-\epsilon},
		\frac{|\log \alpha|}{(c_i^*(\boldsymbol{f})-\epsilon) I_i(f_i)-\epsilon}
		\right\}.
	\end{equation*}
	Letting $\epsilon \to 0$, it yields
	\begin{equation*}
		\Exp_{\boldsymbol{f}}[T^R]
		\lesssim
		\frac{|\log \beta|}{y_A(\boldsymbol{f}) J_A^*(\boldsymbol{f})}
		\sim
		\frac{|\log \alpha|}{x_A(\boldsymbol{f}) I_A^*(\boldsymbol{f})},
	\end{equation*}
	which matches the UALB.
	
	In order to prove \eqref{tsh2}, we note that on the event $\{T^R > n\}$, there exists either $j \notin A$ such that $\Lambda_j^R(n) > -a$, or $i \in A$ such that $\Lambda_i^R(n) < b$. Hence, by the union bound,
	\begin{equation}\label{nty}
		\Pro_{\boldsymbol{f}}(T^R > n)
		\leq
		\sum_{j \notin A} \Pro_{\boldsymbol{f}}(-\Lambda_j^R(n) < a)
		+
		\sum_{i \in A} \Pro_{\boldsymbol{f}}(\Lambda_i^R(n) < b).
	\end{equation}
	
	For $n > N_\epsilon(a,b)$, we have
	\[
	a < n \big((c_j^*(\boldsymbol{f})-\epsilon) J_j(f_j)-\epsilon\big),
	\quad
	b < n \big((c_i^*(\boldsymbol{f})-\epsilon) I_i(f_i)-\epsilon\big),
	\]
	and therefore
	\begin{align}\label{series_bound10}
		\sum_{n > N_\epsilon(a,b)} \Pro_{\boldsymbol{f}}(T^R > n)
		\leq &
		\sum_{j \notin A} \sum_{n > N_\epsilon(a,b)} \Pro_{\boldsymbol{f}} \left( 
		-\frac{1}{n} \Lambda^R_j(n) <
		(c_j^*(\boldsymbol{f})-\epsilon) J_j(f_j) - \epsilon
		\right) \nonumber \\
		&+
		\sum_{i \in A} \sum_{n > N_\epsilon(a,b)} \Pro_{\boldsymbol{f}} \left( 
		\frac{1}{n} \Lambda^R_i(n) <
		(c_i^*(\boldsymbol{f})-\epsilon) I_i(f_i) - \epsilon
		\right).
	\end{align}
	By assumption \eqref{amp1}, both series converge, which proves \eqref{tsh2}. 
\end{IEEEproof}

Next, we provide sufficient conditions under which a policy $(R,T^R,\Delta^R)$ attains the infimum \eqref{J}. Since the test statistics $\Lambda_i^R$ depend on both the sampling rule $R$ and the estimator, we decompose the asymptotic optimality condition into three components: a condition on the sampling frequencies, a condition on the asymptotic behavior of the test statistics, and a condition on the convergence of the estimator.

The first condition requires that, in the long run, the sampling rule allocates to each source $i \in [M]$ a frequency that is at least $c_i^*(\boldsymbol{f})$, as defined in \eqref{c_star}. The second condition imposes a growth requirement on the
test statistics in the spirit of a one-sided law of large
numbers. For each source and each $\epsilon>0$, it requires
that the probabilities of meeting the stated sampling-frequency
bound while falling below the corresponding lower bound on
the rate of evidence accumulation be summable. Thus, by the
Borel--Cantelli lemma, such shortfalls occur only finitely
often almost surely. The third condition controls the rate at which the estimator $\hat{\boldsymbol{f}}_n$ converges to the true densities. In particular, for each $i \in [M]$ and $\epsilon > 0$, define
\begin{equation}\label{vva}
V_{\epsilon}(f_i) :=\{h \in \Phi_i : \|f_i -h\|_1 < \epsilon\}
\end{equation}
\begin{equation}\label{tta}
	\tau^R_i(\epsilon)
	:=
	\inf\left\{
	n \in \bN : \hat{f}_{i,m} \in V_{\epsilon}(f_i)
	\; \text{ for all } m \ge n
	\right\},
\end{equation}
which is the first time after which the estimator remains within $\epsilon$ distance from the true density. The requirement is essentially that $\Exp_{\boldsymbol{f}}[\tau^R_i(\epsilon)] < \infty$ for every $\epsilon > 0$, which is stronger than almost sure uniform convergence and is commonly referred to as \emph{quick} convergence; see \cite[Subsec.~2.4.3]{tartakovsky2014sequential}.

\begin{theorem}\label{theo:prob_ao}  
	Let $A \subseteq [M]$, $\boldsymbol{f}\in \boldsymbol{\Phi}_A$, and let $R$ be a sampling rule that satisfies \eqref{samp_const} with $T=T^R$ defined in \eqref{stp0m} with thresholds $a,b$ satisfying \eqref{ab}. If
	\begin{enumerate}[(i)]
		\item for each $i \in [M]$, the functional $c_i^R$ satisfies
		\begin{equation}\label{stt1}
			c_i^R(\boldsymbol{g}) \geq c_i^*(\boldsymbol{g}),
		\end{equation}
		for every $D \subseteq [M]$ and every $\boldsymbol{g} \in \boldsymbol{\Phi}_D$;
		
		\item for every $\epsilon > 0$, the test statistics satisfy
		\begin{equation}\label{stt2}
			\begin{aligned}
				&\sum_{n=1}^{\infty} \Pro_{\boldsymbol{f}} \left( 
				\frac{1}{n} \Lambda^R_i(n) <
				(c_i^*(\boldsymbol{f})-\epsilon) I_i(f_i) - \epsilon,\,
				\pi^R_i(n) \geq c_i^*(\boldsymbol{f})-\epsilon 
				\right) < \infty, \quad \forall i \in A,\\
				&\sum_{n=1}^{\infty} \Pro_{\boldsymbol{f}} \left( 
				-\frac{1}{n} \Lambda^R_j(n) <
				(c_j^*(\boldsymbol{f})-\epsilon) J_j(f_j) - \epsilon,\,
				\pi^R_j(n) \geq c_j^*(\boldsymbol{f})-\epsilon 
				\right) < \infty, \quad \forall j \notin A;
			\end{aligned}
		\end{equation}
		
		\item for the estimator $\{\hat{\boldsymbol{f}}_n : n \in \bN\}$, for every $\rho \in (0,1]$, $\epsilon > 0$,
		\begin{equation}\label{stt3}
			\sum_{n=1}^{\infty} 
			n\, \Pro_{\boldsymbol{f}}\left( \hat{f}_{i,n} \notin V_{\epsilon}(f_i),\,
			\pi^R_i(n) \geq \rho
			\right) < \infty, 
			\quad \forall i \in [M],
		\end{equation}
	\end{enumerate}
	then $R$ is asymptotically optimal under $\Pro_{\boldsymbol{f}}$.
\end{theorem}

\begin{IEEEproof}
	In view of Theorem \ref{theo:AO}, it suffices to verify \eqref{amp1}. We prove the condition for $i \in A$, as the proof for $j \notin A$ is analogous. For any $\epsilon>0$,
	\begin{align*}
		&\Pro_{\boldsymbol f}\left(
		\frac{1}{n}\Lambda_i^R(n)
		<
		(c_i^*(\boldsymbol f)-\epsilon)I_i(f_i)-\epsilon
		\right) \\
		&\leq
		\Pro_{\boldsymbol f}\left(
		\frac{1}{n}\Lambda_i^R(n)
		<
		(c_i^*(\boldsymbol f)-\epsilon)I_i(f_i)-\epsilon,\,
		\pi_i^R(n)\geq c_i^*(\boldsymbol f)-\epsilon
		\right) +
		\Pro_{\boldsymbol f}\left(
		\pi_i^R(n)<c_i^*(\boldsymbol f)-\epsilon
		\right).
	\end{align*}
	The first term is summable by \eqref{stt2}. It remains to show that
	\begin{equation}\label{pi_summable}
		\sum_{n=1}^\infty
		\Pro_{\boldsymbol f}\left(
		\pi_i^R(n)<c_i^*(\boldsymbol f)-\epsilon
		\right)<\infty.
	\end{equation}
	
	We write
	\begin{align*}
		\pi_i^R(n)
		&=
		\frac{1}{n}\sum_{m=1}^n
		\left(R_i(m)-c_i^R(\hat{\boldsymbol f}_{m-1})\right)
		+
		\frac{1}{n}\sum_{m=1}^n
		c_i^R(\hat{\boldsymbol f}_{m-1}) \\
		&\geq
		\frac{1}{n}\sum_{m=1}^n
		\left(R_i(m)-c_i^R(\hat{\boldsymbol f}_{m-1})\right)
		+
		\frac{1}{n}\sum_{m=1}^n
		c_i^*(\hat{\boldsymbol f}_{m-1}),
	\end{align*}
	where the inequality follows from \eqref{stt1}. Hence,
	\begin{align*}
		&\Pro_{\boldsymbol f}\left(
		\pi_i^R(n)<c_i^*(\boldsymbol f)-\epsilon
		\right)\\
		&\leq \Pro_{\boldsymbol f}\left(\frac{1}{n}
		\sum_{m=1}^n
		\left(R_i(m)-c_i^R(\hat{\boldsymbol f}_{m-1})\right)
		<-\frac{\epsilon}{2}
		\right)+
		\Pro_{\boldsymbol f}\left(\frac{1}{n}
		\sum_{m=1}^n
		c_i^*(\hat{\boldsymbol f}_{m-1})
		<c_i^*(\boldsymbol f)-\frac{\epsilon}{2}
		\right)\\
		&\leq
		\Pro_{\boldsymbol f}\left(
		\sum_{m=1}^n
		\left(R_i(m)-c_i^R(\hat{\boldsymbol f}_{m-1})\right)
		<-\frac{\epsilon}{2}n
		\right)+
		\Pro_{\boldsymbol f}\left(
		\sum_{m=1}^n
		\left(c_i^*(\hat{\boldsymbol f}_{m-1})-c_i^*(\boldsymbol f)\right)
		<-\frac{\epsilon}{2}n
		\right).
	\end{align*}
	
	The first probability is exponentially decaying by Azuma-Hoeffding inequality, since by \eqref{bern_1}
	\[
	\left\{
	R_i(m)-c_i^R(\hat{\boldsymbol f}_{m-1}) : m\geq1
	\right\}
	\]
	is a martingale difference sequence, uniformly bounded by $1$. Therefore, it remains to control the second probability. By assumption of \textit{continuity of the infimum map} and the connection of $x_A(\boldsymbol{f})$, $y_A(\boldsymbol{f})$ with $I_i(f_i)$ and $J_i(f_i)$, $i \in [M]$, it follows that $c_i^*$ is continuous at $\boldsymbol f$, and thus there exists $\delta>0$ such that for any $\boldsymbol{g}\in \boldsymbol{\Phi}$
	\[
	\max_{j\in[M]}\|f_j(x)-g_j(x)\|_1<\delta
	\quad \Rightarrow \quad
	|c_i^*(\boldsymbol g)-c_i^*(\boldsymbol f)|<\epsilon/4.
	\]
	Let
	\[
	\widetilde{\tau}:=\max_{j\in[M]}\tau^R_j(\delta).
	\]
	Then, for all $m>\widetilde{\tau}$,
	\[
	|c_i^*(\hat{\boldsymbol f}_{m-1})-c_i^*(\boldsymbol f)|<\epsilon/4.
	\]
	Since for all $m \leq \widetilde{\tau}$ we have $|c_i^*(\hat{\boldsymbol f}_{m-1})-c_i^*(\boldsymbol f)|\leq 1$, it holds
	\[
	\Pro_{\boldsymbol f}\left(
	\sum_{m=1}^n
	\left|c_i^*(\hat{\boldsymbol f}_{m-1})-c_i^*(\boldsymbol f)\right|
	>\frac{\epsilon}{2}n
	\right)
	\leq
	\Pro_{\boldsymbol f}\left(
	\widetilde{\tau}+\frac{\epsilon}{4}n>\frac{\epsilon}{2}n
	\right).
	\]
	Thus, it suffices to show that, for every $j\in[M]$,
	\[
	\sum_{n=1}^\infty
	\Pro_{\boldsymbol f}\left(
	\tau^R_j(\delta)>\frac{\epsilon}{4}n
	\right)<\infty,
	\]
	which follows by Proposition \ref{tpro}.
\end{IEEEproof}

An important implication of conditions \eqref{stt1} and
\eqref{stt3} of Theorem~\ref{theo:prob_ao} is that the
estimator $\hat{\boldsymbol{f}}_n$ converges quickly to
the true density vector $\boldsymbol{f}$, in the sense
that the random times $\tau_i^R(\epsilon)$, $i\in[M]$,
defined in \eqref{tta}, have summable tail probabilities
and hence finite expectations.

\begin{proposition}\label{tpro}
	Let $A \subseteq [M]$, $\boldsymbol{f}\in \boldsymbol{\Phi}_A$, and let $R$ be a sampling rule. If assumptions (i) and (iii) of Theorem \ref{theo:prob_ao} hold then 
	\begin{equation}\label{tec}
		\sum_{n=1}^{\infty} \Pro_{\boldsymbol{f}}\big(\tau^R_i(\epsilon) > c\, n\big) < \infty,
	\end{equation}
	for any $\epsilon > 0$ and $c > 0$.
\end{proposition}

\begin{IEEEproof}
The proof is presented in Appendix \ref{app0}.
\end{IEEEproof}

As a corollary to Theorem \ref{theo:prob_ao}, we specify a particular choice of the functional $c_i^R$, $i \in [M]$, that yields asymptotic optimality.

\begin{corollary}
	Let $A \subseteq [M]$, $\boldsymbol{f} \in \boldsymbol{\Phi}_A$, and let $R$ be a sampling rule. If, for each $i \in [M]$,
	\begin{equation*}
		c_i^R(\boldsymbol{g}) = c_i^*(\boldsymbol{g}),
	\end{equation*}
	for every $D \subseteq [M]$ and every $\boldsymbol{g} \in \boldsymbol{\Phi}_D$, and if conditions \eqref{stt2} and \eqref{stt3} of Theorem~\ref{theo:prob_ao} are satisfied, then $R$ is asymptotically optimal.
\end{corollary}

\begin{IEEEproof}
	By Theorem~\ref{theo:prob_ao}, it suffices to verify that $R$ satisfies the sampling constraint \eqref{samp_const} with $T=T^R$. By Proposition~\ref{ck}, it is enough to show that \eqref{sum_less_K} holds. Evaluating at $\boldsymbol{g}=\boldsymbol{f}$ and using the definition \eqref{c_star}, we obtain
	\begin{equation*}
		\begin{aligned}
			\sum_{i=1}^{M} c_i^*(\boldsymbol{f})
			&= x_A(\boldsymbol{f}) \sum_{i \in A} \frac{I_A^*(\boldsymbol{f})}{I_i(f_i)}
			+ y_A(\boldsymbol{f}) \sum_{i \notin A} \frac{J_A^*(\boldsymbol{f})}{J_i(f_i)} \\
			&= x_A(\boldsymbol{f})\, \hat{K}_A(\boldsymbol{f})
			+ y_A(\boldsymbol{f})\, \check{K}_A(\boldsymbol{f})
			\le K,
		\end{aligned}
	\end{equation*}
	where the last inequality follows from the definitions of $x_A(\boldsymbol{f})$ and $y_A(\boldsymbol{f})$ in Theorem~\ref{theo:lowerbound}. 
\end{IEEEproof}

\section{Suitable Test Statistics}\label{test_st}
In this section, we introduce the generalized log-likelihood ratio
(GLLR) and the adaptive log-likelihood ratio (ALLR) and establish
conditions under which each statistic satisfies the growth
requirement \eqref{stt2} of Theorem~\ref{theo:prob_ao}.

To this end, we impose regularity conditions on the density
classes $\Phi_i$, $i\in[M]$, similar in spirit to those verified
by Pavlov \cite{pavlov91} for multidimensional exponential families
under suitable compactness assumptions. For a regular minimal
exponential family restricted to a compact subset of its natural
parameter space, Proposition~\ref{lem:homeomorphism} shows that
the map from parameters to densities is a homeomorphism when the
density class is equipped with the $\sL^1$ topology. This correspondence motivates formulating our regularity
conditions directly in terms of densities and their $\sL^1$
neighborhoods. We can therefore impose these conditions on
non-parametric density classes $\Phi_i$, $i\in[M]$, without
assuming a finite-dimensional parametrization. Proposition~\ref{prop:expf} verifies these conditions for the exponential families.

\subsection{Regularity Conditions}\label{g_ass}
Analogously to the density regularity conditions in \cite{pavlov91}, we impose the following regularity conditions, which will be used in the analysis of the asymptotic properties of the proposed test statistics. For each \(i \in [M]\) and each \(g \in \Phi_i\), define
\begin{equation*}
	p_{\epsilon}^{+}(g,X_i)
	:=
	\sup_{h \in V_{\epsilon}(g)} h(X_i),
\end{equation*}
where
\begin{equation*}
	V_{\epsilon}(g)
	:=
	\left\{
	h \in \Phi_i :
	\|g-h\|_1<\epsilon
	\right\}.
\end{equation*}
We assume that the following conditions hold for every \(i \in [M]\):

\begin{enumerate}[(a)]
	
	\item The class of densities \(\Phi_i\) is compact with respect to the $\sL^1$ norm.
	
	\item Every \(g \in \Phi_i\) satisfies $g(x)>0$ for all $x\in\bS$.
	
	\item The Kullback--Leibler divergence satisfies
	\begin{equation}
		\label{pos}
		D(f_i,g)>0,
		\qquad
		      \forall\, g\in\Phi_i,\; \; g\neq f_i.
	\end{equation}
	
	\item The mapping $g\mapsto D(f_i,g)$ is continuous on \(\Phi_i\).
	
	\item There exist \(\mathfrak{q}>2\), \(\epsilon>0\), and \(c>0\) such that
	\begin{equation}
		\label{moment}
		\Exp_{\boldsymbol f}\!\left[
		\sup_{h\in V_{\epsilon}(g)}
		|\log h(X_i)|^{\mathfrak{q}}
		\right]
		\le c, \qquad \forall\; g, \, f_i \in \Phi_i.
	\end{equation}
	
	\item There exist $\epsilon>0$ and $c>0$ such that
	\begin{equation}\label{eq:ratio}
		\Exp_{f_i}\!\left[
		\frac{p_\epsilon^+(g,X_i)}{f_i(X_i)}
		\right]
		=
		\int_{\bS} p_\epsilon^+(g,x)\,d\nu_i(x)
		\leq c, \qquad \forall\; g, \, f_i \in \Phi_i.
	\end{equation}
	
	\item The function
	\begin{equation*}
		H_{\boldsymbol f}(g,\epsilon)
		:=
		\Exp_{\boldsymbol f}\!\left[
		\sup_{h\in V_{\epsilon}(g)}
		\left|
		\log h(X_i)-\log g(X_i)
		\right|
		\right]
	\end{equation*}
	converges to $0$ as \(\epsilon\to0\), uniformly over \(g\in\Phi_i\).
\end{enumerate}

Condition (a) ensures that $\Phi_i$ is compact and hence
totally bounded in the $\sL^1$ norm. Condition (b) ensures
that the log-likelihoods are finite and well defined.
Conditions (c) and (d) ensure that the Kullback--Leibler
divergence uniquely identifies the true density, up to
almost-everywhere equality, and varies continuously over
$\Phi_i$. Conditions (e) and (f) provide uniform moment
and likelihood-ratio bounds that, together with the other
assumptions, are used to derive summable tail bounds
for the proposed test statistics. Finally, Condition (g)
controls the expected local oscillation of the log-density
uniformly over $\Phi_i$, for each fixed true density, and
enables likelihood comparisons within small $\sL^1$
neighborhoods. In Proposition~\ref{prop:expf} of
Appendix~\ref{appf}, we verify these conditions for regular minimal exponential families.

A non-parametric example is obtained by mixing Gaussian densities
with fixed variance and means in a compact interval. Each component
density is determined by its mean, while the mixing distribution
can be any probability measure on the interval. The following
proposition shows that regularity conditions (a)--(g) also hold
for mixtures of densities from the exponential families.

\begin{proposition}\label{mixt}
	For each $i\in[M]$, let
	$\cE_i=\{g_\theta:\theta\in\Theta_i\}$ be an exponential
	family with densities of the form
	\[
	g_\theta(x)
	=
	h_i(x)\exp\{\theta^\top T_i(x)-b_i(\theta)\},
	\qquad \theta\in\Theta_i,
	\]
	where $\Theta_i$ is a compact subset of the natural parameter
	space. Let $\mathcal P(\Theta_i)$ denote the set of Borel
	probability measures on $\Theta_i$, and define the class
	of mixture densities
	\[
	\Phi_i
	:=
	\{q_{i,P}:P\in\mathcal P(\Theta_i)\},
	\qquad
	q_{i,P}(x)
	:=
	\int_{\Theta_i}g_\theta(x)\,P(d\theta).
	\]
	Then $\Phi_i$ satisfies regularity conditions (a)--(g).
\end{proposition}

\begin{IEEEproof}
The proof is presented in Appendix \ref{appf}.
\end{IEEEproof}

\subsection{Generalized Log-Likelihood Ratio}
To define the Generalized Log-Likelihood Ratio (GLLR), we first consider the log-likelihood of source $i$ under a candidate density $g \in \Phi_i$, based on the observations collected up to time $n$. For any sampling rule $R$, define
\begin{equation}\label{llr}
	l_i(n;g) := l_i(n-1;g) + \log\bigl(g(X_i(n))\bigr) R_i(n),
\end{equation}
with $l_i(0;g) := 0$. Here, $n$ denotes the global sampling time, and should not be confused with the total number of observations collected from source $i$ up to time $n$, which is given by $N_i^R(n)$ in \eqref{nrt}. The GLLR statistic is defined as
\begin{equation}\label{gslr}
	\Lambda_i^R(n) := l_i^1(n) - l_i^0(n), \qquad n \in \bN,
\end{equation}
where
\begin{equation}\label{suprem}
	\begin{aligned}
		l_i^1(n) := \sup_{g \in \Phi_i^1} l_i(n;g),\qquad l_i^0(n) := \sup_{g \in \Phi_i^0} l_i(n;g).
	\end{aligned}
\end{equation}

Large positive values of $\Lambda_i^R(n)$ indicate that $i \in A$, while large negative values indicate that $i \notin A$. The following result shows that the GLLR satisfies condition \eqref{stt2} of Theorem~\ref{theo:prob_ao}.

\begin{theorem}\label{tGLLR}
	Let $A \subseteq [M]$, $\boldsymbol{f} \in \boldsymbol{\Phi}_A$, and let $R$ be a sampling rule. If $\{\Lambda_i^R(n),\, i \in [M]\}$ are the GLLR statistics, then for every $\rho \in (0,1]$ and $\epsilon > 0$, the sequences
	\begin{equation}\label{thm51}
		\begin{aligned}
			&\Pro_{\boldsymbol{f}} \left( 
			\frac{1}{n} \Lambda^R_i(n) <
			\rho I_i(f_i) - \epsilon,\,
			\pi_i^R(n) \geq \rho
			\right), \quad \forall \, i \in A, \\
			&\Pro_{\boldsymbol{f}} \left( 
			-\frac{1}{n} \Lambda^R_j(n) <
			\rho J_j(f_j) - \epsilon,\,
			\pi_j^R(n) \geq \rho
			\right), \quad \forall \, j \notin A,
		\end{aligned}
	\end{equation}
	are exponentially decaying.
\end{theorem}

\begin{IEEEproof}
We cover the null class $\Phi_i^0$ by finitely
many neighborhoods and bound the likelihood over each
neighborhood using the corresponding upper density envelope.
This reduces the GLLR deviation event to a finite union of
deviation events for sums of log ratios between the true
density and these envelopes. Lemma~\ref{aux1} provides
exponential bounds for these events, and a union bound
yields the desired exponential decay. The full proof is given in
Appendix~\ref{app0}.
\end{IEEEproof}

We observe that, for each $i \in [M]$, by choosing $\rho = c_i^*(\boldsymbol{f}) - \epsilon > 0$ for sufficiently small $\epsilon > 0$, Theorem~\ref{tGLLR} implies condition \eqref{stt2} of Theorem~\ref{theo:prob_ao}. 

\subsection{Adaptive Log-Likelihood Ratio}
For the definition of the Adaptive Log-Likelihood Ratio (ALLR), we consider, in addition to \eqref{llr}–\eqref{suprem}, a log-likelihood constructed using estimated densities. Specifically, for any sampling rule $R$, define
\begin{equation}
	\hat{l}_i(n) := \hat{l}_i(n-1) + \log\bigl(\hat{f}_{i,n-1}(X_i(n))\bigr) R_i(n),
\end{equation}
with $\hat{l}_i(0) := 0$, where $\hat{f}_{i,n-1}$ is the estimator of $f_i$ based on observations up to time $n-1$. Unlike the GLLR, the ALLR replaces the unknown densities with sequential estimates, adapting gradually to model uncertainty. For each $i \in [M]$, the ALLR statistic is defined as
\begin{equation}\label{allr}
	\Lambda_i^R(n) :=
	\begin{cases}
		\hat{l}_i(n) - l_i^0(n), & \text{if } l_i^0(n) \leq \min\{l_i^1(n), \hat{l}_i(n)\},\\[4pt]
		-\bigl(\hat{l}_i(n) - l_i^1(n)\bigr), & \text{if } l_i^1(n) \leq \min\{l_i^0(n), \hat{l}_i(n)\},\\[4pt]
		0, & \text{otherwise}.
	\end{cases}
\end{equation}

Thus, the ALLR compares the estimated log-likelihood with the best alternatives under the regular and anomalous hypotheses. Large positive values of $\Lambda_i^R(n)$ provide evidence that $i \in A$, while large negative values indicate that $i \notin A$.

The following theorem establishes condition \eqref{stt2}
of Theorem~\ref{theo:prob_ao} for the ALLR statistics.
Since these statistics incorporate sequential density
estimates, their analysis requires control of the
convergence of these estimates. Proposition~\ref{tpro}
shows that conditions \eqref{stt1} and \eqref{stt3} imply
the summability condition \eqref{tec} for the random
times $\tau_i^R(\delta)$ defined in \eqref{tta}.
We therefore state the theorem directly in terms of
this summability condition.

\begin{theorem}\label{thm_allr}
	Let $A \subseteq [M]$, $\boldsymbol{f} \in \boldsymbol{\Phi}_A$, and let $R$ be a sampling rule. Suppose that:
	\begin{enumerate}[(i)]
		\item $\{\Lambda_i^R(n),\, i \in [M]\}$ are the ALLR statistics;
		
		\item for the estimator $\{\hat{\boldsymbol{f}}_n : n \in \bN\}$, for any $\delta > 0$, and $c > 0$,
		\begin{equation}\label{asmp}
			\sum_{n=1}^{\infty}
			\Pro_{\boldsymbol{f}}\!\left(
			\tau^R_i(\delta) > c\, n
			\right) < \infty,
			\quad \forall i \in [M].
		\end{equation}
	\end{enumerate}
	Then, for every $\rho \in (0,1]$ and $\epsilon > 0$,
	\begin{equation}
		\begin{aligned}
			&\sum_{n=1}^{\infty} 
			\Pro_{\boldsymbol{f}} \left(
			\frac{1}{n} \Lambda^R_i(n) <
			\rho I_i(f_i) - \epsilon,\,
			\pi_i^R(n) \ge \rho
			\right) < \infty, \quad \forall\, i \in A, \\
			&\sum_{n=1}^{\infty} 
			\Pro_{\boldsymbol{f}} \left(
			-\frac{1}{n} \Lambda^R_j(n) <
			\rho J_j(f_j) - \epsilon,\,
			\pi_j^R(n) \ge \rho
			\right) < \infty, \quad \forall\, j \notin A.
		\end{aligned}
	\end{equation}
\end{theorem}

\begin{IEEEproof}
Using the definition of the ALLR and the bounds established
in the proof of Theorem~\ref{tGLLR}, we reduce the argument
to controlling the estimation error
$l_i(n;f_i)-\hat l_i(n)$. This error decomposes into a
martingale term and a sum of Kullback--Leibler divergences.
Establishing summable tail bounds for these two terms
yields \eqref{stt2}. The full proof is given in
Appendix~\ref{app0}.
\end{IEEEproof}

By choosing $\rho = c_i^*(\boldsymbol{f}) - \epsilon > 0$ for each $i \in [M]$, for sufficiently small $\epsilon > 0$, Theorem~\ref{thm_allr} implies condition \eqref{stt2} of Theorem~\ref{theo:prob_ao}.

\subsection{Error Control and Threshold Selection}

In addition to condition \eqref{stt2} of Theorem~\ref{theo:prob_ao}, we must verify that, for a suitable test statistic, the stopping and decision rules introduced in Subsection~\ref{st_dec} satisfy the error constraint \eqref{error_constr} and the asymptotic scaling \eqref{ab}.

\subsubsection{GLLR}
For GLLR, this can be established under parametric assumptions on $\Phi_i$, $i \in [M]$. In particular, when each $\Phi_i$ is a subset corresponding to a multidimensional exponential family, Chan and Lai \cite[Theorem~1 and Subsection~3.1]{chan_lai00} provide a first-order asymptotic approximation of the error probabilities for GLLR-type tests.

\begin{proposition}\label{prop:error_contro_gllr}
	Suppose that, for each $i \in [M]$, the sets $\Phi_i^0$ and $\Phi_i^1$ correspond to parametric subsets of a multidimensional exponential family. Let $R$ be an arbitrary sampling rule and $\{\Lambda_i^{R}(n), i \in [M]\}$ be the GLLR statistics. Then the policy $(R, T^{R}, \Delta^{R})$ satisfies the error constraint \eqref{error_constr} asymptotically when the thresholds $a,b \to \infty$ are chosen so that
	\begin{align} \label{thresholds_gi}
		C\, a^{q/2} e^{-a} \sim \beta,
		\qquad
		C\, b^{q/2} e^{-b} \sim \alpha,
	\end{align}
	as $\alpha, \beta \to 0$. In particular, these thresholds satisfy the asymptotic scaling \eqref{ab}.
\end{proposition}

Although the asymptotic optimality theory developed in this paper is formulated under the regularity conditions of Subsection~\ref{g_ass}, Proposition~\ref{prop:expf} shows that these assumptions are satisfied by regular minimal exponential families. Consequently, combining Proposition~\ref{prop:error_contro_gllr} with Theorem \ref{theo:prob_ao} establishes the asymptotic optimality of the proposed GLLR policy for multidimensional exponential-family models.

In Proposition \ref{prop:error_contro_gllr}, $C$ and $q$ are constants. The constant $C$ is given in \cite[Theorem~1]{chan_lai00} and depends on the dimension of the parameter, the exponential-family model, and geometric properties of the parameter space $\boldsymbol{\Phi}$. The constant $q$ is the maximum absolute difference between the dimensions of the manifolds $\Phi_i^1$ and $\Phi_i^0$, i.e.,
\[
q := \max_{i \in [M]} \bigl| \dim(\Phi_i^1) - \dim(\Phi_i^0) \bigr|.
\]
Related asymptotic approximations, with less precise constants, were obtained by Woodroofe \cite[Sections~7.2 and~8.2]{woodroofe1982nonlinear}. In the one-dimensional exponential-family case, the constant $C$ can be evaluated explicitly using Lorden’s approximation; see \cite[(5.67)]{tartakovsky2014sequential}.\\

\subsubsection{ALLR}

The following result shows that the ALLR yields nonasymptotic error control with thresholds that also satisfy the asymptotic scaling \eqref{ab}, without imposing a parametric structure on $\Phi_i$. In contrast to the GLLR, this provides guarantees for nonparametric hypothesis classes.

\begin{proposition}\label{prop:error_contro}
	Let $R$ be an arbitrary sampling rule, and $\{\Lambda_i^{R}(n), i \in [M]\}$ be the ALLR statistics. If the thresholds of the stopping rule are chosen as
	\begin{equation}\label{thresholds}
		\begin{aligned}
			a = |\log(\beta)| + \log M,\quad
			b = |\log(\alpha)| + \log M,
		\end{aligned}
	\end{equation}
	then the policy $(R,T^R,\Delta^R)$ satisfies the error constraint \eqref{error_constr}, and the asymptotic scaling \eqref{ab}.
\end{proposition}

\begin{IEEEproof}
	We prove the bound for the probability of missed detection, as the bound for false alarm is analogous. A missed detection occurs if there exists $i \in A$ such that
	\[
	B_i := \{\Lambda_i^R(T^R) \le -a\}.
	\]
	By the union bound,
	\[
	\Pro_{\boldsymbol f}\big(|A \setminus \Delta^R(T^R)| \ge 1\big)
	\le \sum_{i\in A} \Pro_{\boldsymbol f}(B_i).
	\]
	
	If $\Lambda_i^R(T^R) \ge 0$, then $\Pro_{\boldsymbol f}(B_i)=0$. Otherwise, by \eqref{allr},
	\[
	\Lambda_i^R(T^R) = -\big(\hat l_i(T^R)-l_i^1(T^R)\big),
	\]
	and hence
	\[
	\Pro_{\boldsymbol f}(B_i)
	\le
	\Pro_{\boldsymbol f}\!\left(
	\sup_{n\in\bN}\{\hat l_i(n)-l_i^1(n)\}\ge a
	\right).
	\]
	Since $i\in A$, \eqref{suprem} implies $l_i^1(n)\ge l_i(n;f_i)$ for all $n$, and therefore
	\[
	\Pro_{\boldsymbol f}(B_i)
	\le
	\Pro_{\boldsymbol f}\!\left(
	\sup_{n\in\bN}\{\hat l_i(n)-l_i(n;f_i)\}\ge a
	\right)
	\]
	
	We now show that
	\[
	M_i(n):=\exp\{\hat l_i(n)-l_i(n;f_i)\}, \quad n\ge 0,
	\]
	is a nonnegative martingale. Indeed,
	\[
	\Exp_{\boldsymbol f}\!\left[M_i(n)\mid \cF_{n-1}^R\right]
	=
	M_i(n-1)\,
	\Exp_{\boldsymbol f}\!\left[
	\left(\frac{\hat f_{i,n-1}(X_i(n))}{f_i(X_i(n))}\right)^{R_i(n)}
	\Bigm| \cF_{n-1}^R
	\right].
	\]
	
	Since $R_i(n)$ are $\cF_{n-1}^R$-measurable, we consider two cases. If $R_i(n)=0$, the conditional expectation equals $1$. If $R_i(n)=1$, then using the independence of $X_i(n)$ from $\cF_{n-1}^R$, and the fact that $\hat f_{i,n-1}$ is $\cF_{n-1}^R$-measurable
	\[
	\Exp_{\boldsymbol f}\!\left[
	\frac{\hat f_{i,n-1}(X_i(n))}{f_i(X_i(n))}
	\Bigm| \cF_{n-1}^R
	\right]
	=
	\int \hat f_{i,n-1}(x)\,\dd x
	=1.
	\]
	Thus, $\{M_i(n):n\ge0\}$ is a nonnegative martingale. By Ville's inequality,
	\[
	\Pro_{\boldsymbol f}\!\left(
	\sup_{n\in\bN}\{\hat l_i(n)-l_i(n;f_i)\}\ge a
	\right)
	\le e^{-a}.
	\]
	
	Consequently,
	\[
	\Pro_{\boldsymbol f}\big(|A \setminus \Delta^R(T^R)| \ge 1\big)
	\le |A|e^{-a}\le Me^{-a}.
	\]
	Choosing $a=|\log\beta|+\log M$ yields
	\[
	\Pro_{\boldsymbol f}\big(|A \setminus \Delta^R(T^R)| \ge 1\big)\le \beta.
	\]
\end{IEEEproof}

The above results highlight a distinction between the GLLR and the ALLR at the level of theoretical guarantees. The GLLR can be applied in nonparametric settings; however, rigorous error control guarantees are available only under additional parametric assumptions, typically when $\Phi_i$ corresponds to an exponential family. In contrast, the ALLR admits nonasymptotic error control without requiring a specific parametric structure on $\Phi_i$ thereby ensuring applicability to a broader class of problems. Also the error bound of the ALLR is independent of the structure of space $\Phi_i$, in contrast to \eqref{thresholds_gi} where $C$ depends on $\Phi_i$.

\section{Suitable Estimators}\label{estim}
In this section, we present two density estimators, the
maximum likelihood estimator (MLE) and the kernel density
estimator (KDE), and establish conditions under which
they satisfy \eqref{stt3} of Theorem~\ref{theo:prob_ao}.

\subsection{Maximum Likelihood Estimator}\label{mlee}
The maximum likelihood estimator of $f_i$, based on the
observations collected up to time $n$, is defined by
\begin{equation*}
	\hat f_{i,n}
	\in
	\arg\max_{g\in\Phi_i} l_i(n;g),
\end{equation*}
where $l_i(n;g)$ is the log-likelihood for source $i$
under density $g$, defined in \eqref{llr}. To ensure that a maximizer exists for every realization of the observations, we assume that, for each $i\in[M]$:
\begin{enumerate}[(i)]
	\item
	the class $\Phi_i$ is compact with respect to the
	$\sL^1$ norm;
	
	\item
	for every $x\in\bS$, the map
	\[
	g\mapsto\log g(x)
	\]
	is upper semicontinuous with respect to the
	$\sL^1$ norm.
\end{enumerate}
Under these assumptions, $g\mapsto l_i(n;g)$ is upper
semicontinuous on the compact set $\Phi_i$. Hence, by
the Weierstrass extreme value theorem for upper
semicontinuous functions, it attains its maximum. Lemma \ref{lem:homeomorphism}, implies that the above assumptions are satisfied by every regular minimal exponential family with compact parameter space.

\begin{theorem}\label{mle}
	Let $A\subseteq[M]$,
	$\boldsymbol f\in\boldsymbol{\Phi}_A$,
	and let $R$ be a sampling rule. Suppose that the
	regularity conditions in Subsection~\ref{g_ass}
	and the upper semicontinuity assumption (ii) in
	Subsection~\ref{mlee} hold. For each $i\in[M]$,
	let $\{\hat f_{i,n}:n\in\bN\}$ be a sequence of
	MLEs defined above. Then, for every $i\in[M]$,
	$\epsilon>0$, and $\rho\in(0,1]$, the sequence
	\begin{equation*}
		\left\{
		\Pro_{\boldsymbol f}\left(
		\hat f_{i,n}\notin V_\epsilon(f_i),\,
		\pi_i^R(n)\ge\rho
		\right) \, :\, n \in \bN
		\right\}
	\end{equation*}
	decays exponentially.
\end{theorem}

\begin{IEEEproof}
	We cover the compact set $\Phi_i\setminus V_\epsilon(f_i)$
	by finitely many suitable neighborhoods and bound the
	log-likelihood over each using its upper density envelope.
	This reduces the MLE deviation event to a finite union
	of events whose probabilities decay exponentially.
	The full proof is given in Appendix~\ref{appA}.
\end{IEEEproof}

\subsection{Kernel Density Estimator}
In this subsection, we assume that $\bS=\bR$. A non-parametric
estimate of the density $f_i$ can be obtained via the kernel
density estimator (KDE). Let $K:\bR\to\bR$ be a probability
density, which we denote as the kernel. The kernel density
estimator based on the observed samples is defined as
\[
\hat f_{i,n}(x)
=
\frac{1}{N_i^R(n)\, \eta_{i,n}}
\sum_{k=1}^n
R_i(k)\,
K\!\left(
\frac{x-X_i(k)}{\eta_{i,n}}
\right),
\qquad x\in\bR,
\]
where $\eta_{i,n}:=h_{N_i^R(n)}>0$ is a function of
$N_i^R(n)$, where $h:\bN\to\bR_+$. The sequence
$\{\eta_{i,n}:n\in\bN\}$ is assumed to satisfy
\begin{equation}\label{eqq}
	\eta_{i,n}\to0,
	\quad
	N_i^R(n)\eta_{i,n}\to\infty,
	\quad
	\mbox{as } N_i^R(n)\to\infty.
\end{equation}
We assume that the kernel $K$ satisfies the following conditions:
\begin{enumerate}[(1)]
	\item
	$\int_{\bR}uK(u)\,\dd u=0$,
	
	\item
	$\int_{\bR}u^2K(u)\,\dd u<\infty$,
	
	\item
	$K\in\sL^p$, i.e.,
	$\int_{\bR}(K(u))^p\,\dd u<\infty$,
	for some $p>15/2$.
\end{enumerate}

For the following theorem, we assume that, for the same
$p>15/2$ as in condition (3), the true density
$f_i\in\Phi_i$ satisfies:
\begin{enumerate}[(a)]
	\item
	$f_i\in\sL^{p/2}$, i.e.,
	$\int_{\bR}f_i^{p/2}(x)\,\dd x<\infty$,
	
	\item
	$f_i$ is twice continuously differentiable on $\bR$,
	and $f_i''\in\sL^p$, i.e., $\int_{\bR}|f_i''(x)|^p\,\dd x<\infty$.
\end{enumerate}
The assumptions (a)--(b) are expressed for the true density
$f_i$, which is unknown. Hence, they are assumed to hold
locally on $\Phi_i$, i.e., in a neighborhood of every
$f_i\in\Phi_i$.

\begin{theorem}\label{kde}
Assume that, for each \(i\in[M]\),
\(\{\hat f_{i,n}:n\in\bN\}\) is the sequence of KDEs defined above,
with bandwidth sequence \(\{\eta_{i,n}\}\) satisfying \eqref{eqq}.
Suppose further that the kernel assumptions \((1)\)--\((3)\) and
the density assumptions \((a)\)--\((b)\) hold.  Then for every $\rho \in (0,1]$, $\epsilon > 0$,
\begin{equation*}
	\sum_{n=1}^{\infty} 
	n\, \Pro_{\boldsymbol f}\left(
	\hat f_{i,n}\notin V_\epsilon(f_i),\,
	\pi_i^R(n)\ge\rho
	\right) < \infty, 
	\quad \forall i \in [M].
\end{equation*}
\end{theorem}

\begin{IEEEproof}
The proof follows the arguments in
\cite[Theorem~24.1]{van2000asymptotic}, adapted to the present
sequential setting with observation control. We first reduce the
$\sL^1$ deviation to an $\sL^p$ deviation. We then partition the
event of interest according to the number of samples collected
from source $i$ and decompose
the KDE error into its stochastic and bias components. The bias
is controlled using Taylor's expansion and Minkowski's
inequality, while the stochastic component is bounded using the
martingale version of Rosenthal's inequality. The bandwidth
choice, together with $p>15/2$, makes the resulting weighted
probability bounds summable. The detailed proof is provided in
Appendix~\ref{appA}.
\end{IEEEproof}

\section{Simulation Study - Advantage over agnostic sampling}\label{simul1}
In this section, we present simulation results that illustrate the asymptotic theory developed in the previous sections. Specifically, we compare the performance of the GLLR and the ALLR combined with each of the proposed estimators (MLE, and KDE), under the asymptotically optimal sampling frequencies suggested by our theory. We also compare their performance with that of an agnostic sampling rule that samples each source uniformly with frequency $K/M$.

The simulation setup consists of $M=10$ sources, and at each time instant only half of them can be observed, i.e., $K=5$. We consider both a parametric and a nonparametric framework. In each case, the objective is to identify the subset of anomalous sources as quickly as possible. We assume that the anomalous set consists of the first $|A|$ sources, i.e., $A=\{1,\ldots,|A|\}$, and perform $11$ simulation experiments, corresponding to the possible values $|A|\in\{0,\ldots,10\}$. The thresholds $a$ and $b$ in the stopping rule \eqref{stp0m} are selected via simulation so that both error probabilities in \eqref{error_constr} are equal to $10^{-2}$. For each estimated expected stopping time, the Monte Carlo standard error is equal to $1$. The simulated procedure is summarized in Algorithm \ref{alg1}.

\begin{algorithm}
	\caption{Policy for anomaly identification.}
	\begin{enumerate}[(i)]
		
		\item Repeat the following steps (1)-(5), until the stopping condition \eqref{stp0m} is satisfied.
		
		\begin{enumerate}[(1)]
			\item Based on the already collected measurements, estimate the density $\boldsymbol{f}$ by $\hat{\boldsymbol{f}}_{n}$, using a suitable estimator.
			
			\item Compute the desired minimum long-run sampling frequencies $(c^{*}_{1}(\hat{\boldsymbol{f}}_{n}),\ldots,c^{*}_{M}(\hat{\boldsymbol{f}}_{n}))$ according to \eqref{c_star}.
			
			\item Choose sampling probabilities $c^R_{i}(\hat{\boldsymbol{f}}_{n})$ such that $c^R_{i}(\hat{\boldsymbol{f}}_{n}) \geq c^{*}_{i}(\hat{\boldsymbol{f}}_{n})$, and $\sum_{i=1}^{M} c^R_{i}(\hat{\boldsymbol{f}}_{n}) \leq K$.
			
			\item Perform sampling based on a probabilistic sampling rule presented in Section \ref{prob_rules}.
			
			\item Update the LLRs of the sources we have received samples.
			
		\end{enumerate}
		
		\item Decide as anomalous the sources with positive LLRs.
	\end{enumerate}
	\label{alg1}
\end{algorithm} 

\subsection{Parametric framework}
In the parametric framework, the measurements from each source $i \in [M]$ follow a Gaussian density,
\begin{equation*}
	X_{i}(n) \sim \cN(\mu_i, \sigma^{2}_{i}), \qquad n \in \bN.
\end{equation*}
If $\mu_i \in [0,0.5]$ and $\sigma_i \in [0,0.2]$ then source $i$ is characterized as regular, whereas if $\mu_i \in [0.55,1]$ and $\sigma_i \in [0.25,0.4]$ the source is characterized as anomalous. The true but unknown values of $(\mu_i,\sigma_{i})$ for each $i \in [M]$, for this simulation, are
\begin{equation*}
	(\mu_i, \sigma_i) =
	\begin{cases}
		(0.7+0.02i, \ 0.28+0.01i), & i \in A, \\
		(0.2+0.02i, \ 0.08+0.005i), & i \notin A.
	\end{cases}
\end{equation*}
For each source $i\in[M]$, the parameters are estimated by their empirical counterparts based on the observations collected from that source. Specifically,
\begin{equation}\label{mle_c}
\hat{\mu}_{i}(n) := \frac{1}{N_i^{R}(n)}\sum_{m=1}^{n} X_{i}(m)R_i(m), \quad \hat{\sigma}_{i}(n):= \sqrt{\frac{1}{N_i^{R}(n)} \sum_{m=1}^{n} (X_{i}(m) - \hat{\mu}_{i}(n))^2 R_i(m)}. 
\end{equation}

When $i \in A$, the difficulty of the underlying testing problem corresponding to source $i$ is determined by the KL divergence of $\cN(0.7+0.02i,0.28+0.01i)$ versus $\cN(0.5,0.2)$, where the latter is the regular density with the smallest KL distance from $\cN(0.7+0.02i,0.28+0.01i)$, and KL divergence \eqref{inf1} is
\begin{equation*}
I_i(\mu_{i},\sigma_{i}) = \dD(\cN(0.7+0.02i,0.28+0.01i)||\cN(0.5,0.2)).
\end{equation*}   
When $i \notin A$, the difficulty of the underlying testing problem corresponding to source $i$ is determined by the KL divergence of $\cN(0.2+0.02i, 0.08+0.005i)$ versus $\cN(0.55,0.25)$, where the latter is the anomalous density with the smallest KL distance from $\cN(0.2+0.02i, 0.08+0.005i)$, and the KL divergence \eqref{inf2} is
\begin{equation*}
J_i(\mu_{i},\sigma_{i}) = \dD(\cN(0.2+0.02i, 0.08+0.005i)||\cN(0.55,0.25)).
\end{equation*}   
Table \ref{ttb} summarizes the values of $I_i$, $J_i$ for each $i \in [M]$.
\begin{table}[h!]
	\begin{center}
		\begin{tabular}{ |c||c|c|c|c|c|c|c|c|c|c| } 
			\hline
			$i $ & 1 & 2 & 3 & 4 & 5 & 6 & 7 & 8 & 9 & 10\\ 
			\hline
			\hline
			$I_i(\mu_{i},\sigma_{i})$ & 0.78 & 0.93 & 1.10 & 1.29 & 1.48  & 1.69 & 1.91 & 2.15 & 2.40 & 2.66\\  
			\hline
			$J_i(\mu_{i},\sigma_{i})$ & 1.50 & 1.35 & 1.21 & 1.07 & 0.95 & 0.84 & 0.73 & 0.63 & 0.54 & 0.46\\ 
			\hline
		\end{tabular}
	\end{center}
	\caption{KL divergences.}
	\label{ttb}
\end{table}

When $A=\{1,\ldots,k\}$, with $k\in{0,\ldots,M}$, the difficulty of the anomaly identification problem is determined by the collection of KL divergences
\begin{equation*}
S_k := \{I_1(\mu_{1},\sigma_{1}),\ldots,I_k(\mu_{k},\sigma_{k}),J_{k+1}(\mu_{k+1},\sigma_{k+1}),\ldots,J_M(\mu_{M},\sigma_{M})\},\quad 0<k<M,
\end{equation*}
and $S_0:=\{J_i(\mu_{i},\sigma_{i}) \,:\, i \in [M]\}$, and $S_M:=\{I_i(\mu_{i},\sigma_{i}) \, :\, i \in [M]\}$. Also, 
$$S_k := S_{k-1} \cup \{I_k(\mu_{k},\sigma_{k})\} \setminus \{J_k(\mu_{k},\sigma_{k})\},$$ and thus as $k$ increases from $0$ to $3$ a larger KL divergence is replaced by a smaller one, whereas from $4$ to $10$ a smaller KL divergence is replaced by a larger one. 

\begin{figure}[htb!]
	\centering
	\centerline{\includegraphics[width=9 cm]{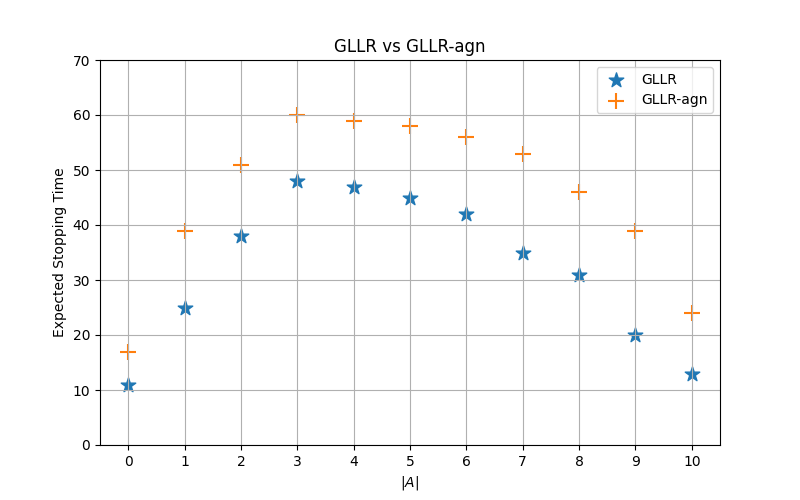}}
	\caption{Expected stopping time vs number of anomalies $|A|$, GLLR vs GLLR agnostic}
	\label{fi1}
\end{figure}

\begin{figure}[htb!]
	\centering
	\includegraphics[width=0.45\textwidth]{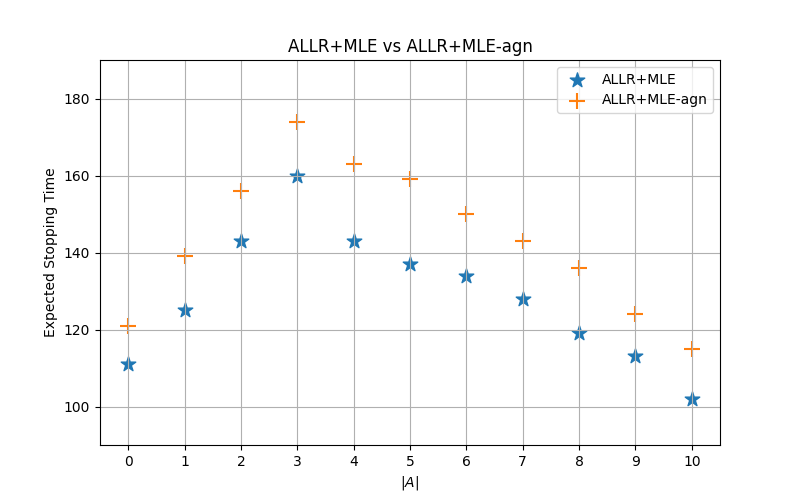}
	\includegraphics[width=0.45\textwidth]{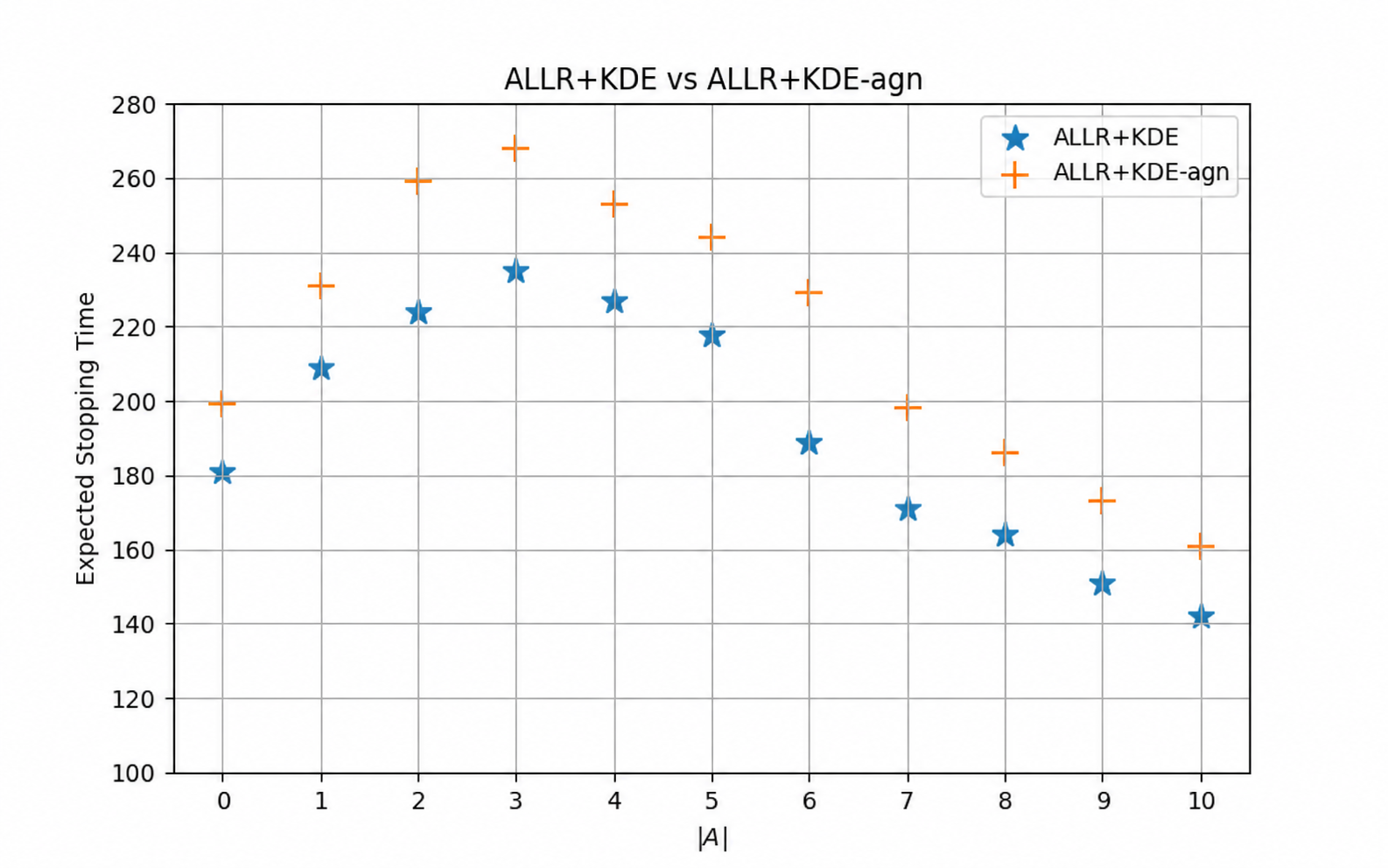}
	\caption{Expected stopping time versus the number of anomalies $|A|$ for ALLR+MLE, and ALLR+KDE vs their agnostic counterparts.}
	\label{fi2}
\end{figure}

Since the expected stopping time is approximately inversely proportional to the KL divergences that characterize the difficulty of the underlying testing problem, we expect the overall expected stopping time $E[T]$ to increase as $k$ increases from $0$ to $3$, and to decrease as $k$ increases from $4$ to $10$. This behavior follows from the evolution of the sets $S_k$, where initially larger KL divergences are replaced by smaller ones, and subsequently smaller KL divergences are replaced by larger ones. The predicted behavior is observed in Figures~\ref{fi1} and~\ref{fi2}.

Figure~\ref{fi1} compares the performance of the GLLR combined with the asymptotically optimal sampling frequencies prescribed by condition~\eqref{stt1} against the corresponding agnostic sampling rule that samples each source with probability $K/M$. Similarly, Figure~\ref{fi2} compares the performance of the ALLR combined with the MLE, and KDE under the asymptotically optimal sampling frequencies against their agnostic counterparts. For the MLE, the parameter estimates are computed according to \eqref{mle_c}. For the KDE, we employ a Gaussian kernel with bandwidth $h_{i,n}=(N_i^R(n))^{-0.2}$.

In all cases, the expected stopping time attains its maximum around $|A|=3$, corresponding to the most difficult anomaly configuration. Moreover, the sampling rules based on the asymptotically optimal frequencies consistently outperform their agnostic counterparts. Among the adaptive procedures, ALLR+MLE achieves the smallest expected stopping times, followed by ALLR+KDE, reflecting the superior accuracy of the MLE in the Gaussian setting.

In the long run, for each $i \in [M]$, the estimators
\begin{equation*}
\hat{\boldsymbol{\mu}}(n):=(\hat{\mu}_{1}(n),\ldots,\hat{\mu}_{M}(n)), \quad \hat{\boldsymbol{\sigma}}(n):=(\hat{\sigma}_{1}(n),\ldots,\hat{\sigma}_{M}(n)),
\end{equation*}
converge to the actual parameters 
\begin{equation*}
\boldsymbol{\mu}:=(\mu_{1},\ldots,\mu_{M}), \quad \boldsymbol{\sigma}:=(\sigma_{1},\ldots,\sigma_{M}),
\end{equation*}
and the minimum sampling frequencies that guarantee asymptotic optimality are given by  
\begin{equation*}
	c_i^*(\boldsymbol{\mu},\boldsymbol{\sigma}) := 
	\begin{cases}
		x_{A}(\boldsymbol{\mu},\boldsymbol{\sigma}) \; I_A^*(\boldsymbol{\mu},\boldsymbol{\sigma})/ I_i(\mu_{i},\sigma_{i}),  \quad  \text{if}  \; \;  i \in A,\\
		y_{A}(\boldsymbol{\mu},\boldsymbol{\sigma}) \; J_A^*(\boldsymbol{\mu},\boldsymbol{\sigma})/J_i(\mu_{i},\sigma_{i}),  \quad  \text{if}  \; \;  i \notin A,
	\end{cases}
\end{equation*}
where $x_{A}(\boldsymbol{\mu},\boldsymbol{\sigma})$, $y_{A}(\boldsymbol{\mu},\boldsymbol{\sigma})$ are given in Theorem \ref{theo:lowerbound}. We assume $r=1$, and according to Table \ref{ttb} the minimum values of the respective KL divergences are $I_A^*(\boldsymbol{\mu},\boldsymbol{\sigma})=0.78$, $J_A^*(\boldsymbol{\mu},\boldsymbol{\sigma})=0.46$.
\begin{table}[h!]
	\centering
	\begin{tabular}{|c||c|c|c|c|c|c|c|c|c|c|}
		\hline
		$i$ & 1 & 2 & 3 & 4 & 5 & 6 & 7 & 8 & 9 & 10\\
		\hline
		$c_i^*(\boldsymbol{\mu},\boldsymbol{\sigma})$ &
		0.497 &
		0.417 &
		0.353 &
		0.301 &
		0.262 &
		0.462 &
		0.531 &
		0.616 &
		0.718 &
		0.843 \\
		\hline
	\end{tabular}
	\caption{Optimal sampling frequencies $c_i^*(\boldsymbol{\mu},\boldsymbol{\sigma})$ for $|A|=5$, $A=\{1,2,3,4,5\}$.}
	\label{tab:cstar5}
\end{table}

Table~\ref{tab:cstar5} reports the asymptotically optimal sampling frequencies for the case $|A|=5$, i.e., $A={1,2,3,4,5}$. We observe that sources associated with smaller KL divergences are assigned larger sampling frequencies. This is consistent with the asymptotic theory, according to which more samples should be allocated to the most difficult sources.

\subsection{Non-Parametric framework}
In the nonparametric framework, the observations from source $i \in [M]$ are generated according to an unknown density $f_i$,
\begin{equation*}
	X_{i}(n) \sim f_i, \quad n \in \mathbb{N}.
\end{equation*}
We consider two classes of densities, which define the two hypotheses,
\begin{enumerate}[(i)]
	\item the class $\Phi_i^0:=\mathcal S_{\uparrow}$ contains all densities on $[0,1]$ that are nondecreasing, uniformly bounded by $2$, uniformly equicontinuous, and satisfy
	\[
	g(1)-g(0)\ge 0.01,
	\qquad g\in\mathcal S_{\uparrow},
	\]
	
	\item the class $\Phi_i^1:=\mathcal S_{\downarrow}$ contains all densities on $[0,1]$ that are nonincreasing, uniformly bounded by $2$, uniformly equicontinuous, and satisfy
	\[
	g(0)-g(1)\ge 0.01,
	\qquad g\in\mathcal S_{\downarrow}.
	\]
\end{enumerate}
If $f_i \in S_{\uparrow}$ then source $i$ is characterized as regular, whereas if $f_i \in \mathcal{S}_\downarrow$ the source is characterized as anomalous. The above classes satisfy the assumptions imposed in Section~\ref{estim}. The monotonicity constraints ensure separation of the two hypotheses. The true but unknown densities $f_i$, used for this simulation, are
\begin{equation*}
	f_i(x) =
	\begin{cases}
		-0.8x + 1.4, & i \in A, \\
		1.6x + 0.2, & i \notin A.
	\end{cases}
\end{equation*}

\begin{figure}[htb!]
	\centering
	\centerline{\includegraphics[width=9 cm]{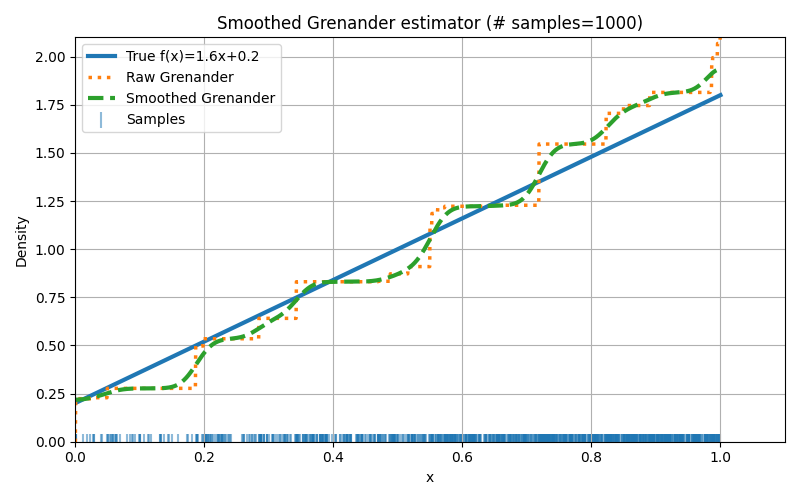}}
	\caption{Smoothed Grenander estimator for the density $f(x)=1.6x+0.2$.}
	\label{fig1}
\end{figure}

For each $i\in[M]$, we estimate $f_i$ using the Grenander estimator, which is the nonparametric maximum likelihood estimator for monotone densities \cite{gren56}. The Grenander estimator is obtained by computing the least concave majorant (respectively, the greatest convex minorant) of the empirical distribution function and then taking its left derivative. Since the resulting estimator is piecewise constant, we employ a smoothed version of the Grenander estimator in order to reduce discontinuities and improve numerical stability. Figure~\ref{fig1} illustrates a typical realization of the smoothed Grenander estimator for the density $f(x)=1.6x+0.2$.

When $i \in A$, the difficulty of the underlying testing problem corresponding to source $i$ is determined by the KL divergence of $f_i(x)=-0.8x+1.4$ versus $g^1_i(x)=0.01x+0.995$, where the latter is the regular density with the smallest KL distance from $f_i(x)$, and the KL divergence \eqref{inf1} is equal to
\begin{equation*}
I_i(f_i) = \dD(f_i||g^1_i)=0.027.
\end{equation*}

When $i \notin A$, the difficulty of the underlying testing problem corresponding to source $i$ is determined by the KL divergence of $f_i(x)=1.6x+0.2$ versus $g^0_i(x)=-0.01x+1.005$, where the latter is the anomalous density with the smallest KL distance from $f_i(x)$, and the KL divergence \eqref{inf2} is equal to
\begin{equation*}
J_i(f_i) = \dD(f_i||g^0_i)=0.116.
\end{equation*}

In case $A=\{1,\ldots,k\}$, $k \in \{0,\ldots,M\}$ the overall difficulty of the anomaly detection problem is determined by the KL divergences
\begin{equation*}
S_k := \{I_1(f_1),\ldots,I_k(f_k),J_{k+1}(f_{k+1}),\ldots,J_M(f_M)\},
\end{equation*}

\begin{figure}[htb!]
	\centering
	\centerline{\includegraphics[width=9 cm]{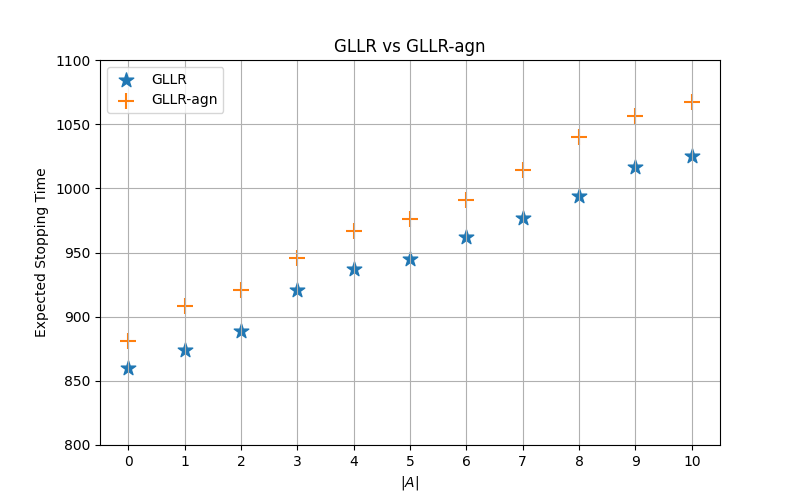}}
	\caption{Expected stopping time vs number of anomalies $|A|$, GLLR vs GLLR agnostic}
	\label{fi3}
\end{figure}

\begin{figure}[htb!]
	\centering
	\includegraphics[width=0.45\textwidth]{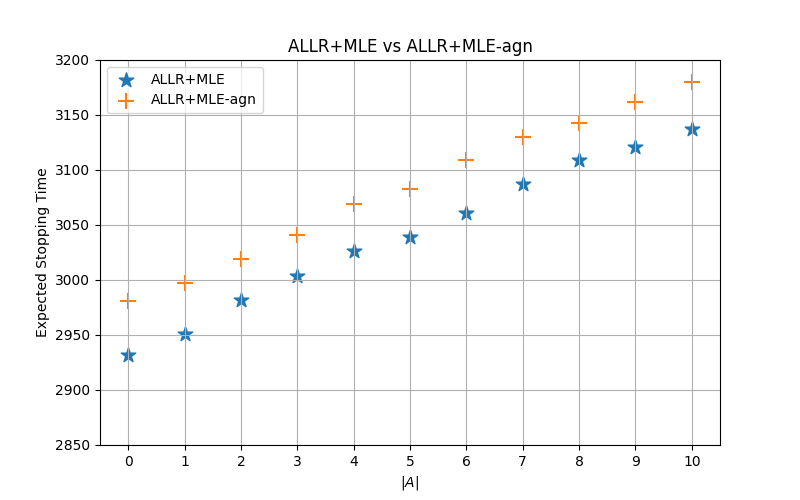}
	\includegraphics[width=0.45\textwidth]{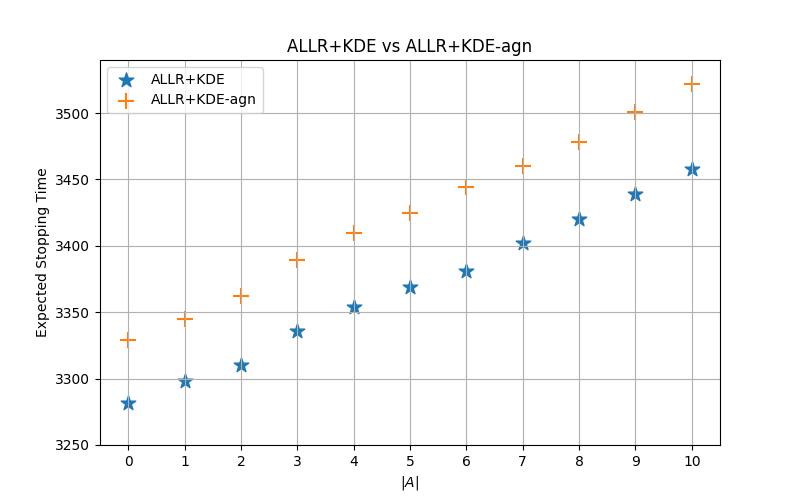}
	\caption{Expected stopping time versus the number of anomalies $|A|$ for ALLR+MLE,  ALLR+KDE vs their agnostic counterparts.}
	\label{fi4}
\end{figure}

Since $I_i(f_i) < J_i(f_i)$, an anomalous source is harder to classify than a regular source. Consequently, as the number of anomalies $k$ increases, the overall difficulty of the anomaly identification problem increases, and we expect the expected stopping time $E[T]$ to increase monotonically. This behavior is clearly observed in Figures~\ref{fi3} and~\ref{fi4}.

Figure~\ref{fi3} compares the performance of the GLLR combined with the asymptotically optimal sampling frequencies prescribed by condition~\eqref{stt1} against the corresponding agnostic sampling rule that samples each source with probability $K/M$. To compute the GLLR, we evaluate the likelihood ratio at the Grenander MLE. Similarly, Figure~\ref{fi4} compares the performance of ALLR combined with the MLE, and KDE under the asymptotically optimal sampling frequencies against their agnostic counterparts. For the ALLR+MLE procedure, the density is estimated by the smoothed Grenander estimator. For the KDE, we employ a Gaussian kernel with bandwidth $h_{i,n}=(N_i^R(n))^{-0.2}$ and we restrict and normalize the estimate at $[0,1]$.

In all cases, the expected stopping time increases with the number of anomalies, reflecting the fact that anomalous sources are associated with smaller KL divergences. Moreover, the sampling rules based on the asymptotically optimal frequencies consistently outperform their agnostic counterparts. Among the adaptive procedures, ALLR+MLE achieves the smallest expected stopping times, followed by ALLR+KDE. This ordering is expected, since the Grenander estimator is the nonparametric MLE for monotone densities and is therefore well matched to the underlying model.

In the long run, for each $i \in [M]$, the estimator $\hat{\boldsymbol{f}}_{n}$ converges to the actual densities $\boldsymbol{f}$. We note that, in the present simulation setup, all anomalous (resp. regular) sources share the same KL divergence. Consequently, $I_A^*(\boldsymbol f)=I$, $J_A^*(\boldsymbol f)=J$, and
\[
\hat K_A(\boldsymbol f)=|A|,
\qquad
\check K_A(\boldsymbol f)=M-|A|,
\qquad
s_A^*(\boldsymbol f)=I/J.
\]
For $K=5$, $M=10$, and assuming $r=1$, the minimum sampling frequencies that guarantee asymptotic optimality are given by
\begin{enumerate}[(i)]
	\item if $|A|=0$
	\begin{equation*}
		c_i^*(\boldsymbol{f})= 0, \quad \forall \; i\in A, \qquad  c_j^*(\boldsymbol{f})=\frac{K}{(M-|A|)}\wedge 1, \quad \forall \; j \notin \,A,
	\end{equation*}
	\item if $0<|A|<M$
	\begin{equation*}
		c_i^*(\boldsymbol{f})=\frac{KJ}{|A|J+(M-|A|)I}\wedge \frac{J}{I}\wedge 1,\;\; \forall \; i \in A, \qquad c_j^*(\boldsymbol{f})=\frac{KI}{|A|J+(M-|A|)I}\wedge \frac{I}{J}\wedge 1,\;\; \forall \; j\notin A,
	\end{equation*}
	\item if $|A|=M$
	\begin{equation*}
		c_i^*(\boldsymbol{f})= \frac{K}{|A|} \wedge 1, \quad \forall \; i\in A, \qquad  c_j^*(\boldsymbol{f})=0, \quad \forall \; j \notin \,A.
	\end{equation*}
\end{enumerate}

\section{Simulation Study - Robustness in model misspecification}\label{simul2}
We compare GLLR with ALLR combined with KDE under model misspecification and sampling constraints. When the parametric model is correctly specified, GLLR can benefit from the efficiency of the MLE. Under misspecification, however, the fitted parametric density may remain systematically different from the true density, whereas KDE can approximate it more accurately. This can improve the estimation of KL numbers and lead to a more efficient allocation of samples across sources.

We fix constants $B,C\geq 1$ and $0<c< 1$, we denote by
\begin{equation*}
	\cH
	:=
	\left\{
	\cN(\mu,\sigma^2):
	-B\leq\mu\leq B,\;
	c\leq\sigma\leq C
	\right\}
\end{equation*}
the set of Gaussian densities with means in $[-B,B]$ and standard deviations in $[c,C]$, and by $\overline{\cH}$ the closed convex hull of $\cH$ in $\sL^1$. Since $\cH$ is compact in $\sL^1$, its closed convex hull $\overline{\cH}$ is also compact in $\sL^1$. For a fixed $0<\epsilon\leq B$, we consider the classes
\begin{equation}\label{eq:robust-alternative-class}
	\cH_1
	:=
	\left\{
	f\in\overline{\cH}:
	\epsilon\leq\Exp_f[X]\leq B,\;
	c^2\leq\operatorname{Var}_f(X)\leq C^2
	\right\}
\end{equation}
and
\begin{equation}\label{eq:robust-null-class}
	\cH_0
	:=
	\left\{
	\cN(\mu,\sigma^2):
	-B\leq\mu\leq0,\;
	c\leq\sigma\leq C
	\right\}.
\end{equation}
 Both $\cH_1$, $\cH_0$ are closed subsets of the compact set $\overline{\cH}$ and thus themselves compact.

We consider two sources, i.e., $M=2$, and we impose the following composite hypotheses testing problem on each source,
\begin{equation}\label{eq:robust-source-hypotheses}
	\begin{aligned}
		H_i^1 &: f_i\in\cH_1,\\
		H_i^0 &: f_i\in\cH_0,
	\end{aligned}
\end{equation}
where $i\in[2]$. We can sample only one source at a time, i.e., $K=1$. The underlying true densities are
\begin{equation}\label{eq:true-mixture-densities}
	f_i
	:=
	\frac12\cN(m_i-d,s^2)
	+
	\frac12\cN(m_i+d,s^2),
	\qquad i\in[2],
\end{equation}
where $s \in (0,1)$ and $d >0$ satisfy
\begin{equation*}
	c\leq s,
	\qquad
	s^2+d^2+1\leq C^2,
	\qquad
	1+d\leq B.
\end{equation*}
The source means are
\begin{equation}\label{eq:true-source-means}
	m_i
	:=
	\begin{cases}
		1, & i=1,\\
		\epsilon, & i=2.
	\end{cases}
\end{equation}
where $\epsilon \in (0,1)$. Each density in \eqref{eq:true-mixture-densities} has mean and variance
\begin{equation*}
	\Exp_{f_i}[X]=m_i,
	\qquad
	\operatorname{Var}_{f_i}(X)=s^2+d^2.
\end{equation*}
The parameter restrictions ensure that both Gaussian components
belong to $\cH$ and that the resulting densities satisfy
$f_1,f_2\in\cH_1$. Thus, both sources are anomalous.

\subsection{Misspecified parametric model.}
Due to the complexity of the alternative class $\cH_1$,
a practitioner may approximate it by the Gaussian subclass
$\cH_1^{\mathrm{s}}$, while retaining the original null class.
Thus, the original and surrogate problems have the same
Gaussian null, and the model misspecification concerns
only the alternative. Specifically, we define
\begin{equation}\label{eq:surrogate-gaussian-classes}
	\begin{aligned}
		\cH_1^{\mathrm{s}}
		&:=
		\left\{
		\cN(\mu,\sigma^2):
		\epsilon\leq\mu\leq B,\;
		c\leq\sigma\leq C
		\right\},
	\end{aligned}
\end{equation}
and the corresponding surrogate hypotheses are
\begin{equation}\label{eq:surrogate-source-hypotheses}
	\begin{aligned}
		H_i^{1,\mathrm{s}} &: f_i\in\cH_1^{\mathrm{s}},\\
		H_i^{0,\mathrm{s}} &: f_i\in\cH_0.
	\end{aligned}
\end{equation}
where $i \in [2]$. Set $\cH_1^{\mathrm{s}}$ is compact in $\sL^1$. The observations are generated from the true mixture densities in \eqref{eq:true-mixture-densities}. Their means are given by \eqref{eq:true-source-means}, and their common variance is $s^2+d^2$.

We compare the performance of GLLR with that of ALLR combined
with KDE, both operating under the simplified Gaussian hypotheses.
GLLR uses the MLE to exploit the assumed Gaussian structure,
which can yield efficient parameter estimation when the model
is correctly specified. ALLR combined with KDE estimates the
underlying density without imposing a Gaussian form. Under
model misspecification, this flexibility can improve the
estimation of the true information numbers used to determine
the optimal sampling frequencies and thereby yield a smaller
first-order asymptotic approximation to the expected stopping
time.

\subsection{Computation of optimal sampling frequencies}

\subsubsection{GLLR for the simplified problem}

GLLR uses the MLE to estimate the Gaussian mean and variance from the observations sampled from each source. The Gaussian MLEs are given by the empirical mean and standard deviation,
\begin{equation*}
	\begin{aligned}
		\widehat{\mu}_i(n)
		&:=
		\frac{1}{N_i^R(n)}
		\sum_{m=1}^{n}X_i(m)R_i(m),\\
		\widehat{\sigma}_i(n)
		&:=
		\left[
		\frac{1}{N_i^R(n)}
		\sum_{m=1}^{n}
		\bigl(X_i(m)-\widehat{\mu}_i(n)\bigr)^2R_i(m)
		\right]^{1/2}.
	\end{aligned}
\end{equation*}
Almost surely, as $n \to \infty$, the empirical estimators converge to the true values
\begin{equation}\label{eq:gaussian-estimator-limits}
	\begin{aligned}
		\bigl(\widehat{\mu}_i(n),\widehat{\sigma}_i(n)\bigr)
		&\longrightarrow
		\left(m_i,\sqrt{s^2+d^2}\right)
	\end{aligned}
\end{equation}
where $i \in [2]$. Therefore, for any Gaussian density $\cN(\mu,\sigma^2)$, the KL divergence from the fitted Gaussian $\cN(m_i,s^2+d^2)$ is
\begin{equation*}
	\dD\!\left(
	\cN(m_i,s^2+d^2)
	\,\middle\|\,
	\cN(\mu,\sigma^2)
	\right)
	=
	\frac12
	\left[
	\log\frac{\sigma^2}{s^2+d^2}
	+
	\frac{s^2+d^2+(m_i-\mu)^2}{\sigma^2}
	-1
	\right],
\end{equation*}
The minimizing mean over $[-B,0]$ is $\mu=0$, and the minimizing variance is  $s^2+d^2+m_i^2$, and thus
\begin{equation*}
	\operatorname*{arg\,min}_{g\in\cH_0}
	\dD\!\left(\cN(m_i,s^2+d^2)\middle\|g\right)
	=
	\cN(0,s^2+d^2+m_i^2),
	\qquad i\in[2].
\end{equation*}
The corresponding information numbers, under the assumption of Gaussian model, are therefore
\begin{equation}\label{ig}
	\begin{aligned}
		I_1^{\mathrm G}
		&:=
		\inf_{g\in\cH_0}
		\dD\!\left(\cN(1,s^2+d^2)\middle\|g\right)
		=
		\frac12\log\left(1+\frac{1}{s^2+d^2}\right),\\
		I_2^{\mathrm G}(\epsilon)
		&:=
		\inf_{g\in\cH_0}
		\dD\!\left(\cN(\epsilon,s^2+d^2)\middle\|g\right)
		=
		\frac12\log\left(1+\frac{\epsilon^2}{s^2+d^2}\right).
	\end{aligned}
\end{equation}
Consequently, the sampling frequencies that satisfy
\begin{equation*}
	c_1^{\mathrm G}(\epsilon)I_1^{\mathrm G}
	=
	c_2^{\mathrm G}(\epsilon)I_2^{\mathrm G}(\epsilon),
	\qquad
	c_1^{\mathrm G}(\epsilon)+c_2^{\mathrm G}(\epsilon)=1,
\end{equation*}
are 
\begin{equation}\label{eq:gaussian-sampling-frequencies}
		c_1^{\mathrm G}(\epsilon)=
		\frac{I_2^{\mathrm G}(\epsilon)}
		{I_1^{\mathrm G}+I_2^{\mathrm G}(\epsilon)}, \qquad
		c_2^{\mathrm G}(\epsilon)=
		\frac{I_1^{\mathrm G}}
		{I_1^{\mathrm G}+I_2^{\mathrm G}(\epsilon)}.
\end{equation}
We note that as $\epsilon \to 0$, the sampling frequencies $c_1^{\mathrm G}(\epsilon) \to 0$, and $c_2^{\mathrm G}(\epsilon) \to 1$.

\subsubsection{ALLR+KDE for the simplified problem}

Although the MLE incorporates the assumed parametric structure when
estimating the unknown parameters, the KDE estimates the underlying density without imposing a Gaussian form. Under the KDE consistency conditions, and since each
source is sampled infinitely often, we have
\begin{equation}\label{eq:kde-consistency-mixtures}
	\left\|\hat f_{i,n}-f_i\right\|_1
	\rightarrow 0,\; \mbox{ a.s.}, \qquad i \in [2],
\end{equation}
where each $f_i$ is defined in \eqref{eq:true-mixture-densities}. For each source, the Gaussian null density minimizing the KL divergence from the true density is
\begin{equation*}
	\operatorname*{arg\,min}_{g\in\cH_0}
	\dD(f_i\|g)
	=
	\cN(0,\nu_i^2),
	\qquad i\in[2],
\end{equation*}
where 
\begin{equation}\label{vi}
\nu_i^2:=m_i^2+s^2+d^2, \quad i\in[2].
\end{equation}
To show this, let us consider a  $g=\cN(\mu,\sigma^2)$, where $\mu\in[-B,0]$ and $\sigma\in[c,C]$. We have
\begin{equation*}
	\begin{aligned}
		D(f_i,g)&=
		\int_{\mathbb R}f_i(x)\log f_i(x)\,dx -\int_{\mathbb R}f_i(x)\log g(x)\,dx\\
		&=\int_{\mathbb R}f_i(x)\log f_i(x)\,dx +\frac12\log(2\pi\sigma^2) + \frac{\Exp_{f_i}[(X-\mu)^2]}{2\sigma^2}.
	\end{aligned}
\end{equation*}
Since 
\begin{equation*}
	\Exp_{f_i}[X]=m_i,
	\qquad
	\operatorname{Var}_{f_i}(X)=s^2+d^2,
\end{equation*}
we obtain
\begin{equation*}
	\begin{aligned}
		\Exp_{f_i}\!\left[(X-\mu)^2\right]
		&=
		\Exp_{f_i}\!\left[(X-m_i+m_i-\mu)^2\right]\\
		&= s^2+d^2+(m_i-\mu)^2,
	\end{aligned}
\end{equation*}
which gives
\begin{equation*}
	D(f_i,g)
	=
	\int_{\mathbb R}f_i(x)\log f_i(x)\,dx
	+ \frac12\log(2\pi\sigma^2) + \frac{s^2+d^2+(m_i-\mu)^2}{2\sigma^2}.
\end{equation*}
The expression is minimized over $\mu\in[-B,0]$ at $\mu=0$, and then the minimizing variance is \eqref{vi}. In order to compute the corresponding information numbers, we set 
\begin{equation*}
g_i := \cN(0,\nu_i^2), \quad f_0:= \frac12\cN(-d,s^2) + \frac12\cN(d,s^2),
\end{equation*}
and we note that $f_i(x)=f_0(x-m_i)$, $i \in [2]$. Thus,
\begin{equation*}
	\begin{aligned}
		I_i(\epsilon):= D(f_i,g)
		&=\int_{\mathbb R}
		f_0(x-m_i)
		\log\frac{f_0(x-m_i)}{g_i(x)}\,dx\\
		&=
		\int_{\mathbb R}
		f_0(y)
		\log\frac{f_0(y)}{g_i(y+m_i)}\,dy\\
		&=\int_{\mathbb R}
		f_0(y)\log\frac{f_0(y)}{g_i(y)} dy + \int_{\mathbb R} f_0(y)\log\frac{g_i(y)}{g_i(y+m_i)} dy.
	\end{aligned}
\end{equation*}
For the centered Gaussian density $g_i$,
\begin{equation*}
	\log\frac{g_i(y)}{g_i(y+m_i)}
	=
	\frac{m_i y}{\nu_i^2}
	+
	\frac{m_i^2}{2\nu_i^2},
\end{equation*}
and since $\Exp_{f_0}[Y]=0$, we get
\begin{equation}\label{ie}
I_i(\epsilon)=D(f_0\|g_i)+\frac{m_i^2}{2\nu_i^2}, \quad i \in [2].
\end{equation}
Consequently, the sampling frequencies that satisfy
\begin{equation*}
	c_1(\epsilon)I_1(\epsilon)
	=
	c_2(\epsilon)I_2(\epsilon),
	\qquad
	c_1(\epsilon)+c_2(\epsilon)=1,
\end{equation*} 
are
\begin{equation}\label{eq:allr-sampling-frequencies}
	c_1(\epsilon)
	=
	\frac{I_2(\epsilon)}
	{I_1(\epsilon) +I_2(\epsilon) },
	\qquad
	c_2(\epsilon)
	=
	\frac{I_1(\epsilon)}
	{I_1(\epsilon)+I_2(\epsilon)}.
\end{equation}
In contrast to \eqref{eq:gaussian-sampling-frequencies}, both
\begin{equation}\label{pss}
 \lim_{\epsilon \to 0^+} c_1^{\mathrm A}(\epsilon),\;\; \lim_{\epsilon \to 0^+} c_2^{\mathrm A}(\epsilon) > 0,
\end{equation}
due to the fact that $\lim\limits_{\epsilon \to 0} I_i(\epsilon) \geq D(f_0\|\cN(0,s^2+d^2)) > 0$, for both $i \in [2]$.

\subsection{Effect on the expected stopping time}
We examine how replacing the original problem
\eqref{eq:robust-source-hypotheses} by the Gaussian surrogate
\eqref{eq:surrogate-source-hypotheses} affects the expected stopping time asymptotically. We denote by $T_G^R$ the stopping time obtained using GLLR under the simplified Gaussian model, and by $T_{AK}^R$
the stopping time obtained using ALLR combined with KDE. Both procedures use stopping rule \eqref{stp0m}, with their
respective statistics. The quantities $I_i(\epsilon)$ govern the limiting growth of the ALLR statistic, whereas $I_i^{\mathrm G}(\epsilon)$ govern the limiting growth of the GLLR statistic under the true densities. Although $I_i^{\mathrm G}(\epsilon)$ are defined in \eqref{ig} assuming that the fitted Gaussian is the true density, in this example they also equal the limiting GLLR drifts under the true densities because $f_i$ and $\cN(m_i,s^2+d^2)$ share the same mean and variance. For fixed $\epsilon>0$, we have
\begin{equation}\label{eq:gllr-expected-stopping-time}
	\begin{aligned}
	\Exp_{[2]}[T_G^R]
	&\sim
	\frac{|\log\alpha|}
	{\displaystyle
		\min_{i\in[2]}
		c_i^{\mathrm G}(\epsilon)I^{\mathrm G}_i(\epsilon)}
	=
	\frac{|\log\alpha|}
	{c_1^{\mathrm G}(\epsilon)I^{\mathrm G}_1},
	\qquad \alpha\to0,\\
	\Exp_{[2]}[T_{AK}^R]
	&\sim
	\frac{|\log\alpha|}
	{\displaystyle
		\min_{i\in[2]}c_i(\epsilon)I_i(\epsilon)}
	=
	\frac{|\log\alpha|}
	{c_1(\epsilon)I_1(\epsilon)},
	\qquad \alpha\to0,
	\end{aligned}
\end{equation}
and
\begin{equation}\label{re}
	\begin{aligned}
		\cR(\epsilon)
		:=
		\lim_{\alpha\to0}
		\frac{\Exp_{[2]}[T_G^R]}
		{\Exp_{[2]}[T_{AK}^R]}=
		\frac{c_1(\epsilon)I_1(\epsilon)}
		{c_1^{\mathrm G}(\epsilon)I^{\mathrm G}_1}.
	\end{aligned}
\end{equation}
Letting $\epsilon \to 0$, $c_1(\epsilon)I_1(\epsilon)$ converges to a strictly positive constant, whereas $c_1^{\mathrm G}(\epsilon)I^{\mathrm G}_1 \to 0$. Therefore, $\cR(\epsilon) \to \infty$, which means that for sufficiently small $\epsilon$, the first-order
asymptotic expected stopping time of GLLR can be arbitrarily
larger than that of ALLR combined with KDE.

\subsection{Simulation}

We simulate two independent sources in which both
alternative hypotheses hold. The observations from source $i$ have
density
\begin{equation*}
	f_i
	=\frac12\cN(m_i-d,s^2)+\frac12\cN(m_i+d,s^2),
	\qquad m_1=1,\quad m_2=\epsilon,
\end{equation*}
where $s=0.25$ and $d=\sqrt{1-s^2}$. Thus, both true variances equal
one. We take $B=3$, $c=0.1$, and $C=2$. The shared Gaussian null has
mean in $[-B,0]$ and variance in $[c^2,C^2]$, while the Gaussian
surrogate alternative has mean in $[\epsilon,B]$ and the same variance
bounds.

\begin{figure}[htb!]
	\centering
	\includegraphics[width=10cm]{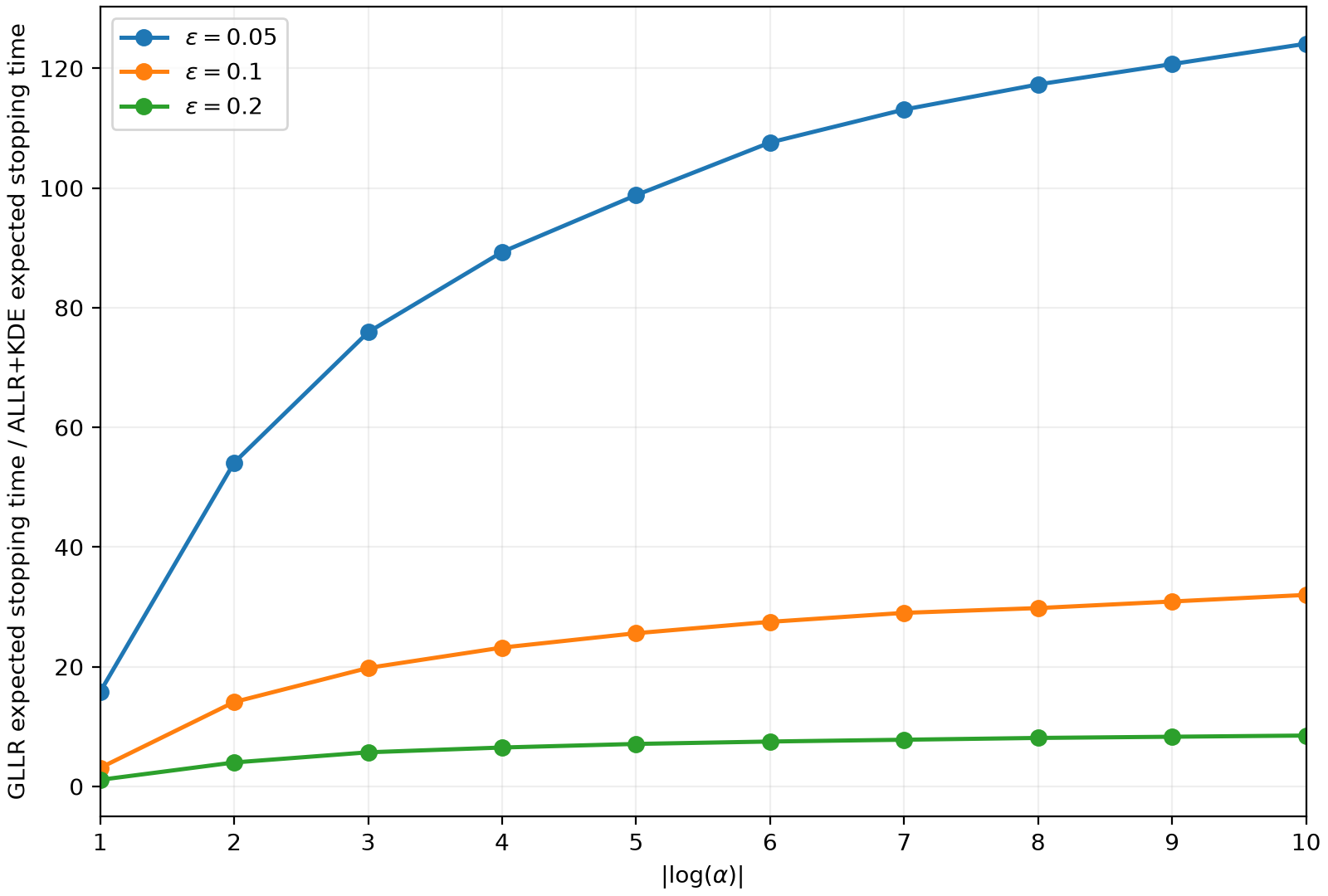}
	\caption{Ratio of expected stopping times versus $|\log(\alpha)|$.}
	\label{figur1}
\end{figure}

\textit{1)}
In the following simulations, we compare
$\Exp_{[2]}[T_G^R]$ and $\Exp_{[2]}[T_{AK}^R]$ for fixed, small
values of $\epsilon$ and decreasing values of $\alpha$. Figure~\ref{figur1} shows the ratio
$\Exp_{[2]}[T_G^R]/\Exp_{[2]}[T_{AK}^R]$ at target error
probabilities $\alpha$ versus $|\log(\alpha)| \in\{1,\ldots,10\}$. For each
target error, stopping thresholds were selected using simulations
so that the estimate of the error probability was as close
as possible to $\alpha$. The Monte Carlo standard error of
each estimated stopping time ratio is below $0.1$, for each value of $\epsilon$. The ratio is larger for smaller $\epsilon$, consistent with $\cR(\epsilon)\to\infty$ as $\epsilon\to0$.

\begin{figure}[htb!]
	\centering
	\includegraphics[width=10cm]{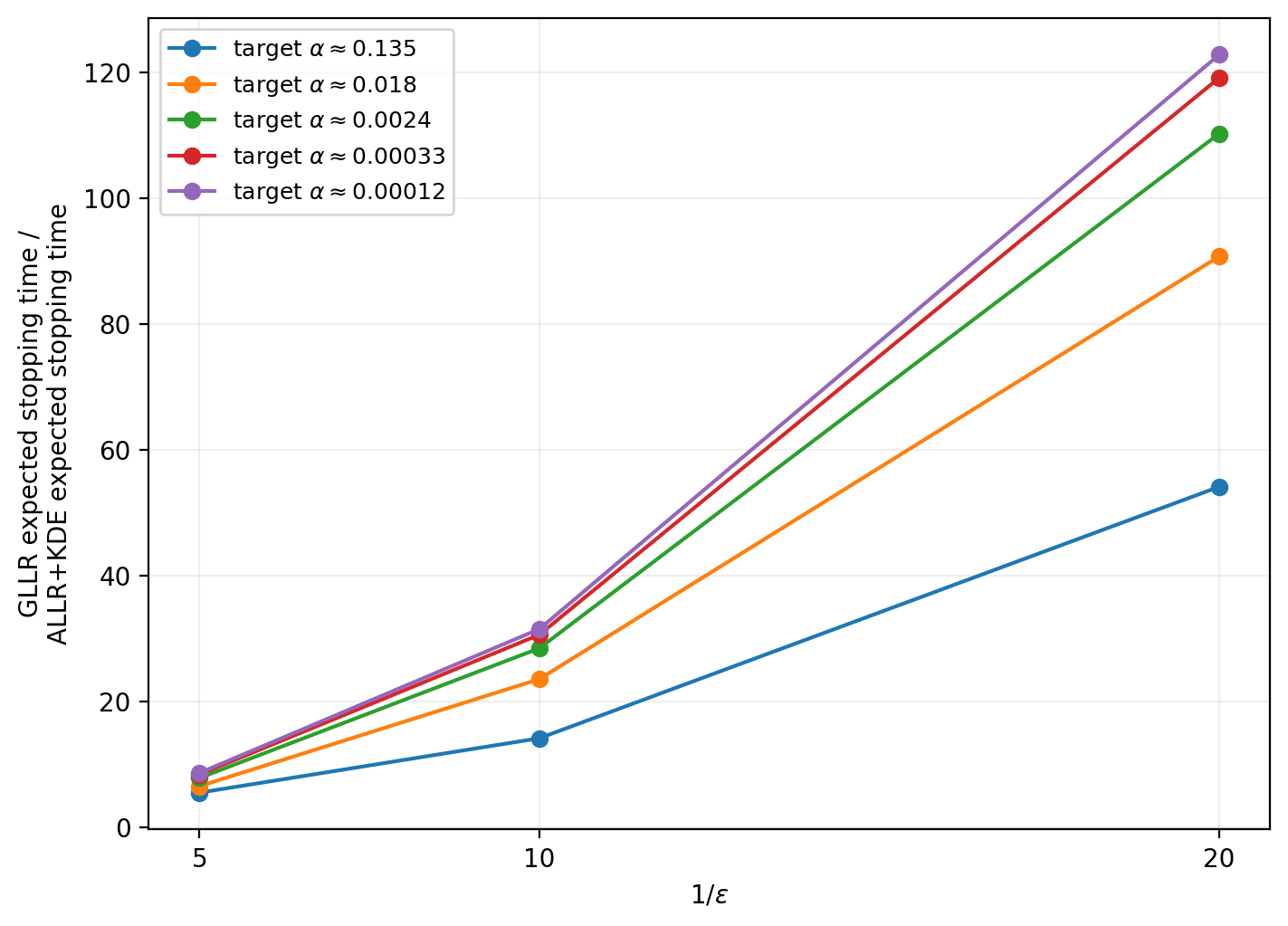}
	\caption{Ratio of expected stopping times versus $1/\epsilon$.}
	\label{figur2}
\end{figure}

\textit{2) }
In the following simulations, we compare
$\Exp_{[2]}[T_G^R]$ and $\Exp_{[2]}[T_{AK}^R]$ for small, fixed
target error probabilities $\alpha$. Figure~\ref{figur2} shows
the estimated ratio of the expected stopping time of the surrogate
GLLR to that of the ALLR with KDE versus $1/\epsilon$, for
$\epsilon\in\{0.05,0.1,0.2\}$. For each fixed $\epsilon$, smaller
target error probabilities yield larger ratios. The Monte Carlo standard error of every estimated ratio is below $0.1$.

\section{Conclusion}\label{concl}

We studied sequential anomaly identification under sampling
constraints and composite hypotheses, with particular emphasis
on the nonparametric setting. We derived a universal asymptotic
lower bound on the expected stopping time and established
sufficient conditions for attaining this bound as the error
probabilities vanish. These conditions identify the sampling
frequencies needed for first-order asymptotic optimality.
We also proposed suitable test statistics and density estimators
and showed, under appropriate assumptions, that they satisfy
the corresponding convergence conditions.

The framework also accommodates generalized error criteria
that tolerate a prescribed number of classification errors.
In \cite{gener24}, we developed the corresponding asymptotic
optimality conditions for simple hypotheses and a more
restrictive sampling constraint. The same methodology applies
directly to the present composite setting by replacing the
Kullback--Leibler divergences between simple hypotheses with
the information quantities defined in \eqref{inf1} and
\eqref{inf2}. Directions for future research include extending
the framework to non-i.i.d. observations from each source
\cite{nitv15} and incorporating switching costs into the
sampling process \cite{lamb21}.

\bibliographystyle{IEEEtran}
\bibliography{biblio_c}

\appendices
\renewcommand{\thetheorem}{ \Alph{section}.\arabic{theorem}}
\renewcommand{\thelemma}{ \Alph{section}.\arabic{lemma}}

\section{}\label{appf}
In this appendix, we verify that regularity conditions (a)--(g)
in Subsection~\ref{g_ass} hold for regular minimal exponential
families restricted to compact parameter sets. We also prove
Proposition~\ref{mixt}, which extends these conditions to
mixtures of densities from these families.

For each \(i\in[M]\), consider the exponential family
\[
g_\theta(x)
=
h_i(x)\exp\{\theta^\top T_i(x)-b_i(\theta)\},
\qquad \theta\in \Theta_i,
\]
where \(h_i:\bS\to(0,\infty)\) and
\(T_i:\bS\to\mathbb{R}^{d_i}\) are measurable,
\[
b_i(\theta)
=
\log\int_{\bS}
h_i(x)\exp\{\theta^\top T_i(x)\}\,d\nu_i(x),
\]
and
\[
\widetilde{\Theta}_i
=
\left\{
\theta\in\mathbb{R}^{d_i}:
0<
\int_{\bS}
h_i(x)\exp\{\theta^\top T_i(x)\}\,d\nu_i(x)
<\infty
\right\}
\]
is the natural parameter space. We impose the following
conditions:
\begin{enumerate}[(1)]
	\item The family is \emph{regular}: the natural
	parameter space \(\widetilde{\Theta}_i\) is a nonempty
	open subset of \(\mathbb{R}^{d_i}\).
	
	\item The parameter set \(\Theta_i\) is a nonempty
	compact subset of \(\widetilde{\Theta}_i\), equipped
	with the Euclidean topology.
	
	\item The family is \emph{minimal}: there do not
	exist a nonzero vector \(v\in\mathbb{R}^{d_i}\)
	and a constant \(c\in\mathbb{R}\) such that
	\[
	v^\top T_i(x)=c,
	\qquad \nu_i\text{-almost everywhere}.
	\]
	which implies that the parametrization is
	identifiable, i.e., for all \(\theta,\eta\in \Theta_i\),
	\[
	g_\theta=g_\eta,
	\quad \nu_i\text{-almost everywhere}
	\quad\Leftrightarrow\quad
	\theta=\eta.
	\]
\end{enumerate}

\begin{proposition}\label{lem:homeomorphism}
	Fix \(i\in[M]\). Let \(\Theta_i\) be a compact subset
	of the open natural parameter space of a regular minimal
	exponential family
	\[
	g_\theta(x)
	=
	h_i(x)\exp\{\theta^\top T_i(x)-b_i(\theta)\},
	\qquad \theta\in\Theta_i,
	\]
	and define
	\begin{equation}\label{eexp}
		\cE_i:=\{g_\theta:\theta\in\Theta_i\}.
	\end{equation}
	Equip \(\Theta_i\) with the Euclidean norm and
	\(\cE_i\) with the \(\sL^1\) norm. Then the mapping
	\[
	\Psi_i:\Theta_i\to \cE_i,
	\qquad
	\Psi_i(\theta)=g_\theta,
	\]
	is a homeomorphism, i.e., it is a bijection and for every sequence
	\(\{\theta_n\}_{n\ge1}\subset\Theta_i\) and every
	\(\theta\in\Theta_i\),
	\[
	\theta_n\rightarrow\theta
	\quad\Leftrightarrow\quad
	\|g_{\theta_n}-g_\theta\|_1\rightarrow0.
	\]
	Moreover, \(\cE_i\) is compact in \(\sL^1\), and both \(\Psi_i\) and \(\Psi_i^{-1}\) are
	uniformly continuous.
\end{proposition}

\begin{IEEEproof}
In order to show that $\Psi_i$ is bijective, it suffices to show that it is surjective and injective. Surjectivity follows from the definition of \(\cE_i\). To prove injectivity, suppose that
\(g_\theta=g_\eta\), \(\nu_i\)-almost everywhere.
Since \(h_i>0\), \(\nu_i\)-almost everywhere, taking
logarithms gives
\[
(\theta-\eta)^\top T_i(x)
=
b_i(\theta)-b_i(\eta),
\qquad \nu_i\text{-almost everywhere}.
\]
The minimality condition implies \(\theta=\eta\).

To prove that $\Psi_i$ is continuous, we consider a convergent sequence \(\theta_n\to\theta\) in \(\Theta_i\). By continuity of the log-partition function \(b_i\),
\[
g_{\theta_n}(x)\rightarrow g_\theta(x),
\qquad \nu_i\text{-almost everywhere}.
\]
Since these functions are probability densities,
Scheff\'e's lemma \cite{scheffe1947convergence} yields
\[
\|g_{\theta_n}-g_\theta\|_1\rightarrow0.
\]
Thus, \(\Psi_i\) is continuous.

Since \(\Psi_i\) is a continuous bijection from the
compact space \(\Theta_i\) onto the Hausdorff space
\(\cE_i\), by \cite[Theorem~26.6]{munkres2000topology} we get that \(\Psi_i^{-1}\) is continuous. Alternatively, we can prove that \(\Psi_i^{-1}\) is continuous as follows. Suppose that
\begin{equation}\label{gl0}
	g_{\eta_n}\rightarrow g_\eta
	\qquad\text{in }\sL^1.
\end{equation}
Since \(\Theta_i\) is compact, in order to prove that \(\eta_{n}\to\eta\) it suffices to show that every convergent subsequence of \(\{\eta_n \,:\, n \in \bN \}\)
has limit \(\eta\). Let \(\{\eta_{n_k} \,:\, k \in \bN \}\) be any convergent subsequence, with
\[
\eta_{n_k}\rightarrow\zeta\in\Theta_i.
\]
By continuity of \(\Psi_i\),
\[
g_{\eta_{n_k}}\rightarrow g_\zeta
\qquad\text{in }\sL^1.
\]
On the other hand, by \eqref{gl0},
\[
g_{\eta_{n_k}}\rightarrow g_\eta
\qquad\text{in }\sL^1.
\]
Since limits in \(\sL^1\) are unique,
\[
g_\eta=g_\zeta,
\qquad \nu_i\text{-almost everywhere}.
\]
Identifiability therefore implies \(\eta=\zeta\).
Consequently, \(\eta_{n}\to\eta\), proving continuity
of \(\Psi_i^{-1}\).

Furthermore, \(\cE_i\) is compact as the continuous
image of \(\Theta_i\). Since both spaces are compact metric spaces, the
Heine--Cantor theorem
\cite[Theorem~4.19]{rudin1976principles}
implies that \(\Psi_i\) and \(\Psi_i^{-1}\) are
uniformly continuous.
\end{IEEEproof}

\begin{proposition}\label{prop:expf}
	Suppose that $\Phi_i=\cE_i$ for every $i\in[M]$,
	where $\cE_i$ is defined in \eqref{eexp}.
	Then regularity conditions (a)--(g) in Subsection~\ref{g_ass}
	are satisfied.
\end{proposition}

\begin{IEEEproof} 
Fix \(i\in[M]\), and let \(f_i=g_{\theta_i}\). To simplify the notation, we use the change of measure
\[
d\mu_i(x)=h_i(x)\,d\nu_i(x).
\]	
	
\begin{enumerate}[(a)]
	\item  Condition (a) follows by Lemma \ref{lem:homeomorphism}.
	
	\item Condition (b) follows from \(h_i(x)>0\) for all \(x\in\bS\) and from the exponential form of \(g_\theta\).
	
	\item Let \(\eta\in\Theta_i\). Then
	\begin{equation}\label{qq}
		\begin{aligned}
			D(g_{\theta_i},g_\eta)
			&=
			\Exp_{g_{\theta_i}}
			\left[
			\log
			\frac{g_{\theta_i}(X_i)}
			{g_\eta(X_i)}
			\right] \\
			&=
			(\theta_i-\eta)^\top\nabla b_i(\theta_i)
			-
			\bigl(b_i(\theta_i)-b_i(\eta)\bigr).
		\end{aligned}
	\end{equation}
	
	Since the exponential family is regular and minimal, the log-partition function
	\(b_i\) is strictly convex. Therefore,
	\[
	D(g_{\theta_i},g_\eta)>0,
	\qquad
	\forall\,\eta\in\Theta_i,\ \eta\neq\theta_i.
	\]
	This proves condition (c).
	
	\item By Lemma \ref{lem:homeomorphism}, we have
	\begin{equation*}
	\|g_{\theta_n}-g_\theta\|_1\rightarrow0 \quad\Rightarrow\quad	 \theta_n\rightarrow\theta.
	\end{equation*}
	Then the continuity of the mapping
	\[
	g_\eta\mapsto D(g_{\theta_i},g_\eta)
	\]
	follows from the continuity of \(b_i\) on \(\Theta_i\) and \eqref{qq}.
	
	\item Fix \(\mathfrak q > 2\) and \(\epsilon >0\). Fix $\eta \in \Theta_i$. We have
	\begin{equation*}
		\log g_{\eta}(X_i)=\eta^\top T_i(X_i)-b_i(\eta),
	\end{equation*}
	and, by convexity,
	\begin{equation*}
		\begin{aligned}
			|\log g_{\eta}(X_i)|^{\mathfrak q}
			&\le
			2^{\mathfrak q-1}
			\left(
			|\eta^\top T_i(X_i)|^{\mathfrak q}
			+
			|b_i(\eta)|^{\mathfrak q}
			\right)\\
			&\le
			2^{\mathfrak q-1}
			\left(
			\|\eta\|_2^{\mathfrak q}
			\|T_i(X_i)\|_2^{\mathfrak q}
			+
			|b_i(\eta)|^{\mathfrak q}
			\right),
		\end{aligned}
	\end{equation*}
	where the last inequality follows from the Cauchy--Schwarz inequality. Since \(\Theta_i\) is compact and \(b_i\) is continuous on \(\Theta_i\), there exist constants
	\(c_1,c_2>0\) such that
	\begin{equation*}
		\sup_{\eta\in\Theta_i}
		|\log g_{\eta}(X_i)|^{\mathfrak q}
		\le
		c_1
		+
		c_2 \sum_{k=1}^{d_i} |T_{i,k}(X_i)|^{\mathfrak q}.
	\end{equation*}
	Consequently, since \(V_\epsilon(g_\eta)\subset\Phi_i\),
	\begin{equation}\label{stp}
		\sup_{h\in V_\epsilon(g_\eta)}
		|\log h(X_i)|^{\mathfrak q}
		\le c_1 + c_2 \sum_{k=1}^{d_i} |T_{i,k}(X_i)|^{\mathfrak q}.
	\end{equation}
	
	It therefore suffices to show that all moments of \(T_i(X_i)\) are uniformly bounded over \(\Theta_i\). Since \(\Theta_i\) is a compact subset of the open set \(\widetilde{\Theta}_i\), there exists \(\delta>0\) such that
	\begin{equation*}
		\theta_i+u\in\widetilde{\Theta}_i,
	\end{equation*}
	for all \(\theta_i \in \Theta_i\) and all \(u \in \bR^{d_i}\) with \(\|u\|_{\infty} < \delta\). Fix \(\delta' < \delta\). Since
	\begin{equation*}
		\exp\left\{
		\delta'
		\sum_{k=1}^{d_i}
		|T_{i,k}(X_i)|
		\right\}
		\le
		\sum_{\sigma\in\{-1,1\}^{d_i}}
		\exp\left\{
		\delta'\,\sigma^\top T_i(X_i)
		\right\},
	\end{equation*}
	we have
	\begin{equation*}
		\begin{aligned}
			\Exp_{g_{\theta_i}}\left[
			\exp\left\{
			\delta'
			\sum_{k=1}^{d_i}
			|T_{i,k}(X_i)|
			\right\}\right]
			&\le
			\sum_{\sigma\in\{-1,1\}^{d_i}}
			\Exp_{g_{\theta_i}}\left[
			\exp\left\{
			\delta'\,\sigma^\top T_i(X_i)
			\right\}\right].
		\end{aligned}
	\end{equation*}
	
	Fix \(\sigma\in\{-1,1\}^{d_i}\). Since \(\|\delta' \sigma\|_\infty<\delta\), it follows that \(\theta_i+\delta'\sigma\in\widetilde{\Theta}_i\). Therefore,
	\begin{equation*}
		\begin{aligned}
			\Exp_{g_{\theta_i}}\left[
			\exp\left\{
			\delta'\,\sigma^\top T_i(X_i)
			\right\}\right]
			&=
			\int
			\exp\left\{
			\delta'\,\sigma^\top T_i(x)
			\right\}
			e^{\theta_i^\top T_i(x)-b_i(\theta_i)}
			\,d\mu_i(x)\\
			&=
			e^{-b_i(\theta_i)}
			\int
			\exp\left\{
			(\theta_i+\delta'\sigma)^\top T_i(x)
			\right\}
			\,d\mu_i(x)\\
			&=
			\exp\left\{
			b_i(\theta_i+\delta'\sigma)-b_i(\theta_i)
			\right\}
			<
			\infty.
		\end{aligned}
	\end{equation*}
	
	Since \(b_i(\theta_i+\delta'\sigma)-b_i(\theta_i)\) is continuous with respect to \(\theta_i\) on the compact set \(\Theta_i\), we conclude that
	\begin{equation}\label{unb}
		\Exp_{g_{\theta_i}}\left[
		\exp\left\{
		\delta'\,\sigma^\top T_i(X_i)
		\right\}\right]
		\le
		\sup_{\zeta\in\Theta_i}
		\exp\left\{
		b_i(\zeta+\delta'\sigma)-b_i(\zeta)
		\right\}
		<
		\infty,
	\end{equation}
	which proves the claim.
	
	\item  We prove condition (f) by constructing an integrable
	upper envelope for the entire family. Since $\Theta_i$ is
	compact and $\widetilde{\Theta}_i$ is open, each point of
	$\Theta_i$ lies in the interior of a closed cube contained
	in $\widetilde{\Theta}_i$. The interiors of finitely many
	such cubes cover $\Theta_i$. Let $\mathcal V_i$ denote the
	set of all their vertices. These vertices need not belong
	to $\Theta_i$, but they belong to $\widetilde{\Theta}_i$,
	so each $v\in\mathcal V_i$ defines a density $g_v$.
	
	Fix $\theta\in\Theta_i$ and $x\in\bS$, and choose a cube
	containing $\theta$. Over this cube, the linear function
	$u\mapsto u^\top T_i(x)$ is maximized at a vertex $v$:
	in coordinate $k$, choose the upper endpoint if
	$T_{i,k}(x)\geq0$ and the lower endpoint otherwise.
	Thus,
	\[
	\theta^\top T_i(x)\leq v^\top T_i(x).
	\]
	Also, by continuity of $b_i$ and compactness of $\Theta_i$,
	the minimum $m_i:=\min_{\theta\in\Theta_i}b_i(\theta)$
	is finite. Combining these two bounds gives
	\[
	g_\theta(x)
	=
	h_i(x)e^{\theta^\top T_i(x)-b_i(\theta)}
	\leq
	h_i(x)e^{v^\top T_i(x)-m_i}
	=
	e^{b_i(v)-m_i}g_v(x).
	\]
	The selected vertex may depend on $\theta$ and $x$, but
	it always belongs to the same finite set $\mathcal V_i$.
	Hence, the sum over all vertices bounds $g_\theta(x)$
	for every $\theta\in\Theta_i$. Moreover, since
	$V_\epsilon(g_\eta)\subseteq\Phi_i$, we obtain
	\[
	p_\epsilon^+(g_\eta,x)
	\leq
	\sup_{\theta\in\Theta_i}g_\theta(x)
	\leq
	\sum_{v\in\mathcal V_i}e^{b_i(v)-m_i}g_v(x).
	\]
	Finally, each $g_v$ integrates to one. Therefore,
	\[
	\Exp_{f_i}\!\left[
	\frac{p_\epsilon^+(g_\eta,X_i)}{f_i(X_i)}
	\right]
	=
	\int_{\bS}p_\epsilon^+(g_\eta,x)\,d\nu_i(x)
	\leq
	\sum_{v\in\mathcal V_i}e^{b_i(v)-m_i}
	<\infty.
	\]
	This bound is independent of $\eta$, $f_i$, and $\epsilon$.
	
	\item For \(\eta\in\Theta_i\), and \(g_\zeta\in V_\epsilon(g_\eta)\),
	\begin{equation*}
		\begin{aligned}
			\left|
			\log g_\eta(X_i)-\log g_\zeta(X_i)
			\right|
			&=
			\left|
			(\eta-\zeta)^\top T_i(X_i)
			-
			b_i(\eta)
			+
			b_i(\zeta)
			\right|\\
			&\le
			\|\eta-\zeta\|_2
			\|T_i(X_i)\|_2
			+
			|b_i(\eta)-b_i(\zeta)|.
		\end{aligned}
	\end{equation*}
	
	Since \(b_i\) is continuously differentiable on the open natural parameter
	space and \(\Theta_i\) is compact, \(b_i\) is Lipschitz on \(\Theta_i\). Hence,
	there exists a constant \(C>0\) such that
	\[
	\left|
	\log g_\eta(X_i)-\log g_\zeta(X_i)
	\right|
	\le
	\|\eta-\zeta\|_2
	\left(
	\|T_i(X_i)\|_2+C
	\right).
	\]
	
	Since \(\|\eta-\zeta\|_2\to0\) as \(\epsilon\to0\), and all moments of
	\(T_{i,k}(X_i)\) are uniformly bounded over \(\Theta_i\) by \eqref{unb}, the dominated
	convergence theorem yields $H_{g_{\theta_i}}(g,\epsilon)\rightarrow0$, as \(\epsilon\rightarrow0\), uniformly over all \(\theta_i \in\Theta_i\).
\end{enumerate}
\end{IEEEproof}

We proceed to the proof of Proposition \ref{mixt}.

\begin{IEEEproof}[Proof of Proposition \ref{mixt}]
Fix $i\in[M]$. By Proposition~\ref{prop:expf}, conditions
(a)--(g) hold for $\cE_i$. We show that they also hold
for the mixture class 
\[
\Phi_i
:=
\{q_{i,P}:P\in\mathcal P(\Theta_i)\}.
\]
The proof for condition (d) is provided last.

\begin{itemize}
	\item [(a)] To prove condition (a) it suffices to show that $\mathcal P(\Theta_i)$ equipped with its weak topology is compact, and that the map $P\mapsto q_{i,P}$ is continuous from $\mathcal P(\Theta_i)$ into $\sL^1$. Then, $\Phi_i$ is compact in $\sL^1$, as the continuous image of a compact space.
	
	Since $\Theta_i$ is compact and $P(\Theta_i)=1$ for every
	$P\in\mathcal P(\Theta_i)$, the family $\mathcal P(\Theta_i)$ is uniformly tight. By Prokhorov's theorem \cite[Theorem~5.1]{billingsley1999convergence},
	every sequence in $\mathcal P(\Theta_i)$ has a weakly
	convergent subsequence with limit in $\mathcal P(\Theta_i)$. Since the weak topology is metrizable, this establishes
	compactness of $\mathcal P(\Theta_i)$.
	
	To show the continuity of the map $P\mapsto q_{i,P}$, it suffices to show that for every sequence $P_n \xrightarrow{w} P$, where $\xrightarrow{w}$ stands for weak convergence, it holds
	\begin{equation}\label{qip}
	q_{i,P_n} \rightarrow q_{i,P}, \quad \mbox{in} \; \sL^1.
	\end{equation}
	We note that for each fixed $x\in\bS$, the exponential-family
	representation
	\[
	g_\theta(x)
	=
	h_i(x)\exp\{\theta^\top T_i(x)-b_i(\theta)\}
	\]
	and the continuity of $b_i$ imply that
	$\theta\mapsto g_\theta(x)$ is continuous. Since $\Theta_i$
	is compact, this map is also bounded on $\Theta_i$,
	with a bound that may depend on $x$. By definition of weak convergence, the  $P_n \xrightarrow{w} P$ implies
	\begin{equation}\label{ptc}
	q_{i,P_n}(x)
	=
	\int_{\Theta_i}g_\theta(x)\,P_n(d\theta)
	\rightarrow
	\int_{\Theta_i}g_\theta(x)\,P(d\theta)
	=
	q_{i,P}(x).
	\end{equation}
	Since $q_{i,P_n}$ and $q_{i,P}$ are probability densities, the pointwise convergence \eqref{ptc} implies the $\sL^1$ convergence \eqref{qip} by Scheff\'e's lemma \cite{scheffe1947convergence}.
	
	\item [(b)] Each $g_\theta(x)$ is strictly positive for every $x\in\bS$. Therefore,
	\[
	q_{i,P}(x)
	=
	\int_{\Theta_i}g_\theta(x)\,P(d\theta)
	>0
	\]
	for every $P\in\mathcal P(\Theta_i)$ and $x\in\bS$.
	
	\item [(c)] If $f_i$ and $g$ are distinct as elements of $\Phi_i$ with $\sL^1$, then $D(f_i,g)>0$, by Pinsker's inequality.
	
	\item [(e)] In order to prove the claim it suffices to show that for the envelope 
	 \begin{equation*}
	  L_i(x):=\sup_{g\in\cE_i}|\log g(x)|,
	 \end{equation*}
	it holds
	\begin{equation}\label{st1}
	 \sup_{f_i\in\cE_i} \Exp_{f_i}[L_i(X_i)^{\mathfrak q}]<\infty.
	\end{equation}
	Indeed, for every  $P\in\mathcal P(\Theta_i)$, Tonelli's theorem gives
	\[
	\begin{aligned}
		\Exp_{q_{i,P}}[L_i(X_i)^{\mathfrak q}]
		&=
		\int_{\Theta_i}
		\Exp_{g_\theta}[L_i(X_i)^{\mathfrak q}]
		\,P(d\theta)\\
		&\leq
		\sup_{f_i\in\cE_i}
		\Exp_{f_i}[L_i(X_i)^{\mathfrak q}].
	\end{aligned}
	\]
	Also, for every mixture $q_{i,P}(x)$, it holds $|\log q_{i,P}(x)|\leq L_i(x)$, and thus for every $\epsilon>0$,
	\begin{equation}\label{egt}
	\begin{aligned}
	\sup_{f_i,g\in\Phi_i}
	\Exp_{f_i}\!\left[
	\sup_{p\in V_\epsilon(g)}
	|\log p(X_i)|^{\mathfrak q}
	\right] &\leq \sup_{f_i\in\Phi_i}
	\Exp_{f_i}\!\left[
	L_i(X_i)^{\mathfrak q}
	\right]\\
	&\leq
	\sup_{f_i\in\cE_i}
	\Exp_{f_i}[L_i(X_i)^{\mathfrak q}]
	<\infty,
	\end{aligned}
	\end{equation}
	which shows (e) for $\Phi_i$. To prove \eqref{st1}, we apply condition (e) to $\cE_i$, and we cover $\cE_i$ by finitely many neighborhoods of the radius specified
	in that condition. Such a finite cover exists by compactness
	of $\cE_i$. Then $L_i(x)$ is bounded by the sum of the corresponding local upper envelopes, each of which is integrable by Proposition~\ref{prop:expf}.
	
	\item [(f)] In order to prove the claim it suffices to show that for the envelope
	\[
	U_i(x):=\sup_{g\in\cE_i}g(x),
	\]
	it holds
	\begin{equation}\label{str2}
	 \int_{\bS}U_i(x)\,d\nu_i(x)<\infty.
	\end{equation}
	Indeed, for every mixture $q_{i,P}(x)\leq U_i(x)$, and thus for all $f_i,g\in\Phi_i$
	and every $\epsilon>0$,
	\[
	\begin{aligned}
		\Exp_{f_i}\!\left[
		\frac{p_\epsilon^+(g,X_i)}{f_i(X_i)}
		\right]
		&=
		\int_{\bS}p_\epsilon^+(g,x)\,d\nu_i(x)\\
		&\leq
		\int_{\bS}U_i(x)\,d\nu_i(x)
		<\infty,
	\end{aligned}
	\]
	where the equality uses the strict positivity of $f_i$. The upper bound is independent of both $f_i$ and $g$. To show \eqref{str2}, we apply condition (f) to $\cE_i$, and we cover $\cE_i$ by finitely many neighborhoods of the radius specified
	in that condition. Such a finite cover exists by compactness
	of $\cE_i$. Then $U_i(x)$ is bounded by the sum of the corresponding local upper envelopes, each of which is integrable by Proposition~\ref{prop:expf}.
	
	\item [(g)] Fix $\epsilon >0$, and define
	  \begin{equation}\label{wex}
	  W_{i,\epsilon}(x)
	  :=
	  \sup_{\substack{
	  		p,q\in\Phi_i\\
	  		\|p-q\|_1<\epsilon
	  }}
	  |\log p(x)-\log q(x)|.
	  \end{equation}
	 We observe that for each fixed $f_i \in \Phi_i$,
	\[
	\begin{aligned}
		&\sup_{g\in\Phi_i}
		\Exp_{f_i}\!\left[
		\sup_{p\in V_\epsilon(g)}
		|\log p(X_i)-\log g(X_i)|
		\right]\\
		&\qquad\leq
		\Exp_{f_i}[W_{i,\epsilon}(X_i)].
	\end{aligned}
	\]
	Thus, in order to prove the claim it suffices to show that
	\begin{equation}\label{tsh3}
	W_{i,\epsilon}(x)\rightarrow0
	\qquad\text{as }\epsilon\downarrow0,
	\end{equation}
	for $\nu_i$-almost every $x$. Then by dominated convergence theorem, we deduce that
	\begin{equation}\label{eq1}
		\Exp_{f_i}[W_{i,\epsilon}(X_i)]\rightarrow0
		\qquad\text{as }\epsilon\downarrow0,
	\end{equation}
	because for every mixture $q_{i,P}(x)$, it holds $|\log q_{i,P}(x)|\leq L_i(x)$, and thus
	\[
	W_{i,\epsilon}(x) \leq 2 L_i(x), \quad \nu_i \; \mbox{almost everywhere,}
	\]
	where $\sup_{f_i\in\Phi_i} \Exp_{f_i}[L_i(X_i)]<\infty$  by  \eqref{egt}. 
	
	To this end, for each $P\in\mathcal P(\Theta_i)$ we write
	\[
	q_{i,P}(x)=h_i(x)F_{i,P}(T_i(x)),
	\]
	where
	\[
	F_{i,P}(t)
	:=
	\int_{\Theta_i}
	\exp\{\theta^\top t-b_i(\theta)\}\,P(d\theta),
	\qquad t\in\bR^{d_i}.
	\]
	The function $F_{i,P}$ is strictly positive because its
	integrand is strictly positive. It is also continuous in $t$
	by dominated convergence, since the integrand is continuous
	and bounded on $\Theta_i\times K$ for every compact set
	$K\subseteq\bR^{d_i}$.  
	
	Fix $\theta_{i,0}\in\Theta_i$, and let $S_{T_i}$ denote
	the support of the distribution of $T_i(X_i)$ under
	$g_{\theta_{i,0}}$. If $q_{i,P}=q_{i,Q}$
	$\nu_i$-almost everywhere, then, since $h_i$ is strictly
	positive,
	\[
	F_{i,P}(T_i(X_i))=F_{i,Q}(T_i(X_i))
	\]
	almost surely under $g_{\theta_{i,0}}$. It follows that
	$F_{i,P}(t)=F_{i,Q}(t)$ for every $t\in S_{T_i}$.
	Indeed, if the two functions differed at some $t\in S_{T_i}$,
	their continuity would imply that they differ throughout
	an open neighborhood of $t$. By the definition of support,
	this neighborhood has positive probability under the law of
	$T_i(X_i)$, contradicting the almost sure equality.
	
	In order to prove \eqref{tsh3}, it suffices to show that for each $t\in S_{T_i}$,
	\begin{equation}\label{eq:mixture-pointwise-oscillation}
		\lim_{\epsilon\downarrow0}
		\sup_{\substack{
				P,Q\in\mathcal P(\Theta_i)\\
				\|q_{i,P}-q_{i,Q}\|_1<\epsilon
		}}
		|\log F_{i,P}(t)-\log F_{i,Q}(t)|
		=0.
	\end{equation}
	We proceed by contradiction. Fix $t\in S_{T_i}$ and suppose that there
	exist $\eta>0$ and sequences $P_n,Q_n\in\mathcal P(\Theta_i)$
	such that
	\[
	\|q_{i,P_n}-q_{i,Q_n}\|_1\rightarrow0,
	\qquad
	|\log F_{i,P_n}(t)-\log F_{i,Q_n}(t)|\geq\eta
	\quad\text{for every }n.
	\]
	By compactness of $\mathcal P(\Theta_i)$ in the weak topology,
	there exists a subsequence $\{P_{n_k}\}$ converging weakly to
	some $P\in\mathcal P(\Theta_i)$. Applying compactness again
	to $\{Q_{n_k}\}$, we obtain a further subsequence
	$\{Q_{n_{k_\ell}}\}$ converging weakly to some
	$Q\in\mathcal P(\Theta_i)$. Along these same indices,
	$P_{n_{k_\ell}}$ still converges weakly to $P$.
	Relabeling this common subsequence by $m$, we thus have
	$P_m\xrightarrow{w} P$ and $Q_m\xrightarrow{w} Q$.
	
	In part (a), we proved the continuity of the map $P\mapsto q_{i,P}$ from $\mathcal P(\Theta_i)$ equipped with its weak topology into $\sL^1$. Thus,
	\[
	\|q_{i,P_m}-q_{i,P}\|_1\rightarrow0,
	\qquad
	\|q_{i,Q_m}-q_{i,Q}\|_1\rightarrow0.
	\]
	Since $\|q_{i,P_m}-q_{i,Q_m}\|_1\to0$ by construction,
	the triangle inequality gives
	\[
	\|q_{i,P}-q_{i,Q}\|_1
	\leq
	\|q_{i,P}-q_{i,P_m}\|_1
	+\|q_{i,P_m}-q_{i,Q_m}\|_1
	+\|q_{i,Q_m}-q_{i,Q}\|_1
	\rightarrow0.
	\]
	Hence, $q_{i,P}=q_{i,Q}$ $\nu_i$-almost everywhere.
	Since $t\in S_{T_i}$, the support argument above yields
	\[
	F_{i,P}(t)=F_{i,Q}(t)>0.
	\]
	For this fixed $t$, the map $\theta\mapsto\exp\{\theta^\top t-b_i(\theta)\}$ is bounded
	and continuous on $\Theta_i$. By weak convergence,
	\[
	F_{i,P_m}(t)\rightarrow F_{i,P}(t),
	\qquad
	F_{i,Q_m}(t)\rightarrow F_{i,Q}(t).
	\]
	Since the two limits coincide and are strictly positive,
	continuity of the logarithm gives
	\[
	|\log F_{i,P_m}(t)-\log F_{i,Q_m}(t)|
	\rightarrow0.
	\]
	This contradicts the lower bound $\eta>0$.
	
	\item [(d)] The log-moment bound \eqref{st1} implies that $D(f_i,p)$
	is finite for every $f_i,p\in\Phi_i$. Whenever
	$p,q\in\Phi_i$ satisfy $\|p-q\|_1<\epsilon$,
	\[
	\begin{aligned}
		|D(f_i,p)-D(f_i,q)|
		&=
		\left|
		\Exp_{f_i}[\log q(X_i)-\log p(X_i)]
		\right|\\
		&\leq
		\Exp_{f_i}[W_{i,\epsilon}(X_i)].
	\end{aligned}
	\]
	where $W_{i,\epsilon}(x)$ is defined in \eqref{wex}. By \eqref{eq1}, the right-hand side tends to zero with $\epsilon \to 0$, proving the continuity required in (d).
\end{itemize}
\end{IEEEproof}

\section{}\label{app0}
In Appendix~\ref{app0}, we prove Theorems~\ref{tGLLR}
and~\ref{thm_allr}, establishing the required asymptotic
properties of the GLLR and ALLR statistics. The proofs
use Lemma~\ref{aux1}, which provides an exponential bound
on the probability that the sampled sum defined below
is negative while the empirical sampling frequency is
at least a fixed positive constant. To state this lemma, we introduce the following notation.

For each $i\in[M]$, $g\in\Phi_i$, $\epsilon>0$, and $\delta \in \bR$
define
\begin{equation*}
	z_i(X_i(m))
	:=
	\log f_i(X_i(m))
	-
	\log p_\epsilon^+(g,X_i(m))
	-
	\delta,
	\qquad
	m\in\bN,
\end{equation*}
and
\begin{equation*}
	S_i(n)
	:=
	\sum_{m=1}^n
	z_i(X_i(m))R_i(m),
	\qquad
	n\in\bN.
\end{equation*}

\begin{lemma}\label{aux1}
	If \(\delta<D(f_i,g)\), then there exists \(\epsilon>0\) such that, for every \(\rho>0\),
	\[
	\Pro_{\boldsymbol f}
	\left(
	S_i(n) \leq 0,\,
	\pi_i^R(n)\ge\rho
	\right), \quad n \in \bN,
	\]
	is exponentially decaying.
\end{lemma}

\begin{IEEEproof}
To prove the claim, it suffices to show that there exist constants \(s>0\) and
\(b>0\) such that
\begin{equation}\label{eq:mgf}
	\Exp_{\boldsymbol f}
	\left[
	e^{-sz_i(X_i(1))}
	\right]
	\le
	e^{-b}.
\end{equation}
Indeed, since \(R_i(n)\in\{0,1\}\) is \(\mathcal{F}_{n-1}^R\)-measurable and
\(X_i(n)\) is independent of \(\mathcal{F}_{n-1}^R\), we have
\begin{equation}\label{eq:condmgf}
	\begin{aligned}
		\Exp_{\boldsymbol f}
		\left[
		e^{-sz_i(X_i(n))R_i(n)}
		\,\middle|\,
		\mathcal{F}_{n-1}^R
		\right]
		&=
		\left(
		\Exp_{\boldsymbol f}
		\left[
		e^{-sz_i(X_i(n))}
		\right]
		\right)^{R_i(n)}\\
		&\le
		e^{-bR_i(n)}, \quad \mbox{a.s.},
	\end{aligned}
\end{equation}
where the last inequality follows from \eqref{eq:mgf}. Consequently, the process
\[
Q(n)
:=
\exp\left\{
-sS_i(n)+bN_i^R(n)
\right\},
\qquad
n\in\bN,
\]
is a nonnegative supermartingale. Therefore,
\[
\Exp_{\boldsymbol f}[Q(n)]
\le
\Exp_{\boldsymbol f}[Q(0)]
=
1.
\]
On the event
\[
\left\{
S_i(n) \leq 0,\,
N_i^R(n)\ge\rho n
\right\},
\]
we have
\[
Q(n)\ge e^{b\rho n}.
\]
Hence,
\begin{equation*}
	\begin{aligned}
		\Pro_{\boldsymbol f}
		\left(
		S_i(n)\leq 0,\,
		N_i^R(n)\ge\rho n
		\right)
		&\le
		\Pro_{\boldsymbol f}
		\left(
		Q(n)\ge e^{b\rho n}
		\right)\\
		&\le
		e^{-b\rho n}
		\Exp_{\boldsymbol f}[Q(n)]\\
		&\le
		e^{-b\rho n},
	\end{aligned}
\end{equation*}
which proves the claim. To prove \eqref{eq:mgf}, we consider the function
\begin{equation*}
	w(t):=
	\Exp_{\boldsymbol f}
	\left[
	e^{-tz_i(X_i(1))}
	\right].
\end{equation*}
By the dominated convergence theorem, in order to prove that \(w(t)\) is twice differentiable on the interval \([0,\gamma)\), for some $\gamma >0$, it suffices to show that
\begin{equation}\label{eq:intcond}
	\Exp_{\boldsymbol f}
	\left[
	z_i^2(X_i(1))
	e^{-tz_i(X_i(1))}
	\right]
	<
	\infty,
	\qquad
	\forall\,t\in[0,\gamma),
\end{equation}
which follows form the fact that for every $t\in[0,a]$, $a\in(0,\gamma)$, and $k\in\{1,2\}$, it holds
\[
|z_i(X_i(1))|^k e^{-tz_i(X_i(1))}
\le
e^a+
z_i^2(X_i(1))
\left(1+e^{-az_i(X_i(1))}\right), \quad \mbox{a.s.}
\]

By Hölder's inequality,
\begin{equation*}
	\begin{aligned}
		\Exp_{\boldsymbol f}
		\left[
		z_i^2(X_i(1))
		e^{-tz_i(X_i(1))}
		\right]
		&\le
		\Exp_{\boldsymbol f}
		\left[
		|z_i(X_i(1))|^{\mathfrak q}
		\right]^{2/\mathfrak q}
		\Exp_{\boldsymbol f}
		\left[
		e^{-\frac{\mathfrak q}{\mathfrak q-2}tz_i(X_i(1))}
		\right]^{(\mathfrak q-2)/\mathfrak q}.
	\end{aligned}
\end{equation*}
Since \(\log\) is an increasing function, the first term on the right-hand side is bounded by regularity condition (e). For the second term, it suffices to show that there exists a constant \(C>0\) such that
\begin{equation}\label{eq:expbound}
	\Exp_{\boldsymbol f}
	\left[
	e^{-uz_i(X_i(1))}
	\right]
	\le
	C,
	\qquad
	\forall\,u\in[0,1).
\end{equation}
Indeed, for every \(u\in[0,1)\),
\begin{equation*}
	\begin{aligned}
		\int
		e^{-uz_i(x)}
		f_i(x)\,d\nu_i(x)
		&=
		\int_{\{z_i(x)\ge0\}}
		e^{-uz_i(x)}
		f_i(x)\,d\nu_i(x)
		+
		\int_{\{z_i(x)<0\}}
		e^{-uz_i(x)}
		f_i(x)\,d\nu_i(x)\\
		&\le
		1
		+
		\int_{\{z_i(x)<0\}}
		e^{-z_i(x)}
		f_i(x)\,d\nu_i(x)\\
		&\le
		\Exp_{\boldsymbol f}
		\left[
		e^{-z_i(X_i(1))}
		\right]
		+1=
		e^\delta
		\Exp_{\boldsymbol f}
		\left[
		\frac{
			p_\epsilon^+(g,X_i(1))
		}{
			f_i(X_i(1))
		}
		\right]+1,
	\end{aligned}
\end{equation*}
and \eqref{eq:expbound} follows from regularity condition (f). Therefore, \(w(t)\) is twice differentiable for all $t \in [0,\gamma)$, where $\gamma = (\mathfrak q-2)/\mathfrak q$, and
\begin{equation*}
	w'(t)
	=
	\Exp_{\boldsymbol f}
	\left[
	-z_i(X_i(1))
	e^{-tz_i(X_i(1))}
	\right], \quad  	w''(t)
	=
	\Exp_{\boldsymbol f}
	\left[
	z_i^2(X_i(1))
	e^{-tz_i(X_i(1))}
	\right].
\end{equation*}
Therefore, $w(0)=1$ and
\begin{equation*}
	\begin{aligned}
		w'(0)
		&=
		\Exp_{\boldsymbol f}
		\left[
		-z_i(X_i(1))
		\right]\\
		&=
		\Exp_{\boldsymbol f}
		\left[
		\sup_{h\in V_\epsilon(g)}
		\log h(X_i(1))
		-
		\log f_i(X_i(1))
		\right]
		+\delta\\
		&=
		\Exp_{\boldsymbol f}
		\left[
		\sup_{h\in V_\epsilon(g)}
		\Bigl(
		\log h(X_i(1))
		-
		\log g(X_i(1))
		\Bigr)
		\right]
		+
		\Exp_{\boldsymbol f}
		\left[
		\log g(X_i(1))
		-
		\log f_i(X_i(1))
		\right]
		+\delta\\
		&\le
		H_{\boldsymbol f}(g,\epsilon)-D(f_i,g)+\delta.
	\end{aligned}
\end{equation*}
By regularity condition (g), there exists \(\epsilon>0\) sufficiently small such that
\[
H_{\boldsymbol f}(g,\epsilon)
<
\frac{D(f_i,g)-\delta}{2},
\]
and therefore
\[
w'(0)
<
-
\frac{D(f_i,g)-\delta}{2}
<
0.
\]
Moreover,
\[
w''(t)\le D,
\qquad
\forall\,t\in\left[0,(\mathfrak q-2)/\mathfrak q\right),
\]
where \(D>0\) is a constant. Hence, \(w'\) is Lipschitz continuous on $\left[0,(\mathfrak q-2)/\mathfrak q\right)$. Since \(w(0)=1\) and the right derivative \(w'_+(0)<0\), there exist constants \(s>0\) and \(b>0\) such that
\[
w(s)=e^{-b}<1,
\]
which proves \eqref{eq:mgf}.
\end{IEEEproof}

\begin{IEEEproof}[Proof of Theorem \ref{tGLLR}]
We prove the result for $i \in A$, as the case $i \notin A$ follows analogously. Fix $i \in A$, $\rho \in (0,1]$, and $\epsilon \in (0, I_i(f_i))$. By the definition of $l_i^1(n)$ in \eqref{suprem}, we have
\begin{equation}\label{inq1}
	l_i^1(n) \geq l_i(n; f_i), \quad \forall n \in \bN.
\end{equation}
Since $\Phi_i^0$ is compact, there exists a finite $\delta$-cover $\{g_j\}_{j=1}^{L(\delta)} \subset \Phi_i^0$, and
\begin{equation*}
	l_i^0(n)
	\leq 
	\max_{1\le j\le L(\delta)}
	\;
	\sup_{g\in V_{\delta}(g_j)}
	l_i(n;g).
\end{equation*}
To prove the claim, it suffices to show that for each $j \in [L(\delta)]$,
\begin{equation*}
	\Pro_{\boldsymbol f}
	\left(
	l_i(n;f_i)
	-
	\sup_{g\in V_\delta(g_j)}
	l_i(n;g)
	<
	\left(
	\rho I_i(f_i)-\frac{\epsilon}{2}
	\right)n,
	\,
	\pi_i^R(n)\ge\rho
	\right)
\end{equation*}
is exponentially decaying. Fix $j \in [L(\delta)]$. Since \(N_i^R(n)\ge \rho n\), and
\[
\sup_{g\in V_\delta(g_j)}
l_i(n;g)
\le
\sum_{m=1}^{n}
\log p_\delta^+(g_j,X_i(m))R_i(m),
\]
we have
\begin{equation*}
	\begin{aligned}
		&\Pro_{\boldsymbol f}
		\left(
		l_i(n;f_i)
		-
		\sup_{g\in V_\delta(g_j)}
		l_i(n;g)
		<
		\left(
		\rho I_i(f_i)-\frac{\epsilon}{2}
		\right)n,
		\,
		\pi_i^R(n)\ge\rho
		\right)\\
		&\le
		\Pro_{\boldsymbol f}
		\left(
		l_i(n;f_i)
		-
		\sum_{m=1}^{n}
		\log p_\delta^+(g_j,X_i(m))R_i(m)
		-
		N_i^R(n)I_i(f_i)
		<
		-\frac{\epsilon}{2}n,
		\,
		\pi_i^R(n)\ge\rho
		\right)\\
		&=
		\Pro_{\boldsymbol f}
		\left(
		\sum_{m=1}^{n}
		\Bigl(
		\log f_i(X_i(m))
		-
		\log p_\delta^+(g_j,X_i(m))
		-
		I_i(f_i)
		\Bigr)
		R_i(m)
		<
		-\frac{\epsilon}{2}n,
		\,
		\pi_i^R(n)\ge\rho
		\right)\\
		&\le
		\Pro_{\boldsymbol f}
		\left(
		\sum_{m=1}^{n}
		\Bigl(
		\log f_i(X_i(m))
		-
		\log p_\delta^+(g_j,X_i(m))
		-
		I_i(f_i)
		+
		\frac{\epsilon}{2}
		\Bigr)
		R_i(m)
		<
		0,
		\,
		\pi_i^R(n)\ge\rho
		\right),
	\end{aligned}
\end{equation*}
where, in the last inequality, we used the fact that $N_i^R(n) \leq n$. By definition
\[
I_i(f_i)-\frac{\epsilon}{2}
<
D(f_i,g_j),
\]
and thus for $\delta = I_i(f_i)-\frac{\epsilon}{2}$ the claim follows by Lemma \ref{aux1}.\\
\end{IEEEproof}

We next prove Proposition~\ref{tpro}, which will be used
in the proof of Theorem~\ref{thm_allr}.

\begin{IEEEproof}[Proof of Proposition \ref{tpro}]
	By definition \eqref{tta} of $\tau^R_i(\epsilon)$, we have
	\[
	\{\tau^R_i(\epsilon) > c n\}
	\subseteq
	\{ \exists\, m \ge c n : \hat f_{i,m} \notin V_\epsilon(f_i)\}
	\]
	and by the union bound,
	\[
	\Pro_{\boldsymbol f}(\tau^R_i(\epsilon) > c n)
	\le
	\sum_{m=\lfloor c n \rfloor}^{\infty}
	\Pro_{\boldsymbol f}\left(\hat f_{i,m} \notin V_\epsilon(f_i)\right).
	\]
	Therefore,
	\[
	\sum_{n=1}^\infty \Pro_{\boldsymbol f}(\tau^R_i(\epsilon) > c n)
	\le
	\sum_{n=1}^\infty
	\sum_{m=\lfloor c n \rfloor}^{\infty}
	\Pro_{\boldsymbol f}\big(\hat f_{i,m} \notin V_\epsilon(f_i)\big).
	\]
	For fixed $m$, the number of integers $n$ such that $\lfloor c n \rfloor = m$ is at most $\lceil 1/c \rceil$. Hence,
	\begin{equation*}
		\begin{aligned}
			\sum_{n=1}^\infty \Pro_{\boldsymbol f}(\tau^R_i(\epsilon) > c n)
			&\le
			\sum_{n=1}^\infty
			\sum_{m=n}^{\infty} \lceil 1/c \rceil
			\Pro_{\boldsymbol f}\big(\hat f_{i,m} \notin V_\epsilon(f_i)\big)\\
			&\le \lceil 1/c \rceil \sum_{n=1}^\infty n \Pro_{\boldsymbol f}\big(\hat f_{i,n} \notin V_\epsilon(f_i)\big).
		\end{aligned}
	\end{equation*}
	We now split
	\begin{equation*}
		\begin{aligned}
			\sum_{n=1}^\infty n \Pro_{\boldsymbol f}\big(\hat f_{i,n} \notin V_\epsilon(f_i)\big) \leq &\sum_{n=1}^\infty n \Pro_{\boldsymbol f}\left(\hat f_{i,n} \notin V_\epsilon(f_i),\ \pi_i^R(n) \ge \rho\right)\\
			&+\sum_{n=1}^{\infty} n\Pro_{\boldsymbol f}(\pi_i^R(n) < \rho).
		\end{aligned}
	\end{equation*}
	The first term is summable by condition \eqref{stt3} of Theorem \ref{theo:prob_ao}. The second term is summable by Lemma \ref{exp_d}, because by \eqref{c_star} and the fact that $x_{A}(\boldsymbol f), y_{A}(\boldsymbol f) > 0$ it follows that $c_i^*(\boldsymbol g)>0$ for every $i\in[M]$ and $\boldsymbol g\in\boldsymbol\Phi$, and therefore condition \eqref{stt1} implies \eqref{inf_cond}.
\end{IEEEproof}

\begin{IEEEproof}[Proof of Theorem \ref{thm_allr}]
	We prove the result for $i \in A$, as the case $i \in A^c$ follows analogously. Fix $i \in A$, $\rho \in (0,1]$, and $\epsilon > 0$. We begin with the decomposition
	\begin{align*}
		&\Pro_{\boldsymbol f}\left(
		\frac{1}{n}\Lambda_i^R(n)
		<
		\rho I_i(f_i) - \epsilon,\,
		\pi_i^R(n) \geq \rho
		\right)\\
		&\leq
		\Pro_{\boldsymbol f}\left(
		\frac{1}{n}\Lambda_i^R(n)
		<
		\rho I_i(f_i) - \epsilon,\,
		l_i^1(n)-l_i^0(n) > 0,\,
		\hat l_i(n)-l_i^0(n) > 0,\,
		\pi_i^R(n) \geq \rho
		\right)\\
		&\quad+
		\Pro_{\boldsymbol f}\left(
		\hat l_i(n)-l_i^0(n) \leq 0,\,
		\pi_i^R(n) \geq \rho
		\right)
		+
		\Pro_{\boldsymbol f}\left(
		l_i^1(n)-l_i^0(n) \leq 0,\,
		\pi_i^R(n) \geq \rho
		\right).
	\end{align*}
	
	On the event $\{l_i^1(n)-l_i^0(n)>0,\ \hat l_i(n)-l_i^0(n)>0\}$, we have $\Lambda_i^R(n)=\hat l_i(n)-l_i^0(n)$ by \eqref{allr}. Hence,
	\begin{align*}
		\Pro_{\boldsymbol f}\left(
		\frac{1}{n}\Lambda_i^R(n)
		<
		\rho I_i(f_i) - \epsilon,\,
		\pi_i^R(n) \geq \rho
		\right) &\leq
		2\,\Pro_{\boldsymbol f}\left(
		\hat l_i(n)-l_i^0(n) <
		(\rho I_i(f_i)-\epsilon)n,\,
		\pi_i^R(n) \geq \rho
		\right)\\
		&\quad+
		\Pro_{\boldsymbol f}\left(
		l_i^1(n)-l_i^0(n) \leq 0,\,
		\pi_i^R(n) \geq \rho
		\right).
	\end{align*}
	
	By Theorem~\ref{tGLLR}, the last term is exponentially decaying. It remains to show that
	\begin{equation*}
		\sum_{n=1}^{\infty}
		\Pro_{\boldsymbol{f}}\left(
		\hat{l}_{i}(n)-l^{0}_{i}(n)< (\rho I_i(f_{i})-\epsilon) n,\,
		\pi^{R}_{i}(n) \geq \rho
		\right) < \infty.
	\end{equation*}
	
	Since
	\begin{equation*}
		\begin{aligned}
			\hat{l}_{i}(n)-l^{0}_{i}(n) < (\rho I_i(f_{i})-\epsilon) n
			\;\Leftrightarrow\;
			\hat{l}_{i}(n) -l_i(n;f_i)
			+ l_i(n;f_i) -l^{0}_{i}(n)
			< (\rho I_i(f_{i})-\tfrac{\epsilon}{2}) n - \tfrac{\epsilon}{2} n,
		\end{aligned}
	\end{equation*}
	we obtain
	\begin{equation*}
		\begin{aligned}
			\Pro_{\boldsymbol{f}}\left(
			\hat{l}_{i}(n)-l^{0}_{i}(n)< (\rho I_i(f_{i})-\epsilon) n,\,
			\pi^{R}_{i}(n) \geq \rho
			\right) &\leq
			\Pro_{\boldsymbol{f}}\left(
			\hat{l}_{i}(n)-l_{i}(n;f_i)< -\tfrac{\epsilon}{2} n,\,
			\pi^{R}_{i}(n) \geq \rho
			\right)\\
			&\;\;+
			\Pro_{\boldsymbol{f}}\left(
			l_{i}(n;f_i)-l^{0}_{i}(n)< (\rho I_i(f_{i})-\tfrac{\epsilon}{2}) n,\,
			\pi^{R}_{i}(n) \geq \rho
			\right).
		\end{aligned}
	\end{equation*}
	
	In the proof of Theorem~\ref{tGLLR}, it was shown that the second term is exponentially decaying. Hence, it suffices to show that 
	\begin{equation}\label{tsh1}
		\sum_{n=1}^{\infty}
		\Pro_{\boldsymbol{f}}\left(
		\hat{l}_{i}(n)-l_{i}(n;f_i)< -\tfrac{\epsilon}{2} n,\,
		\pi^{R}_{i}(n) \geq \rho
		\right) < \infty.
	\end{equation}
	
	We decompose
	\[
	l_i(n;f_i)-\hat l_i(n)
	=
	\bar{\Lambda}_i^R(n) + \bar{Q}_i^R(n),
	\]
	where
	\begin{align*}
		\bar{\Lambda}_i^R(n)
		&:= \sum_{m=1}^n
		\left(
		\log\frac{f_i(X_i(m))}{\hat f_{i,m-1}(X_i(m))}
		- D(f_i,\hat f_{i,m-1})
		\right) R_i(m),\\
		\bar{Q}_i^R(n)
		&:= \sum_{m=1}^n D(f_i,\hat f_{i,m-1}) R_i(m).
	\end{align*}
	
	Thus,
	\begin{equation*}
	  \begin{aligned}
		\Pro_{\boldsymbol f}\left(
		l_i(n;f_i)-\hat l_i(n) \geq \tfrac{\epsilon}{2} n,\,
		\pi_i^R(n) \geq \rho
		\right) &\leq \Pro_{\boldsymbol f}\left(
		l_i(n;f_i)-\hat l_i(n) \geq \tfrac{\epsilon}{2} n \right)\\
		&\leq\Pro_{\boldsymbol f}\left(
		\bar{\Lambda}_i^R(n) \geq \tfrac{\epsilon}{4} n
		\right)
		+
		\Pro_{\boldsymbol f}\left(
		\bar{Q}_i^R(n) \geq \tfrac{\epsilon}{4} n\right).
	  \end{aligned}	
	\end{equation*}
	
	Since $\hat{f}_{i,n-1}, R_i(n) \in \cF^{R}_{n-1}$, $\{\bar{\Lambda}_i^R(n) \,:\, n \in \bN \}$ is a $\cF^{R}_n$-martingale. By regularity condition (e) and Rosenthal’s inequality  \cite[Theorem 2.12]{hall14}, there exists $C_0>0$ such that
	\[
	\Exp_{\boldsymbol f}[|\bar{\Lambda}_i^R(n)|^{\mathfrak q}] \le C_0 n^{\mathfrak q/2}.
	\]
	By Markov's inequality
	\[
	\Pro_{\boldsymbol f}\left(
	\bar{\Lambda}_i^R(n) \geq \tfrac{\epsilon}{4} n
	\right)
	\le
	\frac{C_0}{(\epsilon/4)^{\mathfrak q}} \frac{1}{n^{\mathfrak q/2}},
	\]
	and since $\mathfrak q >2$
	\[
	\sum_{n=1}^{\infty}
	\Pro_{\boldsymbol f}\left(
	\bar{\Lambda}_i^R(n) \geq \tfrac{\epsilon}{4} n
	\right)
	< \infty.
	\]
	For the second term, by regularity condition (d) there exists $\delta>0$ such that for any $h \in \Phi_i$ with $\|h-f_i\|_{1}<\delta$, it holds that $D(f_i,h)<\epsilon/8$. By definition of $\tau^R_i(\delta)$ in \eqref{tta}, we note that $D(f_i,\hat f_{i,m-1})<\epsilon/8$ for all $m \geq \tau^R_i(\delta)$. Therefore,
	\[
	\bar{Q}_i^R(n)
	\le
	\sum_{m=1}^{\tau^R_i(\delta)} D(f_i,\hat f_{i,m-1})
	+
	\frac{\epsilon}{8} n.
	\]
	Hence,
	\[
	\Pro_{\boldsymbol f}\left(
	\bar{Q}_i^R(n) \geq \tfrac{\epsilon}{4} n\right)
	\le
	\Pro_{\boldsymbol f}\left(
	\tau^R_i(\delta) \geq \tfrac{\epsilon}{8B} n \right),
	\]
	where $B:=\sup_{g \in \Phi_i} D(f_i,g)$ and by regularity conditions (a) and (d) it follows that $B<\infty$. The summability follows from \eqref{asmp}.
\end{IEEEproof}

\section{}\label{appA}
This appendix provides the proofs of Theorems~\ref{mle}
and~\ref{kde}, which establish that the MLE and KDE,
respectively, satisfy condition \eqref{stt3} of
Theorem~\ref{theo:prob_ao} under the stated assumptions.

\begin{IEEEproof}[Proof of Theorem \ref{mle}]
	We fix $i\in[M]$ and $\epsilon>0$. Let $V_\epsilon(f_i)$ denote the $\epsilon$-neighborhood of $f_i$ defined in \eqref{vva}. Since $\Phi_i\setminus V_\epsilon(f_i)$ is compact as a closed subset of the compact set $\Phi_i$, it admits a finite $\epsilon'$-cover $\{g_j\}_{j=1}^{L(\epsilon')}$ for any $\epsilon' > 0$. For each \(j\in[L(\epsilon')]\), define
	\begin{equation*}
		S_i^j(n)
		:=
		\sum_{m=1}^{n}
		z_i^j(X_i(m))R_i(m),
	\end{equation*}
	where
	\begin{equation*}
		z_i^j(X_i(m))
		:=
		\log f_i(X_i(m))
		-
		\log p_{\epsilon'}^{+}(g_j,X_i(m)),
		\qquad
		m\in\bN.
	\end{equation*}
	We note that by definition \eqref{llr}, for each $j\in[L(\epsilon')]$,
	\[
	l_i(n;f_i)-l_i(n;g)\ge S_i^j(n), \quad \forall \, g\in V_{\epsilon'_j}(g_j).
	\]
	Consequently, if $S_i^j(n)>0$ for every $j\in[L(\epsilon')]$,
	then every density outside $V_\epsilon(f_i)$ has
	strictly smaller likelihood than $f_i$. Since
	$\hat f_{i,n}$ is a maximizer over $\Phi_i$, this implies
	$\hat f_{i,n}\in V_\epsilon(f_i)$.
	Thus, by the union bound,
	\begin{equation*}
		\Pro_{\boldsymbol f}
		\left(
		\hat f_{i,n}\notin V_\epsilon(f_i),\,
		\pi_i^R(n)\ge\rho
		\right)
		\le
		\sum_{j=1}^{L(\epsilon')}
		\Pro_{\boldsymbol f}
		\left(
		S_i^j(n)\le0,\,
		\pi_i^R(n)\ge\rho
		\right).
	\end{equation*}
	By regularity condition (c), \(D(f_i,g_j)>0\) for every \(j\in[L(\epsilon')]\). Hence, Lemma~\ref{aux1} implies that each probability on the right-hand side is exponentially decaying. Since the sum is finite, the left-hand side is also exponentially decaying.
\end{IEEEproof}

The following lemma allows us to control $\sL^1$ estimation
errors through bounds on $\sL^p$ errors, for $p>1$, and
will be used in the proof of Theorem~\ref{kde}. 
 
\begin{lemma}\label{llp}
	Let $(\bS,\cS,\nu_i)$ be a $\sigma$-finite measure space,
	let $p>1$, and let $f_i$ be a density with respect to $\nu_i$.
	For every $\epsilon>0$, there exists a measurable set
	$B_\epsilon\subseteq\bS$ such that
	\begin{equation*}
		0<\nu_i(B_\epsilon)<\infty
	\end{equation*}
	and, for every density $g$ with respect to $\nu_i$,
	\begin{equation*}
		\|g-f_i\|_1>\epsilon
		\quad\Rightarrow\quad
		\|g-f_i\|_p
		>
		\frac{\epsilon}
		{4\nu_i(B_\epsilon)^{1-1/p}}.
	\end{equation*}
\end{lemma}

\begin{IEEEproof}
	Since $f_i$ is a density and $\nu_i$ is $\sigma$-finite, there
	exists a measurable set $B_\epsilon\subseteq\bS$ such that
	\begin{equation*}
		0<\nu_i(B_\epsilon)<\infty
		\qquad\text{and}\qquad
		\int_{B_\epsilon^c}f_i(x)\,\nu_i(dx)
		<
		\frac{\epsilon}{4}.
	\end{equation*}
	Let $g$ be any density with respect to $\nu_i$. Since both $g$
	and $f_i$ integrate to one,
	\begin{align*}
		\int_{B_\epsilon^c}g(x)\,\nu_i(dx)
		&=
		1-\int_{B_\epsilon}g(x)\,\nu_i(dx)\\
		&=
		\int_{B_\epsilon^c}f_i(x)\,\nu_i(dx)
		+
		\int_{B_\epsilon}
		\bigl(f_i(x)-g(x)\bigr)\,\nu_i(dx)\\
		&\leq
		\int_{B_\epsilon^c}f_i(x)\,\nu_i(dx)
		+
		\int_{B_\epsilon}
		|g(x)-f_i(x)|\,\nu_i(dx).
	\end{align*}
	Moreover, because $g$ and $f_i$ are nonnegative,
	\begin{align*}
		\|g-f_i\|_1
		&=
		\int_{B_\epsilon}|g(x)-f_i(x)|\,\nu_i(dx)
		+
		\int_{B_\epsilon^c}|g(x)-f_i(x)|\,\nu_i(dx)\\
		&\leq
		\int_{B_\epsilon}|g(x)-f_i(x)|\,\nu_i(dx)
		+
		\int_{B_\epsilon^c}g(x)\,\nu_i(dx)
		+
		\int_{B_\epsilon^c}f_i(x)\,\nu_i(dx)\\
		&\leq
		2\int_{B_\epsilon}|g(x)-f_i(x)|\,\nu_i(dx)
		+
		2\int_{B_\epsilon^c}f_i(x)\,\nu_i(dx).
	\end{align*}
	Since $\nu_i(B_\epsilon)<\infty$, by Holder's inequality,
	\begin{equation*}
		\int_{B_\epsilon}|g(x)-f_i(x)|\,\nu_i(dx)
		\leq
		\nu_i(B_\epsilon)^{1-1/p}\|g-f_i\|_p.
	\end{equation*}
	Therefore,
	\begin{equation*}
		\|g-f_i\|_1
		\leq
		2\nu_i(B_\epsilon)^{1-1/p}\|g-f_i\|_p
		+
		\frac{\epsilon}{2}.
	\end{equation*}
	If $\|g-f_i\|_1>\epsilon$, it follows that
	\begin{equation*}
		2\nu_i(B_\epsilon)^{1-1/p}\|g-f_i\|_p
		>
		\frac{\epsilon}{2}.
	\end{equation*}
\end{IEEEproof}

\begin{IEEEproof}[Proof of Theorem \ref{kde}]
	We fix $i \in [M]$, and $\epsilon > 0$. By Lemma \ref{llp}, there is $\epsilon'>0$ that depends only on $\epsilon$, $p$ and $f_i$, such that
	\begin{equation*}
	 \Pro_{\boldsymbol f}\left(\int_{\bR} |f_i(x)-\hat f_{i,n}(x)| dx
	 \ge \epsilon,\,\pi_i^R(n)\ge \rho\right) \leq \Pro_{\boldsymbol f}\left(\int_{\bR} |f_i(x)-\hat f_{i,n}(x)|^p dx
	 \ge \epsilon',\,\pi_i^R(n)\ge \rho\right).
	\end{equation*}
	Thus, in order to prove the claim it suffices to show that
	\begin{equation}\label{cl1}
		\sum_{n=1}^{\infty}
		n
		\Pro_{\boldsymbol f}\left(
		\int_{\bR}
		|f_i(x)-\hat f_{i,n}(x)|^p dx
		\ge \epsilon',\,
		\pi_i^R(n)\ge \rho
		\right)
		<\infty, \quad \forall \, i \in [M].
	\end{equation}
	By union bound
	\begin{equation}\label{cl2}
	 \begin{aligned}	
	 &\Pro_{\boldsymbol f}\left(
	 \int_{\bR}
	 |f_i(x)-\hat f_{i,n}(x)|^p dx
	 \ge \epsilon',\,
	 N_i^R(n)\ge \rho\, n
	 \right) \\
	 &\;\;\leq \sum_{m=\lceil \rho n\rceil}^{n} \Pro_{\boldsymbol f}\left(
	 \int_{\bR}
	 |f_i(x)-\hat f_{i,n}(x)|^p dx
	 \ge \epsilon',\,
	 N_i^R(n)=m
	 \right)
	 \end{aligned}
	\end{equation}
	and by Markov's inequality
	\begin{equation*}
	 \Pro_{\boldsymbol f}\left(
	 \int_{\bR}
	 |f_i(x)-\hat f_{i,n}(x)|^p dx
	 \ge \epsilon',\,N_i^R(n)=m\right) \leq \frac{1}{\epsilon'} \mathbb E_{\boldsymbol f}
	 \left[
	 \int_{\bR}
	 |\hat f_{i,n}(x){-}f_i(x)|^p dx
	 ;
	 N_i^R(n){=}m
	 \right].
	\end{equation*}
	In order to prove the claim, it suffices to show that there exists a constant \(D_0>0\) such that
	\begin{equation}\label{tsh}
		\mathbb E_{\boldsymbol f}
		\left[
		\int_{\bR}
		|\hat f_{i,n}(x)-f_i(x)|^p\,dx
		;
		N_i^R(n)=m
		\right]
		\le
		D_0
		\left(
		\frac{1}{(n h_{m})^{p/2}}
		+
		h_{m}^{2p}
		\right).
	\end{equation}
	According to \cite[Theorem 24.1]{van2000asymptotic} we choose the function $h$ as $h_m:=m^{-1/5}$. Therefore, $\eta_{i,n}:=(N_i^R(n))^{-1/5}$ and on the event \(\{N_i^R(n)=m\}\) we have $\eta_{i,n}=h_m=m^{-1/5}$. Since $\rho n \leq m \leq n$, it holds
	\begin{equation*}
	 n^{-1/5} \leq h_m = m^{-1/5} \leq (\rho n)^{-1/5}.
	\end{equation*}
	and as a result
	\begin{equation*}
		\begin{aligned}
			\left(
			\frac{1}{n^{p/2}h_{m}^{p/2}}
			+
			h_{m}^{2p}
			\right)
			&\le
			\left(
			\frac{1}{n^{2p/5}}
			+
			\frac{1}{\rho^{2p/5}}
			\frac{1}{n^{2p/5}}
			\right)=\left(1 + \frac{1}{\rho^{2p/5}} \right) \frac{1}{n^{2p/5}}.
		\end{aligned}
	\end{equation*}
	Since the upper bound is independent of $m$, by \eqref{cl2} we deduce that
	\begin{equation*}
	  \Pro_{\boldsymbol f}\left(
	  \int_{\bR}
	  |f_i(x)-\hat f_{i,n}(x)|^p dx
	  \ge \epsilon',\,N_i^R(n)\ge \rho\, n \right) \leq \left(1 + \frac{1}{\rho^{2p/5}} \right)(1-\rho) \frac{1}{n^{(2p-5)/5}},
	\end{equation*}
	and thus for $p > 15/2$ we prove claim \eqref{cl1}.
	
	In order to prove \eqref{tsh}, we write
	\[
	\hat f_{i,n}(x)-f_i(x)
	=
	\bigl(
	\hat f_{i,n}(x)-\bar{f}_{n,m}
	\bigr)
	+
	\bigl(
	\bar{f}_{n,m}-f_i(x)
	\bigr),
	\]
	where
	\begin{equation*}
	 \bar{f}_{n,m} := \frac{1}{m h_m} \sum_{k=1}^{n} R_i(k) \Exp_{\boldsymbol f}\left[ K\!\left(
	 \frac{x-X_i(k)}{h_m}
	 \right) \right],
	\end{equation*}
	and by convexity, we have
	\[
	\begin{aligned}
		\mathbb E_{\boldsymbol f}
		\left[
		\left|\hat f_{i,n}(x)-f_i(x)\right|^p;N_i^R(n){=}m
		\right]
		\le &
		2^{p-1}
		\mathbb E_{\boldsymbol f}
		\left[
		\left|
		\hat f_{i,n}(x)-\bar{f}_{n,m}
		\right|^p;N_i^R(n){=}m
		\right] \\
		&+
		2^{p-1}E_{\boldsymbol f}
		\left[
		\left|
		\bar{f}_{n,m}-f_i(x)
		\right|^p;N_i^R(n){=}m \right].
	\end{aligned}
	\]
    By Tonelli's theorem
    \begin{equation*}
      \begin{aligned}
    	\mathbb E_{\boldsymbol f}
    	\left[ \int_{\bR}
    	|\hat f_{i,n}(x){-}f_i(x)|^p dx
    	;N_i^R(n){=}m\right] &= \int_{\bR} \mathbb E_{\boldsymbol f}
    	\left[
    	|\hat f_{i,n}(x){-}f_i(x)|^p
    	;N_i^R(n){=}m\right] dx\\
    	& \leq 2^{p-1} \int_{\bR}
    	\mathbb E_{\boldsymbol f}
    	\left[
    	\left|
    	\hat f_{i,n}(x)-\bar{f}_{n,m}
    	\right|^p;N_i^R(n){=}m
    	\right] dx \\
    	&\quad +
    	2^{p-1} \int_{\bR} E_{\boldsymbol f}
    	\left[
    	\left|
    	\bar{f}_{n,m}-f_i(x)
    	\right|^p;N_i^R(n){=}m \right]  dx,
      \end{aligned}
    \end{equation*}	
	where the first term on the right hand side is referred to as the \textit{stochastic} term, while the second one is the \textit{bias} term. 
	
	We first bound the \textit{bias} term. Since $\{ X_i(n)\, :\, n \in \bN \}$ are iid, we get that on the event $\{N_i^R(n) = m\}$,
	\begin{equation*}
	  \begin{aligned}
	  	E_{\boldsymbol f}
	  	\left[\left|\bar{f}_{n,m}-f_i(x)\right|^p;N_i^R(n){=}m \right] &= E_{\boldsymbol f}
	  	\left[\left|\frac{1}{h_m}
	  	\mathbb E_{\boldsymbol f}\!\left[
	  	K\!\left(
	  	\frac{x-X_i(1)}{h_m}
	  	\right)
	  	\right]-f_i(x)\right|^p;N_i^R(n){=}m \right]\\
	  	&\leq \left|\frac{1}{h_m}
	  	\mathbb E_{\boldsymbol f}\!\left[
	  	K\!\left(
	  	\frac{x-X_i(1)}{h_m}
	  	\right)
	  	\right]-f_i(x)\right|^p.
	  \end{aligned}
	\end{equation*}
	We have
	\begin{equation*}
	 \frac{1}{h_m}
	 \mathbb E_{\boldsymbol f}\!\left[
	 K\!\left(
	 \frac{x-X_i(1)}{h_m}
	 \right)
	 \right] =
	 \frac{1}{h_m}
	 \int_{\bR}
	 K\!\left(
	 \frac{x-y}{h_m}
	 \right)
	 f_i(y)\,dy =
	 \int_{\bR}
	 K(u) f_i(x-h_m u)\,du.
	\end{equation*}
	Since $K$ is a density,
	\begin{equation*}
	  \frac{1}{h_m}
	  \mathbb E_{\boldsymbol f}\!\left[
	  K\!\left(
	  \frac{x-X_i(1)}{h_m}
	  \right)
	  \right] -f_i(x) = \int_{\bR}
	  K(u) (f_i(x-h_m u) - f_i(x)) \,du.
	\end{equation*}
	Using Taylor's expansion with the integral form of the remainder \cite[Theorem 24.1]{van2000asymptotic},
	\[
	f_i(x-h_m u)-f_i(x)
	=
	-h_m u\,f_i'(x)
	+
	h_m^2u^2
	\int_0^1
	f_i''(x-th_m u)(1-t)\,dt,
	\]
	which implies
	\begin{equation*}
	  \begin{aligned}
		\frac{1}{h_m}
		\mathbb E_{\boldsymbol f}\!\left[
		K\!\left(
		\frac{x-X_i(1)}{h_m}
		\right)
		\right] -f_i(x) &= -h_m f_i'(x)
		\int_{\bR} uK(u)\,du +
		h_m^2
		\int_{\bR}
		u^2K(u)
		\int_0^1
		f_i''(x-th_m u)(1-t)\,dt\,du\\
		&= h_m^2
		\int_{\bR}
		u^2K(u)
		\int_0^1
		f_i''(x-th_m u)(1-t)\,dt\,du,
	  \end{aligned}
	\end{equation*}
	where the second equality follows by assumption (1). Considering the $\sL^p$-norm both sides, by Minkowski's inequality \cite[Theorem 202]{hardyinequalities} we have
		\[
	\begin{aligned}
		\left\|
		\frac{1}{h_m}
		\mathbb E_{\boldsymbol f}\!\left[
		K\!\left(
		\frac{\cdot -X_i(1)}{h_m}
		\right)
		\right] -f_i
		\right\|_p
		&\le
		h_m^2
		\int_{\bR}
		u^2 |K(u)|
		\int_0^1
		(1-t)
		\left\|
		f_i''(\cdot-th_m u)
		\right\|_p
		\,dt\,du.
	\end{aligned}
	\]
	By translation invariance of the $\sL^p$ norm,
	\begin{equation*}
		\left\|f_i''(\cdot -th_m u)\right\|_p = \left\|f_i''\right\|_p,
	\end{equation*}
	and thus by assumptions (2) and (b) there exists a constant $D_1 >0$ such that
	\begin{equation}\label{r1}
		\int_{\bR} E_{\boldsymbol f}
		\left[
		\left|
		\bar{f}_{n,m}-f_i(x)
		\right|^p;N_i^R(n){=}m \right]  dx \leq \left\|
		\frac{1}{h_m}
		\mathbb E_{\boldsymbol f}\!\left[
		K\!\left(
		\frac{\cdot -X_i(1)}{h_m}
		\right)
		\right] -f_i
		\right\|_p^p
		\le
		D_1 h_m^{2p}.
	\end{equation}
	Next, we bound the stochastic term. For this, we set
	\begin{equation}\label{yy}
		Y_i(x,k)
		:=
		\frac{1}{h_m}
		K\!\left(
		\frac{x-X_i(k)}{h_m}
		\right)
		-
		\frac{1}{h_m}
		\mathbb E_{\boldsymbol f}\!\left[
		K\!\left(
		\frac{x-X_i(k)}{h_m}
		\right)
		\right],
	\end{equation}
	and
	\[
	S_n:=\sum_{k=1}^n R_i(k)\, Y_i(x,k).
	\]
	Therefore,
	\begin{equation*}
	 \begin{aligned}
	   \mathbb E_{\boldsymbol f}
	   \left[\left| \hat f_{i,n}(x)-\bar{f}_{n,m}\right|^p;N_i^R(n){=}m\right]&=\frac{1}{m^p}\mathbb E_{\boldsymbol f}\left[|S_n|^p;N_i^R(n){=}m\right]\\
	   &\leq \frac{1}{m^p}\mathbb E_{\boldsymbol f}\left[|S_n|^p\right].
	 \end{aligned}
	\end{equation*}
	Since $R_i(n) \in \cF^{R}_{n-1}$, the $\{ S_n \, :\, n \in \bN \}$ is a martingale with respect to $\cF^{R}_{n}$. By Rosenthal's inequality \cite[Theorem 2.12]{hall14}, there exists a constant \(D_2>0\) such that
	\begin{equation}
		\begin{aligned}
			\mathbb E_{\boldsymbol f}[|S_n|^p]
			\le D_2
			\left\{
			n^{p/2}
			\left(
			\mathbb E_{\boldsymbol f}[|Y_i(x,1)|^2]
			\right)^{p/2}
			+
			n\,
			\mathbb E_{\boldsymbol f}[|Y_i(x,1)|^p]
			\right\},
		\end{aligned}
	\end{equation}
	and thus
	\begin{equation}\label{int}
	  \begin{aligned}
	  &\int_{\bR} E_{\boldsymbol f}
	  \left[\left|\bar{f}_{n,m}-f_i(x)\right|^p;N_i^R(n){=}m \right]  dx\\
	  & \quad \le
	  \frac{D_2}{m^p}
	  \left[
	  n^{p/2}
	  \int_{\bR}
	  \left(
	  \mathbb E_{\boldsymbol f}[|Y_i(x,1)|^2]
	  \right)^{p/2}
	  dx
	  +
	  n
	  \int_{\bR}
	  \mathbb E_{\boldsymbol f}[|Y_i(x,1)|^p]
	  dx
	  \right].
	  \end{aligned}
	\end{equation}
	By definition of \(Y_i(x,1)\) in \eqref{yy}, \(\mathbb E_{\boldsymbol f}[|Y_i(x,1)|^2]\) is the variance of
	\[
	\frac{1}{h_m}
	K\!\left(
	\frac{x-X_i(1)}{h_m}
	\right),
	\]
	and thus
		\[
		\mathbb E_{\boldsymbol f}[|Y_i(x,1)|^2]
		\le
		\mathbb E_{\boldsymbol f}\left[
		\frac{1}{h_m^2}
		K^2\!\left(
		\frac{x-X_i(1)}{h_m}
		\right)
		\right]=
		\frac{1}{h_m}
		\int_{\bR}
		K^2(u)\,
		f_i(x-h_m u)\,du.
	    \]
	 Using Minkowski's inequality,
		\[
	\left\|
	\mathbb E_{\boldsymbol f}[|Y_i(x,1)|^2]
	\right\|_{p/2}
	\le
	\frac{1}{h_m}
	\int_{\bR}
	K^2(u)
	\,
	\|f_i(\cdot-h_m u)\|_{p/2}
	\,du,
	\]
	and by translation invariance of the $L^{p/2}$ norm, and assumptions (3) and (a) we deduce that there exists a constant \(D_3>0\) such that
	\begin{equation*}
		\int_{\bR}
		\left(
		\mathbb E_{\boldsymbol f}[|Y_i(x,1)|^2]
		\right)^{p/2}
		dx=\left\|
		\mathbb E_{\boldsymbol f}[|Y_i(x,1)|^2]
		\right\|_{p/2}^{p/2}
		\le
		\frac{D_3}{h_m^{p/2}}.
	\end{equation*}
	For the second integral in \eqref{int}, by convexity we have
	\[
	\begin{aligned}
		\mathbb E_{\boldsymbol f }[|Y_i(x,1)|^p]
		\le
		2^{p-1}
		\mathbb E_{\boldsymbol f }\left[
		\frac{1}{h_m^p}
		\left|
		K\!\left(
		\frac{x-X_i(1)}{h_m}
		\right)
		\right|^p
		\right]+
		2^{p-1}
		\left|
		\mathbb E_{\boldsymbol f }\left[
		\frac{1}{h_m}
		K\!\left(
		\frac{x-X_i(1)}{h_m}
		\right)
		\right]
		\right|^p,
	\end{aligned}
	\]
	and since $p >1$ by Holder's inequality
	\[
	\left|
	\mathbb E_{\boldsymbol f}\left[
	\frac{1}{h_m}
	K\!\left(
	\frac{x-X_i(1)}{h_m}
	\right)
	\right]
	\right|^p
	\le
	\mathbb E_{\boldsymbol f}\left[
	\frac{1}{h_m^p}
	\left|
	K\!\left(
	\frac{x-X_i(1)}{h_m}
	\right)
	\right|^p
	\right],
	\]
	which gives
	\[
	\mathbb E_{\boldsymbol f}[|Y_i(x,1)|^p]
	\le
	2^p
	\mathbb E_{\boldsymbol f}\left[
	\frac{1}{h_m^p}
	\left|
	K\!\left(
	\frac{x-X_i(1)}{h_m}
	\right)
	\right|^p
	\right].
	\]
	Therefore,
	\[
	\begin{aligned}
		\int_{\bR}
		\mathbb E_{\boldsymbol f}[|Y_i(x,1)|^p]
		\,dx
		&\le
		2^p\,
		\mathbb E_{\boldsymbol f}\left[
		\int_{\bR}
		\frac{1}{h_m^p}
		\left|
		K\!\left(
		\frac{x-X_i(1)}{h_m}
		\right)
		\right|^p
		dx
		\right]\\
		&\quad=
		2^p\,
		\int_{\bR}
		\int_{\bR}
		\frac{1}{h_m^p}
		\left|
		K\!\left(
		\frac{x-y}{h_m}
		\right)
		\right|^p
		dx\,f_i(y)\,dy.
	\end{aligned}
	\]
	For fixed \(y\), we do the change of variables $u=(x-y)/h_m$, and
	\[
	\int_{\bR}
	\frac{1}{h_m^p}
	\left|
	K\!\left(
	\frac{x-y}{h_m}
	\right)
	\right|^p
	dx
	=
	\frac{1}{h_m^{p-1}}
	\int_{\bR}
	|K(u)|^p\,du.
	\]
	Hence, by the assumption (3), there exists a constant \(D_4>0\) such that
	\[
	\int_{\bR}
	\mathbb E_{\boldsymbol f }[|Y_i(x,1)|^p]
	\,dx
	\le
	D_4
	\frac{1}{h_m^{p-1}}.
	\]
	As a result, by \eqref{int}
	\begin{equation*}
	 \begin{aligned}
	 \int_{\bR} E_{\boldsymbol f} \left[\left|\bar{f}_{n,m}-f_i(x)\right|^p;N_i^R(n){=}m \right]  dx &\leq \frac{D_2}{m^p}
	 \left[
	 n^{p/2}
	 \frac{D_3}{h_m^{p/2}}
	 +
	 n
	 \frac{D_4}{h_m^{p-1}}
	 \right]\\
	 &\leq \frac{D_2}{\rho^p}
	 \left[
	 \frac{D_3}{n^{p/2}}
	 \frac{1}{h_m^{p/2}}
	 +
	 \frac{D_4}{n^{p-1}}
	 \frac{1}{h_m^{p-1}}
	 \right],
	 \end{aligned}
	\end{equation*}
	where for the last inequality we used the fact that \(\rho n \le m \). Since \(p>2\) and \(1 \leq n h_m \),
	\[
	\frac{1}{(nh_m)^{p-1}}
	\le
	\frac{1}{(nh_m)^{p/2}},
	\]
	and thus there exists \(D_5>0\) such that
	\begin{equation}\label{r2}
	 \int_{\bR} E_{\boldsymbol f} \left[\left|\bar{f}_{n,m}-f_i(x)\right|^p;N_i^R(n){=}m \right]  dx \leq \frac{D_5}{(nh_m)^{p/2}}.
	\end{equation}
	 Combining \eqref{r1} and \eqref{r2}, we show \eqref{tsh}.
\end{IEEEproof}

\end{document}